\documentclass{amsart}
\usepackage[margin=1.2in]{geometry}
\usepackage[english]{babel}
\usepackage{amsmath,amssymb,amsthm,graphicx,xcolor,latexsym}
\usepackage{booktabs}
\usepackage{amsmath}
\usepackage{subcaption}
\usepackage{float}
\usepackage{bbm}
\newcommand{\assign}{:=}

\usepackage{microtype}
\usepackage{xcolor}
\definecolor{myblue}{HTML}{4169E1}
\definecolor{myred}{HTML}{E63946}
\definecolor{mygreen}{HTML}{2ECC71} 
\definecolor{myforestgreen}{RGB}{34, 139, 34}
\newcommand{\Vone}{\mathcal{V}^1([0,T];\RR^m)}
\newcommand{\bmu}{\boldsymbol{\mu}}
\usepackage{tikz-cd}
\usetikzlibrary{cd}
\definecolor{linkblue}{HTML}{1E5EFF}
\usepackage{cite}
\usepackage[
    colorlinks,
    linkcolor=linkblue,
    citecolor=linkblue,
    urlcolor=linkblue,
    bookmarks=false,
    hypertexnames=true
]{hyperref}

\theoremstyle{plain}
\newtheorem{theorem}{Theorem}[section]
\newtheorem{remark}[theorem]{Remark}
\theoremstyle{definition}
\newtheorem{example}[theorem]{Example}
\theoremstyle{plain}
\newtheorem{lemma}[theorem]{Lemma}
\newtheorem{proposition}[theorem]{Proposition}
\newtheorem{Assumption}[theorem]{Assumption}
\newtheorem{corollary}[theorem]{Corollary}
\newtheorem{definition}[theorem]{Definition}
\newtheoremstyle{break}
    {\topsep}
    {\topsep}
    {\itshape}
    {}
    {\bfseries}
    {.}
    {\newline}
    {}
\theoremstyle{break}

\usepackage{bbm}

\newcommand{\cA}{\mathcal{A}}
\newcommand{\cB}{\mathcal{B}}

\newcommand{\cF}{\mathcal{F}}

\newcommand{\cH}{\mathcal{H}}
\newcommand{\cI}{\mathcal{I}}

\newcommand{\cK}{\mathcal{K}}
\newcommand{\cL}{\mathcal{L}}
\newcommand{\cM}{\mathcal{M}}

\newcommand{\cQ}{\mathcal{Q}}

\newcommand{\cS}{\mathcal{S}}
\newcommand{\cT}{\mathcal{T}}

\newcommand{\cV}{\mathcal{V}}
\newcommand{\cW}{\mathcal{W}}

\newcommand{\cY}{\mathcal{Y}}
\newcommand{\cZ}{\mathcal{Z}}

\definecolor{darkred}{rgb}{0.7,0.1,0.1}
\definecolor{darkblue}{rgb}{0.4,0.1,0.8}
\definecolor{darkgreen}{rgb}{0.1,0.7,0.1}
\definecolor{franck}{rgb}{0,0,1}
\definecolor{pagebackground}{rgb}{1,1,1}
\definecolor{pageforeground}{rgb}{0,0,0}

\colorlet{symbols}{blue!90!black}
\colorlet{connection}{red!30!black}
\colorlet{boxcolor}{blue!50!black}

\newcommand{\Lkernel}{L^{\infty,1}(\Delta^2; \mathcal{L}(\RR^d;\RR^m))}
\newcommand{\Ltwokernel}{L^{\infty,2}(\Delta^2; \mathcal{L}(\RR^d;\RR^m))}

\renewcommand{\AA}{\mathbb{A}}

\newcommand{\EE}{\mathbb{E}}
\newcommand{\FF}{\mathbb{F}}

\newcommand{\NN}{\mathbb{N}}

\newcommand{\RR}{\mathbb{R}}

\newcommand{\ZZ}{\cZ}

\newcommand{\bV}{\mathbf{V}}

\newcommand{\bx}{\mathbf{x}}
\newcommand{\bY}{\mathbf{Y}}

\newcommand{\bZ}{\mathbf{Z}}

\definecolor{dg}{rgb}{0, 0.5, 0}
\definecolor{dp}{rgb}{0.50, 0, 0.40}

\newcommand{\vertiii}[1]{{\left\vert\kern-0.25ex\left\vert\kern-0.25ex\left\vert #1 
		\right\vert\kern-0.25ex\right\vert\kern-0.25ex\right\vert}}

\newcommand{\eps}{\varepsilon}
\newcommand{\dd}{\mathop{}\!\mathrm{d}}

\usepackage{shuffle}
\usepackage{mathrsfs}
\usepackage{mathtools}
\usepackage{atbegshi}
\usepackage[nameinlink,capitalize]{cleveref} 

\usepackage[autostyle]{csquotes}
\usepackage[toc,page]{appendix}

\usepackage{graphicx}

\input{tree_code}

\title{Stochastic Taylor expansions for Volterra equations and their geometry}

\author{Christa Cuchiero}
\address{Christa Cuchiero, Department of Statistics and Operations Research, University of Vienna,
Kolingasse 14--16, 1090 Vienna, Austria and Department of Mathematics, University of Oslo, P.O.\ Box 1053 Blindern, 0316 Oslo, Norway}

\author{Paul P. Hager}
\address{Paul P. Hager, Department of Statistics and Operations Research, University of Vienna,
Kolingasse 14--16, 1090 Vienna, Austria}

\author{Julian Pachschwöll}
\address{Julian Pachschwöll, Department of Statistics and Operations Research, University of Vienna,
Kolingasse 14--16, 1090 Vienna, Austria}

\author{Luca Pelizzari}
\address{Luca Pelizzari, Department of Statistics and Operations Research, University of Vienna,
Kolingasse 14--16, 1090 Vienna, Austria}

\date{\today}


\begin{document}

\begin{abstract}
We develop series representations for solutions to general nonlinear controlled and stochastic Volterra equations in terms of two new types of signature, called \emph{geometric and It\^o Volterra signature}. These signatures collect the necessary iterated kernel integrals and satisfy (quasi) shuffle identities, which  allow us to prove universal approximation results for Volterra path functionals.
We then derive finite order expansions with precise remainder terms for both controlled and stochastic Volterra equations with possibly singular kernels, as well as expansions of their expectations. For the latter, we provide sufficient conditions for absolute convergence of the resulting series representation and verify them for polynomial Volterra processes. 
As an application, we propose a weak approximation principle for stochastic Volterra equations based on matching low-order expected signature coordinates of the kernel. Numerical experiments for Volterra--Jacobi processes demonstrate promising performance.

\end{abstract}

\maketitle

\section{Introduction}
Stochastic Taylor expansions for smooth functions of solutions to stochastic differential equations (SDEs) play a fundamental role in numerical analysis, particularly in the development and convergence analysis of higher-order numerical schemes, the derivation of error estimates, and the study of strong and weak approximations, see e.g., \cite{kloedenplaten1992numerical}. 
One of the goals of the current work is to develop similar series expansions for general nonlinear Volterra SDEs, i.e., SDEs driven by  Volterra semimartingales, whose numerical simulation is due to their intrinsic path dependence much more involved than in the standard SDE case with a Markovian driver. 
More precisely, for a Volterra semimartingale of the form
\[
Z_s^{t}=\int_0^s K(t,u) dX_u,
\]
where $X$ is a $d$-dimensional It\^o-semimartingale and $K$ an $m \times d$ matrix-valued, possibly singular or fractional kernel,
we consider  Volterra SDEs of
the form
\begin{equation}\label{eq:volterra_eq}
Y_t = \xi + \sum_{i=1}^m\int_0^t V_i(Y_s)\dd Z_s^{t,i}, \quad 0 \leq t \leq T, \qquad \xi \in \mathbb{R}^k
\end{equation} where  $V_i:\RR^k \rightarrow \RR^k$ are smooth vector-fields. 
Similarly to classical stochastic Taylor expansions, our aim is to   
separate the
nonlinear vector-fields from the  driving signal $Z$ (which in our case involves the kernel)  through
expansions of the form
\begin{equation}\label{eq:taylor_volterra_intro}
   f(Y_t) =f(\xi) + \sum_{n=1}^N\sum_{\mathbf{a}\in {\Lambda}^n} \bV ^{\mathbf{a}}(f)(\xi){\mathbf{Z}}^{\mathbf{a}}_{0,t} + R^{N}(t), \qquad f\in C^{N+1}(\RR^k;\RR).
\end{equation}
Before making all the components in~\eqref{eq:taylor_volterra_intro} precise, let us briefly review  the case $K\equiv \mathrm{Id}$, where Equation~\eqref{eq:volterra_eq} reduces to a standard SDE and where \eqref{eq:taylor_volterra_intro}
coincides with the classical stochastic Taylor-It\^o formulas.

In the simplest case where $X$ is of finite variation, such series representations originate
in the work of Chen~\cite{Chen1957} and were developed for
nonlinear control systems
in~\cite{fliess1981fonctionnelles,sussmann1983lie}.
The standard derivation of these representations follows from repeated applications
of a change-of-variables formula to $f$ and $V$, which for finite variation drivers $X$ leads to 
\begin{equation}\label{eq:classical_expansion}
    f(Y_t)
    =
    f(\xi)
    +\sum_{n=1}^N\sum_{w=i_1\cdots i_n}
    \bV^w(f)(\xi)\,
    \mathbf X_{0,t}^w
    +R^N(t),
\end{equation}
where $\mathbf X$ collects iterated integrals, which is well-known as the path
signature of $X$
\begin{equation}\label{eq:classical_signature}
    \mathbf X_{0,t}^{i_1\cdots i_n}
    =
    \int_{0<r_1<\cdots<r_n<t}
     \dd X_{r_1}^{i_1}\cdots\dd X_{r_n}^{i_n}, \qquad i_1,\dots, i_n \in \{1,\dots,d\}.
\end{equation}
The coefficients are compositions of the differential operators
associated with $V$:
\[
    \bV^i(f)=Df\,V_i,
    \qquad
    \bV^{i_1\cdots i_n}
    =
    \bV^{i_1}\circ\cdots\circ\bV^{i_n}.
\]
Similarly if $X$ is a continuous semimartingale,
 an iterative application of Itô's formula yields an expansion of type \eqref{eq:taylor_volterra_intro}, where
${\mathbf{Z}}^{\mathbf{a}}_{0,t}$
corresponds to iterated It\^o integrals involving the components of $X$ and its quadratic variation; see \cite{Kloeden1991} for the case of Brownian drivers $X$.

The main difficulty in view of  singular kernels $K$ in~\eqref{eq:volterra_eq} arises from the fact that the solution does not necessarily inherit
the regularity of the driver: $Y$ might not be a semimartingale,
and even a smooth driver may not yield a smooth solution.
The classical change-of-variables / It\^o-formula therefore cannot
in general be iterated directly along $t \mapsto Y_t$. Our construction
instead uses a two-parameter representation of the solution 
$$Y_t^{\tau}= \xi+ \sum_{i=1}^m\int_0^tV_i(Y_s^{s})\dd Z^{\tau,i}_s, \qquad 0 \leq t \leq \tau \leq T,$$ with the diagonal  $Y_t^{t}=Y_t$. With a change-of-variables formula in the lower-variable, see Lemma~\ref{lem:ito}, we then achieve the series representation~\eqref{eq:taylor_volterra_intro}. Indeed, iterating this formula produces different upper-parameter evaluations, which correspond to evaluations of the kernel in the first variable. As a consequence, we shall see that compared to the classical case~\eqref{eq:classical_expansion}, we require to enlarge the index set from words, to pairs of words and evaluations maps $\mathbf{a}=(w,\sigma)$. More precisely, in the case where $X$ is a Lipschitz driver, but $K$ general (and potentially singular),  the set $\Lambda^n$ in~\eqref{eq:taylor_volterra_intro} consists of $\mathbf{a}=(w,\sigma)$, where $w=i_1\cdots i_n$ and
\[
    \sigma:\{1,\ldots,n\}\longrightarrow\{2,\ldots,n+1\},
    \qquad \sigma(\ell)>\ell.
\]
A precise definition and algebraic properties of $\Lambda^n$ are presented in Section~\ref{sec:word}. Informally, the word records the components of the driving signal, while
$\sigma(l)$ specifies which vector-field factor, or
the terminal function $f$, is differentiated by means of the change-of-variable formula. This determines the elementary
differential $\bV^{w,\sigma}(f)$, which is introduced and studied in Section~\ref{sec:VF_calculus}. The indices $\mathbf{a}=(w,\sigma)$ record the corresponding iterated Volterra integral, where $\sigma$ tracks the evaluation of the kernels, informally
\begin{equation}\label{eq:generalized_volterra_signature_informal}
    \bZ_{0,t}^{w,\sigma}
    =
    \int_{0<r_1<\cdots<r_n<t}
    \dd Z_{r_1}^{r_{\sigma(1)},i_1}\cdots \dd Z_{r_n}^{r_{\sigma(n)},i_n} ,
    \qquad r_{n+1}:=t.
\end{equation}
We call $\bZ$ the \emph{geometric Volterra signature}, and refer to Section~\ref{sec:smooth_VP}
 for the precise definition in case of Lipschitz drivers $X$. 
 Related objects and variants of Volterra signatures have already appeared in the literature, see e.g.~the expansions in~\cite{feng2023cubature} and the review at the end of this introduction.
Let us just remark that it is an extension of the 
Volterra signature introduced in the pioneering work \cite{HT21} with the difference 
that 
our signature is in contrast to the one in \cite{HT21} indeed \emph{geometric}.
By this we mean that it is a group like element and the classical \emph{shuffle-identity} holds true, see Proposition \ref{prop:shuffle_identity}.  This algebra property together with point separation, then allows to deduce a universal approximation result for functionals of Volterra paths  in terms of linear functions of the geometric Volterra signature, see Corollary \ref{cor:universal}. 

For the case of It\^o-semimartingale drivers, we introduce in Definition
\ref{def:Ito_Volterra signature}
the notion of the \emph{It\^o Volterra signature}, which is analogous to \eqref{eq:generalized_volterra_signature_informal}, however with the integrals understood in It\^o sense and an extended alphabet that involves the quadratic variation terms. In this case we can prove that the so-called \emph{quasi-shuffle-identity} holds true, see Proposition \ref{prop:quasi_shuffle}. This together with a continuity result 
for the It\^o Volterra signature, seen as a map from the space of so-called
$\beta$-Hölder semimartingales
equipped with a 
projective limit topology (see Definition
\ref{def:V_hoelder_space}), then allows to prove an almost sure global universal approximation theorem, see Corollary \ref{cor:SMG_universality}.

Having established these universality results for both the geometric and the It\^o Volterra signature, we then move 
to our main goal, namely 
finite-order
expansions with explicit remainders corresponding to the representation in \eqref{eq:taylor_volterra_intro}. 
In Section~\ref{sec:cde_expansion} 
we work this out 
for controlled differential equations (CDEs), i.e., in the 
case when $X$ is a Lipschitz driver, and thus obtain a deterministic Taylor expansion relying on the geometric Volterra signature. In
Section~\ref{sec:stochastic_expansion} a similar program is carried out for semimartingale drivers, relying on the It\^o Volterra signature. In Section \ref{sec:functional_expansions} 
we further extend this to path functionals using the horizontal and vertical derivatives in the sense of \cite{dupire2019functional, cont2013functional}. Theorem \ref{thm:functional_volterra_expansion} can therefore be seen as 
the It\^o Volterra signature version of the Chen-Fliess approximation for diffusion functionals obtained in \cite{Litterer2011OnAC}.

For time-augmented Brownian motion $X$, which is relevant for standard Volterra SDEs, we apply our expansions to find a series representation of the expectation $\mathbb{E}[f(Y_t)]$. More precisely, we identify a reduced
index sets $\Pi^n \subset \Lambda^n$ in Section~\ref{sec:expected_expansion} and suitable growth, integrability, and kernel-decay
conditions (Assumption~\ref{ass:absolut_convergence_assumption}), for the \emph{absolutely convergent} series representation
\begin{equation}\label{eq:expected_expansion_informal}
    \mathbb E[f(Y_t)]
    =
    f(\xi)
    +\sum_{n\geq1}\sum_{\mathbf a\in\Pi^n}
    \bV^{\mathbf a}(f)(\xi)\,
    \bmu_{0,t}^{\mathbf a},
    \qquad
    \bmu_{0,t}^{\mathbf a}
    :=
    \mathbb E[\bZ_{0,t}^{\mathbf a}],
\end{equation} where the expected signature $\bmu_{0,t}^\mathbf{a}$ is given in closed-form. We verify these conditions for a widely-used class of Volterra SDEs in Section~\ref{sec:fPVP}, namely fractional
polynomial Volterra processes.

Finally, in Section~\ref{sec:application}, we exploit the representation in~\eqref{eq:expected_expansion_informal} to develop a weak kernel approximation methodology for Volterra SDEs. Writing $Y(K)$ for the solution to~\eqref{eq:volterra_eq} corresponding to the kernel $K$ and a fixed driver $X$, we aim to select, from a parametric class of simpler kernels $\{K^\theta:\theta\in\Theta\}$, a parameter $\theta$ that minimizes the weak error
$$
\big|\mathbb{E}[f(Y_T(K))]-\mathbb{E}[f(Y_T(K^\theta))]\big|,
\qquad f\in\cF,
$$
for a suitable class of smooth test functions $\cF$. In the spirit of cubature methods on Wiener space~\cite{Lyons2004}, we propose to do so by matching low-order expected geometric Volterra signature coordinates appearing in~\eqref{eq:expected_expansion_informal}, that is, by minimizing discrepancies $|\bmu_{0,T}^{\mathbf a}(K)-\bmu_{0,T}^{\mathbf a}(K^\theta)|$. For large classes of kernels of interest, these quantities are available in closed form, making the resulting optimization problem particularly tractable. In Section~\ref{sec:numerics_jacobi}, we present first numerical results supporting this procedure for sum-of-exponential approximations of fractional Wright--Fisher (Volterra--Jacobi) processes. The methodology shows highly promising approximation properties, which are further exploited for weak Markovian approximations in concrete applications in the companion paper~\cite{cuchiero2026rough}.

Summarizing, our main contributions are the following:

\begin{itemize}
\item Introduction of the \emph{geometric and It\^o Volterra signature}, for which we show  (quasi-) shuffle-identities and universal approximation properties (see Section 
\ref{sec:generalized_VSig});
\item Deterministic and stochastic Taylor expansions for general Volterra CDEs and SDEs with possibly singular kernels, as represented in \eqref{eq:taylor_volterra_intro}, where the function can depend on the current state, but also on the whole path (see Section \ref{sec:cde_expansion} - \ref{sec:functional_expansions});
\item Absolute convergence analysis for expansions of expected values $\mathbb{E}[f(Y_t)]$ under suitable  assumptions, verified for polynomial Volterra processes (see Section \ref{sec:expected_expansion});
    
\item Development of a weak kernel approximation methodology by matching the expected signature expansions, which we numerically illustrate for fractional Jacobi/Wright-Fisher processes.

\end{itemize}

\subsection*{Related works and literature review}

We start by reviewing the general literature on Volterra SDEs. Early well-posedness studies of Volterra SDEs with sufficiently regular kernels and Brownian drivers
go back to~\cite{berger1980volterra_i,berger1980volterra_ii}, and were later extended to c\`adl\`ag semimartingale
drivers in~\cite{Pr85}. First results for singular kernels were obtained
in~\cite{cochran1995stochastic,coutin2001stochastic}, and we refer to~\cite{abi2021weak,promel2023existence,promel2025pathwise,hamaguchi2025weak} for recent generalizations.
Volterra SDEs play a prominent role in rough and path-dependent
volatility modeling~\cite{BayerFrizGatheral2016,GJR18,guyon2023volatility}, which motivated
the analysis of specific model classes,
including affine and polynomial ones~\cite{jaber2019affinevolterraprocesses,
ackermann2022inhomogeneous,bondi2024affine,abi2024polynomial},
neural Volterra SDEs~\cite{bergerhausen2025}, as well as coefficient specifications of
mean-field~\cite{jie2024one,promel2026mean} and path-dependent type~\cite{viens2019martingale}.
Markovian representations of Volterra processes were studied
in~\cite{carmona1998fractional,carmona2000approximation} and later
systematically developed for Volterra SDEs in~\cite{harms2019affine,abi2019markovian,cuchiero2019markovian,
cuchiero2020generalized,hamaguchi2024markovian,gasteratos2025kolmogorov}.
These representations are typically infinite-dimensional, and finite-dimensional approximations by a finite sum of exponentials have become a standard numerical technique in rough
volatility~\cite{jaber2018multifactorapproximationroughvolatility,abi2019lifting,abi2021weak_1,MR4271188,
alfonsi2024approximation,zhu2021markovian,bayer2022markovianapproximationsstochasticvolterra,
bayer2023weak,bossy2025weak,shinozaki2026higher}. Our kernel matching methodology  contributes to this field by introducing a kernel selection procedure based on the weak-error representation derived from our Taylor expansion.

As noted above, expansions of controlled differential equations have
their roots in Chen's work~\cite{Chen1957}, and were developed for control
systems in~\cite{fliess1981fonctionnelles,sussmann1983lie}; for a historical
survey see~\cite[Section~1.2]{beauchard2023expansions}. In parallel,
Taylor expansions were developed for stochastic differential equations
in both Itô and Stratonovich
form~\cite{kunita1980representation,platenwagner1982taylor,azencott1982taylor,arous1989flots,Kloeden1991,castell1993asymptotic}.
Lyons' rough path theory~\cite{Lyons1998} later provided a pathwise
framework for expanding solutions of differential equations driven by
irregular signals, see also~\cite{frizvictoir2008euler} and developments for fractional Brownian
motion in~\cite{baudoin2007operators,neuenkirch2009trees,baudoinzhang2012taylor}.
Further generalizations include Taylor expansions for
SPDEs~\cite{bayer2008cubature,Jentzen2010} and for path-dependent
functionals~\cite{Litterer2011OnAC,dupiretissot2022functional,cuchiero2025functional}. For our expansion in
Itô-Volterra integrals, we build on the algebra of
iterated Itô integrals~\cite{gaines1994algebra}, following the
presentation in~\cite{ebrahimi2015exponentialcontinuous,ebrahimi2015flows,ebrahimi2024exponential}. Our expansions in Theorems~\ref{thm:main_cde_expansion},
\ref{thm:main_sde_expansion}, and~\ref{thm:functional_volterra_expansion}
extend the classical Taylor expansions for controlled differential
equations, stochastic differential equations, and path-dependent
functionals to the Volterra framework. For the identity kernel
$K\equiv\mathrm{Id}$, we recover the corresponding formulas from the
literature cited above. On the level of expectations, stochastic Taylor expansions play a fundamental role in the weak-error analysis, most notably in cubature methods~\cite{Lyons2004,bayer2008cubature,bayer2013cubature,litterer2012high,cass2014error,passeggeri2020signature,herschell2025new,ferrucci2026high} but also other weak-schemes \cite{kloedenplaten1992numerical,kusuoka2001approximation,milstein1997numerical,rossler2006rooted,rossler2006runge}. Comparable growth and integrability conditions to ensure absolute convergence of~\eqref{eq:expected_expansion_informal} (see Theorem~\ref{thm:expansion_expected}) appear in literature for classical SDEs: most notably already in~\cite[Corollary  1 and Proposition 3]{arous1989flots} or \cite[Theorem 2 and Remark 2]{wagner1982taylor}, but see also~\cite{neuenkirch2009trees,belomestny2009holomorphic,baudoinzhang2012taylor,rossler2004stochastic,feng2020taylor}.

The geometric and It\^o Volterra signatures
appearing in our expansions have precursors
in Volterra rough path theory. Building on the pathwise analysis of
Volterra equations in~\cite{deya2009rough,deya2011rough},
the article~\cite{HT21} developed a rough path approach to
well-posedness for equations with possibly singular kernels and
introduced the Volterra signature, that is 
\eqref{eq:generalized_volterra_signature_informal} with $\sigma(\ell)=\ell+1$. Recently this signature was studied from a learning and
computational point of view in~\cite{HHPT25,ii_part}. The latter Volterra signature suffices to represent solutions of
linear controlled Volterra equations~\cite[Proposition~2.38]{HHPT25} and linear Volterra SDEs ~\cite[Corollary~3.24]{hamaguchi2023variation}, but not for the general nonlinear case where the ``full'' geometric Volterra signature 
with the additional evaluation maps
$\sigma$ in~\eqref{eq:generalized_volterra_signature_informal}
is needed; 
see Remark~\ref{rem:classical_VSIG}.
Note that linear Volterra SDEs with
sufficiently regular kernels can also be expanded in terms of the classical signature of time-augmented
Brownian motion as shown in~\cite[Theorem~4.2]{jaber2024pathdependentprocessessignatures}.
  For nonlinear Volterra integrals equations ($X_t=t$ in our notation), the work~\cite{brunner1982runge} derives Taylor expansions indexed by \emph{Volterra-trees} as part of their Runge-Kutta theory. In the rough setting, higher-order expansions indexed by rooted trees were then introduced in~\cite{HTW23,HTW24}, extending the Volterra rough path theory to more irregular drivers. Later~\cite{BK23}  generalized this to
an abstract theory of branched Volterra rough paths in the spirit of
Gubinelli's branched rough paths~\cite{gubinelli2010ramification}. For smooth drivers $X$,
Remark~\ref{rem:link_branched} illustrates the relation between this
branched and our geometric constructions of iterated Volterra integrals. Another related object is the exponentially fading memory (EFM) signature introduced in~\cite{jaber2025exponentially} and further developed in~\cite{nlend2026time,bloch2026exponentially}, which, at least on finite intervals and Lipschitz drivers, is a special case of~\eqref{eq:generalized_volterra_signature_informal} with $\sigma(\ell)\equiv n+1$ and exponential kernels $K$.

Finally, the work most closely related to ours is that of Feng and
Zhang~\cite{feng2023cubature} on cubature methods for stochastic Volterra
integral equations. Building on the functional Itô calculus for
Volterra processes developed in~\cite{viens2019martingale},
they derive a functional stochastic Taylor expansion for Brownian drivers
and regular kernels~\cite[Section~2.2]{feng2023cubature}, which is then exploited to derive cubature formulas. In particular, a Stratonovich version of the signature~\eqref{eq:generalized_volterra_signature_informal}, with a similar kernel evaluation formalism, already appeared in their expansion, see also~\cite[Remark 3.1]{feng2023cubature}. Our signature may be viewed as the Itô counterpart, which has the advantage that it can also be used for singular kernels, in contrast to the Stratonovich
signature in~\cite{feng2023cubature}. Moreover, while Feng and Zhang use the
functional Itô formula in their derivation, we derive the expansion of
$f(Y_t)$ through a finite-dimensional Itô formula applied
the semimartingales $Y_\cdot^t$, see Lemma~\ref{lem:ito}.

\subsection{Notation}\label{sec:notation} 

For any $0\leq s \leq t \leq T$ and $n \in \NN$ we denote by $\Delta^{n}_{s,t}$ the $n$-simplex $\Delta^{n}_{s,t}:= \{(r_1,\dots,r_n)\in [s,t]^n:r_1<\cdots <r_n \}$. 
If the context is clear we set $\Delta^n := \Delta_{0,T}^n$.
For any normed vector space $(E, |\cdot |)$ and function $f:\Delta^2\rightarrow E$ and $p\geq 1$ we define the $L^{\infty,p}$ semi-norm \begin{equation}\label{eq:kernel_norm}
    \Vert f \Vert_{L^{\infty,p}}:= \sup_{0 \leq \tau \leq T} \left (\int_0^\tau |f(\tau,t)|^p\dd t\right)^{1/p}.
\end{equation}
In particular, we denote by $L^{\infty,p}(\Delta^2,E)$ all functions with finite $L^{\infty,p}$ norm. Moreover, for any $k,p\in \mathbb{N}$ we denote by $\mathrm{Pol}_p(\mathbb{R}^k)$ the space of polynomials $P:\RR^k\rightarrow \RR$ of degree at most $p$. For any two-parameter stochastic process $(A_t^\tau)_{(t,\tau)\in \Delta^2}$, we use the following convention to extend to $(t,\tau)\in [0,T]^2$
\begin{equation}\label{eq:convention}
    A_t^\tau:=A_{t\wedge\tau}^\tau,
    \qquad t,\tau\in[0,T].
\end{equation}
Given a probability space $(\Omega,\cF,\mathbb{P})$, we define for $1\leq p <\infty$ \[
L^p=L^p(\Omega,\cF,\mathbb{P}):=\Big  \{f:\Omega \rightarrow \mathbb{R}: f \text{ measurable and } \Vert f \Vert_p:= \Big(\int_\Omega |f|^p \dd \mathbb{P}\Big)^{1/p}<\infty\Big \}.
\]
Moreover, we define the intersection of $L^p$ spaces endowed with the projective limit topology \begin{equation}\label{eq:projective_Lp}
L^{\infty-}:= \bigcap_{p\geq 1}L^p, \qquad d(f,g):= \sum_{q\geq 1}2^{-q} \frac{\Vert f-g \Vert_{q}}{1+\Vert f-g \Vert_{q}}.
\end{equation}
Lastly, we occasionally write $[n] \assign \{1,\dots,n\}$ for $n \in \NN$.

\section{Geometric and It\^o Volterra signatures}\label{sec:prelim}
In this section, we set the stage for the series representations~\eqref{eq:taylor_volterra_intro}. Following the order in which the different objects appear in this representation, we first introduce the algebraic component, namely the index sets $\Lambda^n$ and the associated algebraic structure. We then turn to the analytic component, given by the elementary vector-field operators $\bV^{\mathbf a}(\cdot)$. Finally, we introduce the iterated-integral components $\bZ^\mathbf{a}$, the geometric Volterra signature, which encodes the analytic and stochastic structure of the underlying Volterra driver. Most of the technical proofs are postponed to Appendix~\ref{sec:proofs}.

\subsection{Words, maps and shuffles}\label{sec:word}

We begin with the notion of words over a given alphabet.

\begin{definition}\label{def:words} 
Consider an alphabet $\cA= \{1,\dots,d\}$ with $d\in \mathbb{N}$ letters, and denote by $\cW_\cA^n$ the set of words of length $n\in \NN$ over the alphabet $\cA$, that is \begin{equation}\label{eq:words}
    \cW_\cA^n:=\left \{w=i_1\cdots i_n: i_1,\dots,i_n \in \cA\right \}.
\end{equation} For \(n=0\) we write \(\cW_\cA^0=\{\varnothing\}\), where \(\varnothing\) denotes the empty word. Finally, we denote by \(\cW_\cA\) the set of all words of any length, and whenever the alphabet is clear from the context we write \(\cW=\cW_\cA\).
\end{definition}

While ordinary words suffice to index classical iterated integrals, our geometric Volterra iterated integrals carry additional information from  kernel evaluations. We encode this information through the following admissible evaluation maps.
\begin{definition}\label{def:admissible-maps}
For any $n \in \NN$, we denote by $\Sigma^n$ the set
$$
\Sigma^n
:=
\left\{
\sigma:\{1,\dots,n\}\to\{2,\dots,n+1\},
\ \sigma(k)>k \text{ for all } k=1,\dots,n
\right\},
$$
and set
$
\Sigma^0=\{\varnothing\}.
$ We call $\Sigma^n$ the set of admissible evaluation maps of length $n$, and set
$
\Sigma:=\bigcup_{n\geq 0}\Sigma^n.
$
\end{definition}
Combining the word and evaluation-map components, we arrive at the basic index set used throughout the paper. 

\begin{definition}\label{def:volterra-words}
For any $n\in \NN_0$ and alphabet $\cA$, we define the sets
$$
\Lambda_\cA^n:=\mathcal W_\cA^n\times \Sigma^n, \qquad \Lambda_\cA:=\bigcup_{n\geq 0}\Lambda_\cA^n.
$$ We also introduce the concatenation product $(w,\sigma)\cdot(v,\tau):=(wv,\sigma\tau)$, where $wv$ is the usual concatenation of words, and for $\sigma \in \Sigma^n$, $\tau \in \Sigma^m$ we define \begin{equation}\label{eq:concat_vwords}
    (\sigma \tau)(e)= \begin{cases}
        \sigma(e) & 1 \leq e \leq n \text{ and } \sigma(e)\leq n \\

        n+\tau(e-n) & n+1\leq e \leq n+m \text{ and } \tau(e-n) \leq m \\

        n+m+1, & \text{otherwise}.
    \end{cases}
\end{equation} Finally, we introduce the linear spaces
$$
\mathcal H_\cA^{n}:=\mathrm{span}_{\mathbb R}(\Lambda_\cA^n),
\qquad
\mathcal H_\cA:=\mathrm{span}_{\mathbb R}(\Lambda_\cA).
$$
Whenever the alphabet is clear from the context we write $\cH^n = \cH_\cA^n$ and $\cH = \cH_\cA$ respectively.
\end{definition}
Next, we define a shuffle-product on the vector space $\cH$, 
which leads to the shuffle identity of the Volterra signature. While shuffles on words are typically introduced via the recursion
$$
(wi)\shuffle (vj):= \big(w\shuffle (vj)\big)i+\big((wi)\shuffle v\big)j,\qquad i,j\in\cA,\; w,v\in\cW_\cA,
$$ 
here we choose the following direct and equivalent definition.

\begin{definition}\label{def:classical_shuffle}
     Define the set \begin{equation*}
    \mathrm{Sh}(n,m):= \{\eta \in \cS_{n+m}: \eta(1)< \cdots <\eta (n), \, \eta(n+1)<\cdots < \eta(n+m)\},
\end{equation*} where $\cS_n$ denotes the permutation group on $\{1,\dots,n\}$. Then, for two words $w=i_1\cdots i_n$ and $v=j_1\cdots j_m$, we define 
\begin{equation*}
    w\shuffle v 
    \assign 
    \sum_{\eta \in \mathrm{Sh}(n,m)} w \shuffle_\eta v
    \assign
    \sum_{\eta \in \mathrm{Sh}(n,m)} k_{\eta^{-1}(1)}\cdots k_{\eta^{-1}(n+m)}, 
    \qquad k_1\cdots k_{n+m}=wv.
\end{equation*} Finally, we  extend $\shuffle$ bi-linearly to $\mathrm{span}_\RR(\cW_\cA)$.
\end{definition}
In order to extend the definition to $\cH$, we also introduce a ``shuffle-product'' $\widehat{\shuffle}$ on maps $\Sigma$.

\begin{definition}\label{def:Sigma_shuffle}
    For any $\sigma\in \Sigma^n$, $\tau \in \Sigma^m$ and $\eta \in \mathrm{Sh}(n,m)$, we define \begin{equation}
    (\sigma \widehat{\shuffle}_\eta \tau)(e)
    :=
    \bigl(\eta\circ(\sigma\tau)\circ\eta^{-1}\bigr)(e),
    \qquad
    \eta(n+m+1):=n+m+1.
\end{equation} where we recall the concatenation $\sigma\tau$ defined in~\eqref{eq:concat_vwords}. The shuffle-product between $\sigma$ and $\tau$ is then defined as \begin{equation}\label{eq:shuffle_maps}
        \sigma \widehat{\shuffle} \tau:= \sum_{\eta \in \mathrm{Sh}(n,m)} \sigma \widehat{\shuffle}_\eta \tau \in \mathrm{span}_\RR(\Sigma).
    \end{equation} 
    Similar to the shuffle on words, we extend $\widehat{\shuffle}$ bi-linearly to $\mathrm{span}_\RR(\Sigma)$.
\end{definition}
Finally, we define the shuffle-product on the space $\cH$ as follows.

\begin{definition}
    \label{def:shuffle_vwords}
    We define the shuffle-product on the space  $\Lambda_\cA$ in Definition~\ref{def:volterra-words} via \begin{equation}\label{eq:shuffle_vords}
        (w,\sigma) \shuffle (v,\tau) := \sum_{\eta \in \mathrm{Sh}(n,m)}(w \shuffle_\eta v, \sigma \widehat{\shuffle}_\eta \tau),
    \end{equation} and again extend it bi-linearly to $\cH_\cA$.
\end{definition}

\begin{lemma}\label{lem:algebra}
 $(\cH_\cA,\shuffle)$ is a graded, commutative unital algebra.     
\end{lemma}
In the context of semimartingale drivers, we will work with Itô integration, and thus, unsurprisingly, the shuffle product will no longer be suitable. Nevertheless, after enhancing the semimartingale with its quadratic variation, the role of shuffles is replaced by the \emph{quasi-shuffle} product; we refer to \cite{ebrahimi2015flows} and the references therein. With enhanced paths $\cZ=(Z,[Z])$ in mind, we extend our alphabet $\cA$ to \begin{equation}\label{eq:enhanced_alpha}
    \mathbb{A}= \cA \cup \{[i,j]:i,j \in \cA \}.
\end{equation} 
We can then similarly define words of length $n$ over the alphabet $\AA$ by  $\cW_\AA^n = \{w=i_1\cdots i_n: i_1,\dots,i_n \in \mathbb{A}\}$ and define the set of all words over this alphabet as $\cW_\AA$.
We will make use of the grading $|i|=1$ and $|[j,k]|=2$ for any $i,j,k \in \cA$, and then $|i_1\cdots i_n| = \sum_{l=1}^n|i_l|$.  
Finally, we extend the set $\Lambda_\cA$ introduced in Definition~\ref{def:volterra-words} to \newcommand{\bsigma}{\boldsymbol{\sigma}}
\newcommand{\btau}{\boldsymbol{\tau}} 
\begin{equation}\label{eq:admissible_maps_ito}
    \Lambda_\AA^n:= \left \{(i_1\cdots i_n,\bsigma) \in \cW_\AA^n \times (\Sigma^n)^2: \sigma_1(k)=\sigma_2(k) \text{ when } |i_k|=1 \right \},
\end{equation} where the pairs $\bsigma =(\sigma_1,\sigma_2)$ will give us the required flexibility in the upper-parameters of the quadratic variation processes; see Definition~\ref{def:Ito_Volterra signature} below. Similar to the above, define the set $\Lambda_\AA := \bigcup_{n \geq 0} \Lambda_\AA^n$.

The natural way to define the quasi-shuffle product, is to introduce a notation $i \diamond j = [i,j]$ if $i,j \in \cA$, and $i \diamond j=0$ otherwise, and then recursively\begin{equation}\label{eq:QSH_def}
    wi \star vj = (w \star vj)i+(wi \star v)j + (w \star v)i\diamond j, \qquad i,j \in \mathbb{A}.
\end{equation}
Similar as for the shuffle product, we can equivalently give a direct definition, see, e.g., \cite[Section 4]{ebrahimi2015flows}.

\begin{definition}\label{def:classical_quasi_shuffle}
     Define the set \begin{equation*}
    \mathrm{QSh}(n,m):= \bigcup_{k=\max(n,m)}^{n+m}\Big\{\eta: [n+m]\rightarrow [k]: \eta \text{ surjective,} \begin{array}{c}\eta(1)<\cdots<\eta(n),\\ \eta(n+1)<\cdots<\eta(n+m) \end{array}\Big \}.
\end{equation*} Then, for two words $w=i_1\cdots i_n \in \cW_{\AA}^n$ and $v=j_1\cdots j_m \in \cW_\AA^m$, define
\begin{equation*}
    w\star v 
    \assign 
    \sum_{\eta \in \mathrm{QSh}(n,m)} w \star_\eta v
    \assign 
    \sum_{\eta \in \mathrm{QSh}(n,m)} \Big (\prod_{j \in \eta^{-1}(1)} k_j\Big )\cdots \Big (\prod_{j\in \eta^{-1}(k)} k_j\Big ),
\end{equation*} where we again write $k_1\cdots k_{n+m}=wv$ and $\prod_{j\in \eta^{-1}(e)}k_j= \begin{cases}
    k_j & \eta^{-1}(e)=\{ j \}, \\ k_j \diamond k_l & \eta^{-1}(e) = \{j,l\}
\end{cases}.$
\end{definition}
Intuitively, we see that the choice $k=n+m$ in the union of $\mathrm{QSh}(m,n)$ corresponds to $\mathrm{Sh}(n,m)$, and for all $\max(n,m) \leq k<n+m$ one necessarily has $\eta(e)=\eta(e')$ for some $1\leq e \leq n$ and $n+1 \leq e' \leq n+m$. For Volterra paths with Lipschitz drivers, the cases $k<m+n$ can be ignored since the corresponding integrals vanish, but for Itô integration such domains generate quadratic variation terms, see Examples~\ref{ex:illustration_shuffle} and \ref{ex:illu_example_continued} below.

Similarly as for the shuffle in Definition~\ref{def:shuffle_vwords}, we can define the quasi-shuffle on the space $\Lambda_\AA$ defined in \eqref{eq:admissible_maps_ito}. For this, let us remark that by definition we have $|\eta^{-1}(e)| \in \{1,2\}$ for all $\eta\in \mathrm{QSh}(m,n)$.
\begin{definition}\label{def:Sigma_quasi_shuffle}
For $\bsigma=(\sigma_1,\sigma_2)\in(\Sigma^n)^2$,
$\btau=(\tau_1,\tau_2)\in(\Sigma^m)^2$, and
$\eta\in\mathrm{QSh}(n,m)$ with image $[k]$, set again $\eta(n+m+1):=k+1$. Then we define
\[
    \eta_-(e):=\min\eta^{-1}(e),
    \qquad
    \eta_+(e):=\max\eta^{-1}(e),
    \qquad e\in[k].
\]
Using the concatenation of evaluation maps
from~\eqref{eq:concat_vwords}, we set
\begin{equation}\label{eq:Sigma_shuffle_general}
    \bsigma\widehat{\star}_\eta\btau
    :=
    \left(
        \eta\circ(\sigma_1\tau_1)\circ\eta_-,
        \eta\circ(\sigma_2\tau_2)\circ\eta_+
    \right)
    \in(\Sigma^k)^2.
\end{equation}
For $(w,\bsigma)\in\Lambda_{\AA}^n$ and
$(v,\btau)\in\Lambda_{\AA}^m$, define
\[
    (w,\bsigma)\star(v,\btau)
    :=
    \sum_{\eta\in\mathrm{QSh}(n,m)}
    \left(
        w\star_\eta v,\,
        \bsigma\widehat{\star}_\eta\btau
    \right),
\]
and extend this operation bi-linearly to $\cH_{\AA}$.
\end{definition}

\subsection{Composition of vector fields\label{sec:VF_calculus}}
Having introduced the index sets $\Lambda_\cA$ and $\Lambda_\AA$, we now turn to the second ingredient in the series representation~\eqref{eq:taylor_volterra_intro}, namely the elementary vector-field operators $f\mapsto \bV^{\mathbf a}(f)$. 
As in classical deterministic and stochastic Taylor expansions, these operators arise by repeatedly applying a change-of-variables formula to the function $f$ and the vector fields $V=(V_1,\dots,V_m)$. 
Let us first briefly recall the classical construction, where we follow the exposition in~\cite[Section~3]{ebrahimi2024exponential}. On the enhanced alphabet $\AA$ we set
\begin{equation}\label{eq:classical_VF_level_1}
    \bV^i(f)(x)
    =
    \sum_{l=1}^k
    V_i^l(x)\partial_{x^l}f(x),
    \qquad
    \bV^{[i,j]}(f)(x)
    =
    \frac12
    \sum_{p,q=1}^k
    V_i^p(x)V_j^q(x)
    \partial_{x^px^q}f(x).
\end{equation}
Iterating the change-of-variables formula then naturally leads to the recursive composition
\begin{equation}\label{eq:classical_VF_iteration}
    \bV^{i_1\cdots i_n}(f)(x)
    :=
    \bV^{i_1}
    \bigl(
        \bV^{i_2\cdots i_n}(f)
    \bigr)(x),
    \qquad
    i_1\cdots i_n\in\cW_\AA^n.
\end{equation}
For Volterra equations, the same principle applies, but the resulting operators are no longer functions of a single state variable. Indeed, different occurrences of the solution are evaluated at different upper parameters, and after one expansion step we are naturally led to functions of two variables. At the first level, for the unique admissible evaluation map $\bsigma$ of length one, we set
\begin{equation}\label{eq:vector_fields_first_level}
\bV^{(i,\bsigma)}(f)(x_1,x_2)
=
\begin{cases}
    \displaystyle
    \sum_{l=1}^k
    V_i^l(x_1)\partial_{x_2^l}f(x_2),
    & |i|=1,
    \\[3mm]
    \displaystyle
    \frac12
    \sum_{p,q=1}^k
    V_{i_1}^p(x_1)V_{i_2}^q(x_1)
    \partial_{x_2^px_2^q}f(x_2),
    & i=[i_1,i_2].
\end{cases}
\end{equation}
In particular, $\bV^{(i,\bsigma)}(f)(x,x)=\bV^i(f)(x)$,
so that the classical operator is recovered on the diagonal. The additional role of the admissible evaluation maps becomes visible from the second expansion step onwards. When differentiating
$\bV^{(i,\bsigma)}(f)(x_1,x_2)$, the new derivative may act either on the vector-field factor depending on $x_1$ or on the terminal function depending on $x_2$. More generally, at level $n$, the elementary differential depends on $(x_1,\dots,x_{n+1})$
and each derivative produced at an earlier expansion step has to remember on which of these factors it acts. This bookkeeping is precisely encoded by the admissible map $\bsigma$: each derivative from the $r$-th letter is assigned by $\bsigma$ to the variable on which it acts. This leads to the following definition, and we illustrate its connection to~\eqref{eq:classical_VF_iteration} right after.
\begin{definition}\label{def:vector_fields_iteration}
   For $V_1,\dots,V_m \in C^\infty(\RR^k;\RR^k)$, we introduce the enhanced vector-fields
   \begin{equation}
       \label{eq:generalized_vfields}
       \mathsf{V}_i^{p}(x):= V_{i}^p(x), \qquad \mathsf{V}_{[i,j]}^{p,q}(x) = \frac{1}{2}V_{i}^p(x)V_{j}^q(x), \qquad p,q=1,\dots,k,
   \end{equation}
   Then, for any $f\in C^\infty(\RR^k;\RR)$ and $\mathbf{a} =(i_1\cdots i_n,\bsigma)\in {\Lambda}_\AA^n$, we define
   \begin{equation}\label{eq:v_field_operator_generalized}
       \mathbf{V}^{\mathbf{a}}(f)(\mathbf{x})
       =
       \sum_{l_1 \in [k]^{|i_1|}}\cdots
       \sum_{l_n \in [k]^{|i_n|}}
       \prod_{q=1}^n
       \bigg(
       \prod_{\substack{(r,h): \bsigma_h(r)=q\\ 1\leq h\leq |i_r|}}
       \partial_{l^h_r}
       \bigg)
       \mathsf{V}_{i_q}^{l_q}(x_q)
       \bigg(
       \prod_{\substack{(r,h): \bsigma_h(r)=n+1\\ 1\leq h\leq |i_r|}}
       \partial_{l^h_r}
       \bigg)
       f(x_{n+1}),
   \end{equation}
   where we write $\mathbf{x}=(x_1,\dots,x_{n+1})$, $[k]=\{1,\dots,k\}$, and recall $|i_r|=1$ for single and $|i_r|=2$ for double letters.
\end{definition}
\begin{remark}
    \begin{itemize}
        \item[(i)]
        The map $\bsigma$ records where each derivative acts. More precisely,
        \[
            \#\{(r,h):1\leq h\leq |i_r|,\ \bsigma_h(r)=q\}
        \]
        is the total number of derivatives acting on the $q$-th variable.
        Since $\bsigma_h(r)>1$, no derivative can act on the newly introduced
        factor $\mathsf{V}_{i_1}(x_1)$, which reflects the
        multiplication by the vector field in the
        change-of-variables formula.

        \item[(ii)]
        Restricted to single letters,
        $\mathbf a=(i_1\cdots i_n,\sigma)\in\Lambda_\cA^n$, the definition~\eqref{eq:v_field_operator_generalized}
        reduces to
        \begin{equation}\label{eq:elementary_vectorfield}
            \bV^{\mathbf{a}}(f)(x_1,\dots,x_{n+1})
            :=
            \sum_{l_1,\dots,l_n=1}^k
            \prod_{q=1}^n
            \left(
                \prod_{r:\sigma(r)=q}
                \partial_{l_r}
            \right)
            V_{i_q}^{l_q}(x_q)
            \left(
                \prod_{r:\sigma(r)=n+1}
                \partial_{l_r}
            \right)
            f(x_{n+1}).
        \end{equation}
    \end{itemize}
\end{remark}
We conclude this section with the following proposition, which shows that the abstract Definition~\ref{def:vector_fields_iteration} appears through a recursion, and indeed generalizes the classical vector-field composition~\eqref{eq:classical_VF_level_1}. This result is a fundamental ingredient for the expansions derived in Section~\ref{sec:Volterra_expansions}, and its proof is presented in Appendix~\ref{sec:proofs}.
\begin{proposition}\label{prop:vector_field_algebra}
    Let $(i_2\cdots i_n,\bsigma)\in \Lambda_\AA^{n-1}$ and $\ell \in \AA$. For any $e \in \{1,\dots,n\}^{|\ell|}$, we define
    \begin{equation}\label{eq:crafting}
         \bsigma_+^{e}(1) = \begin{cases}
             (e_1+1,e_1+1) & |\ell|=1, \\
             (e_1+1,e_2+1) & |\ell|=2,
         \end{cases}
         \qquad
         \bsigma_+^e(k)=\bsigma(k-1)+(1,1), \quad 2\leq k \leq n.
    \end{equation}
    Then, for any $(f,V_1,\dots,V_m)$ as before, we have the recursive identity
    \begin{equation}\label{eq:recursion_formula}
         \bV^{(\ell i_2\cdots i_{n},\bsigma_+^{e})}(f)(x_1,\dots,x_{n+1})
         =
         \sum_{l_1 \in [k]^{|\ell|}}
         \mathsf{V}_{\ell}^{l_1}(x_1)
         \Big(
         \prod_{h=1}^{|\ell|}
         \partial_{x_{e_h+1}^{l_1^h}}
         \Big)
         \bV^{(i_2\cdots i_n,\bsigma)}(f)(x_2,\dots,x_{n+1}).
    \end{equation}
    Moreover, for any $w=i_1\cdots i_n \in \cW_\AA^n$ we have
    \begin{equation}\label{eq:relation_VF_classical}
    \bV^{i_1} \circ \cdots \circ \bV^{i_n}(f)(\xi)
    =
    \sum_{\bsigma:(w,\bsigma) \in \Lambda_\AA^n}
    \mathbf{V}^{(i_1\cdots i_n,\bsigma)}(f)(\xi),
    \qquad
    \xi \in \mathbb{R}^k,
    \end{equation}
    where we write $\bV^\mathbf{a}(f)(\xi):=\mathbf{V}^{\mathbf{a}}(f)(\xi,\dots,\xi)$.
\end{proposition}

\subsection{Iterated Volterra integrals}\label{sec:generalized_VSig}
Finally, we introduce the iterated integrals $\bZ^\mathbf{a}$ for $\mathbf{a}\in \Lambda_\cA$ or $\mathbf{a}\in \Lambda_\AA$. We now distinguish between the path-regularity of the driving path $X$, starting with Lipschitz paths which is later extended semimartingales. We will always consider $m$-dimensional Volterra processes, and fix the alphabet $\cA=\{1,\dots,m\}$. 

\subsubsection{Volterra paths with Lipschitz drivers} 
\label{sec:smooth_VP}
\begin{definition}\label{def:VP}
We denote by $\cV^1([0,T];\RR^m)$ 
the space of Volterra paths  $Z:\Delta^2 \rightarrow \mathbb{R}^m$ with Lipschitz drivers, that is 
\begin{equation}\label{eq:Volterra_path_def}
    Z^\tau_t = \int_0^t K(\tau, s)\dd{X_s}, \qquad (t,\tau) \in \Delta^2,
\end{equation}
for some $X\in\mathrm{Lip}([0,T];\RR^d)$
and $K\in\Lkernel$. 
\end{definition}

Using this notion, we are able to define the geometric Volterra signature.

\begin{definition}\label{def:volterra-signature}
Consider a Volterra path $Z \in \Vone$. For any $(s,t)\in \Delta^2$ we define the geometric Volterra signature as
\begin{equation}\label{eq:VSIG_general}
\mathbf Z_{s,t}: \Lambda_\cA \rightarrow \mathbb{R}, \qquad
(i_1\cdots i_n,\sigma)\mapsto \bZ_{s,t}^{i_1\cdots i_n,\sigma}
:=
\int_{\Delta_{s,t}^n}
\dd Z_{r_1}^{r_{\sigma(1)},i_1}\cdots \dd Z_{r_n}^{r_{\sigma(n)},i_n},
\end{equation} where we use the convention $r_{n+1}:=t$ and linearly extend $\bZ$ to the vector space $\cH_\cA$.

\end{definition}
\begin{remark}
    Writing the Volterra path as in \eqref{eq:Volterra_path_def}, the signature $\bZ$ may equivalently be expressed as \[
\bZ_{s,t}^{w,\sigma}
:=
\int_{\Delta_{s,t}^n}
\prod_{l=1}^n  K^{i_n}(r_{\sigma(l)},r_l)\dd X_{r_l}, \qquad (\sigma,w) \in \Lambda_\cA,
    \] where $K^{i}(\tau,r)\dd X_r:= \sum_{j=1}^dK^{ij}(\tau,r)\dd X^j_r.$
\end{remark}

\begin{remark}\label{rem:classical_VSIG}
    Consider the subset of consecutive admissible evaluation maps  $$\Sigma^{\mathrm{con}}:= \Sigma \cap \{\sigma(k)=k+1\}.$$ Restricted to $\Lambda_\cA^{\mathrm{con}}:= \cW_\cA \times \Sigma^{\mathrm{con}}$, the map $\bZ$ in \eqref{eq:VSIG_general} coincides with the Volterra signature introduced in~\cite{HT21}, that is \begin{equation}\label{eq:classical_VSIG}
    \bZ_{s,t}^{i_1\cdots i_n,\sigma}\equiv \bZ_{s,t}^{i_1\cdots i_n}= \int_{\Delta_{s,t}^n}
\dd Z_{r_1}^{r_{2},i_1}\cdots \dd Z_{r_n}^{t,i_n}.
    \end{equation}
\end{remark}
As mentioned in the introduction, our generalization of the original Volterra signature as introduced in \cite{HT21} comes from expanding equation \eqref{eq:volterra_eq}. 
As it turns out, compared to Volterra signature \eqref{eq:classical_VSIG}, our signature is indeed \emph{geometric}, in a sense we shall describe now. Indeed, as already observed in \cite{HT21}, the following example shows that the Volterra signature \eqref{eq:classical_VSIG} is not group-like, that is, the classical \emph{shuffle-identity} does not hold true.

\begin{example}\label{ex:illustration_shuffle}
    Consider two identical, single letter words $w=v=i$. By definition of the shuffle product (Definition~\ref{def:classical_shuffle}) we have $w\shuffle v = 2 ii$, and the shuffle identity of the signature states $$\left (\mathrm{Sig}(X)_{s,t}^{i}\right )^2= \mathrm{Sig}(X)_{s,t}^{i \shuffle i} = 2\cdot \mathrm{Sig}(X)_{s,t}^{ii}.$$ For the Volterra signature, we observe  \begin{equation}
    \left (\bZ_{s,t}^{i}\right)^2:=\left ( \int_s^t \dd Z^{t,i}_u \right)^2 = \int_{[s,t]^2}\dd Z^{t,i}_{u_1} \dd Z^{t,i}_{u_2} = 2\int_{\Delta_{s,t}^{2}} \dd Z^{t,i}_{u} \dd Z_v^{t,i}.
\end{equation}  where we simply used $[s,t]^2 = \{(u,v) \in [s,t]^2:u\leq v \} \cup \{(u,v) \in [s,t]^2:u\geq v \}$ for the last equality. We can notice that the right-hand side is not equal to $2\cdot \bZ_{s,t}^{ii}$, see also \eqref{eq:classical_VSIG}. However, it is part of the geometric Volterra signature in Definition~\ref{def:volterra-signature}, and $(\bZ_{s,t}^{i})^2=  2 \cdot \bZ_{s,t}^{ ii,\sigma},$ where $\sigma \in \Sigma^2$ and $\sigma(k)\equiv 3$.
\end{example}

In the following, we show more generally that geometric Volterra signature in Definition~\ref{def:volterra-signature} satisfies the shuffle identity with respect to $\shuffle$ from Definition~\ref{def:shuffle_vwords}, the proof can be found in Appendix~\ref{sec:proofs}.

\begin{proposition}\label{prop:shuffle_identity}
    For any Volterra path $Z \in \Vone$, we have \begin{equation}\label{eq:shuffle-identity}
        \bZ_{s,t}^{\mathbf{a}} \bZ_{s,t}^{\mathbf{b}} = \bZ_{s,t}^{\mathbf{a} \shuffle \mathbf{b} }, \quad 0\leq s \leq t \leq T, \quad \mathbf{a},\mathbf{b} \in \Lambda_\cA.
    \end{equation}
\end{proposition}

As a by-product, we can show that linear functionals on time-augmented Volterra signature are universal on the space of Volterra paths. More precisely, for the universality we shall work on \begin{equation}\label{eq:time-augmented}
    \widehat{\cV}:= \Big \{ Z\in \cV^1: Z_t^{\tau,1}=t, \,  \Vert Z \Vert_{{\cV}}:= \sup_{t \in [0,T]} \Vert \dot Z^t \Vert_{\mathbb{L}^1}<\infty \Big \}.
\end{equation} We note that, similar to the analysis in \cite{HHPT25}, one could study universality on the driving noise $X$, but we do not include this here. 
\begin{corollary}\label{cor:universal} Consider a continuous functional $F: \cK \rightarrow \mathbb{R}$ where $\cK$ is a compact subset of $\widehat{\cV}([0,T];\RR^{m})$. Then, for any $\varepsilon>0$, there exists $\mathbf{a}\in \cH_\cA$, such that \begin{equation*}
    \sup_{Z\in \cK }  | F(Z)-\bZ^{\mathbf{a}}_{0,T}| \leq \varepsilon.
\end{equation*}
\end{corollary}
\begin{proof}
    The result follows by the Stone-Weierstrass theorem, if we can show that $\{Z \mapsto \bZ^\mathbf{a}_{0,T}: \mathbf{a}\in \cH_\cA \}$ forms a point-separating subalgebra of $C({\cK};\RR)$. Point-separation follows from injectivity restricted to linear functionals on  $\cW\times \Sigma^{\mathrm{con}}$, see Remark~\ref{rem:classical_VSIG}, which was shown in \cite[Corollary 3.6]{HHPT25}. Moreover, it forms an algebra thanks to the shuffle-identity seen in Proposition~\ref{prop:shuffle_identity}. Thus we are left with showing that $Z \mapsto \bZ^\mathbf{a}_{0,T}$ are indeed continuous maps. But this can be seen directly by definition, that is for $Z,Y \in \widehat{\cV}$ and $\mathbf{a} = (i_1\cdots i_n, \sigma)\in \cW^n_{\cA} \times \Sigma^n$, we have  \begin{align*}
        |\bZ_{0,T}^\mathbf{a}-\bY_{0,T}^\mathbf{a}| & \leq \int_{\Delta_{0,T}^{n}} \Big | \prod_{l=1}^{n}\dot Z^{r_{\sigma(l)},i_l}_{r_l}- \prod_{l=1}^{n}\dot Y^{r_{\sigma(l)},i_l}_{r_l}  \Big | \dd r_1 \cdots \dd r_n \\ & \leq \sum_{k=1}^n \int_{\Delta_{T}^{n}} \big |\prod_{l<k} \dot Z^{r_{\sigma(l)},i_l}_{r_l} \prod_{l>k} \dot Y^{r_{\sigma(l)},i_l}_{r_l} \big | \big | \dot Z_{r_k}^{r_{\sigma(k)},i_k}- \dot Y_{r_k}^{r_{\sigma(k)},i_k} \big | \dd r_1\cdots \dd r_n \\ & \leq \Vert Z- Y \Vert_{\cV}\sum_{k=1}^n  \Vert Z \Vert_{\cV}^{k} \Vert Y \Vert_{\cV}^{(n-k)},
    \end{align*} where we simply used a telescoping sum for the second inequality, and then the definition of the norm in \eqref{eq:time-augmented} for the last inequality. The same can be done for a general element $\mathbf{a} \in \cH_\cA$, which shows the claim.
\end{proof}
\begin{remark}\label{rem:link_branched}
As already emphasized in the introduction, the geometric Volterra signature introduced in
Definition~\ref{def:volterra-signature} is  related to the
tree-indexed Volterra lifts developed in~\cite{HTW23,BK23}. Motivated by the theory of branched rough paths~\cite{gubinelli2010ramification} and the resulting Butcher-series expansions, the branched Volterra signature in~\cite[Definition~2.2]{BK23}, is defined as the collection of iterated Volterra integrals, indexed by rooted trees. More precisely, for ladder trees, that is  $$\cT_\cL = \Big \{\begin{tikzpicture}[scale=0.15,baseline=0.05cm]
    \node at (0,1)  [dot,label={[label distance=-0.2em]left: \scriptsize $i$}] (root) {};
\end{tikzpicture}\, , \begin{tikzpicture}[scale=0.15,baseline=0.05cm]
    \node at (0,0)  [dot,label={[label distance=-0.2em]left: \scriptsize $i$}] (root) {};
    \node at (0,2)  [dot,label={[label distance=-0.2em]left: \scriptsize $j$}] (up1) {};
    \draw[kernel1] (up1) to node [sloped,below] {\small } (root);
\end{tikzpicture}\, ,
\begin{tikzpicture}[scale=0.15,baseline=0.05cm]
    \node at (0,0)  [dot,label={[label distance=-0.2em]left: \scriptsize $i$}] (root) {};
    \node at (0,2)  [dot,label={[label distance=-0.2em]left: \scriptsize $j$}] (up1) {};
    \node at (0,4)  [dot,label={[label distance=-0.2em]left: \scriptsize $k$}] (up2) {};
    \draw[kernel1] (up1) to node [sloped,below] {\small } (root);
    \draw[kernel1] (up2) to node [sloped,below] {\small } (up1);
\end{tikzpicture}\, ,\dots \Big \},$$ the
corresponding coordinates are identified with the Volterra
signature~\eqref{eq:classical_VSIG}. For branching trees, the integration
domain changes according to the leaves, for example
$$  
    \begin{tikzpicture}[scale=0.15,baseline=0.05cm]
        \node at (0,0)
            [dot,label={[label distance=-0.2em]below:
            \scriptsize $k$}] (root) {};
        \node at (-1,2)
            [dot,label={[label distance=-0.2em]left:
            \scriptsize $i$}] (left) {};
        \node at (1,2)
            [dot,label={[label distance=-0.2em]right:
            \scriptsize $j$}] (right) {};
        \draw[kernel1] (left) to (root);
        \draw[kernel1] (right) to (root);
    \end{tikzpicture}
    \longmapsto 
    \int_s^t
    \left(\int_s^u \dd Z_r^{u,i}\right)
    \left(\int_s^u \dd Z_q^{u,j}\right)
    \dd Z_u^{t,k}, \qquad i,j,k \in \cA.
$$
Applying the shuffle identity from
Proposition~\ref{prop:shuffle_identity} to the two inner
integrals gives
\begin{align*}
    \left(\int_s^u \dd Z_r^{u,i}\right)
    \left(\int_s^u \dd Z_q^{u,j}\right)
    &=
    \bZ_{s,u}^{ij,\sigma}
    +
    \bZ_{s,u}^{ji,\sigma}, \qquad \sigma(k)\equiv 3.
\end{align*}
After integrating both terms, we find $$  
    \int_s^t
    \left(\int_s^u \dd Z_r^{u,i}\right)
    \left(\int_s^u \dd Z_q^{u,j}\right)
    \dd Z_u^{t,k} = \bZ_{s,t}^{ijk,\sigma^{+}}+\bZ_{s,t}^{jik,\sigma^{+}}, \qquad \sigma^+(1)=\sigma^+(2)=3,\,  \sigma^+(3)=4.
$$
The last step holds more generally, namely for any 
$(w,\sigma)\in\Lambda_\cA^n$ and $k\in\cA$, we have 
$$
    \int_s^t
    \bZ_{s,u}^{w,\sigma}\dd Z_u^{t,k}
    =
    \bZ_{s,t}^{wk,\sigma^+}, \qquad \sigma^+(\ell)=\begin{cases}
        \sigma(\ell) & \ell \leq n \\ n+2 & \ell=n+1
    \end{cases}.
$$
Combining these identities with the recursive nature of the tree-indexed integrals in~\cite[Definition~2.2]{BK23},
one could in principle establish the corresponding representation
for general trees by induction. A precise formulation would require
introducing the precise tree and forest notation, and we do not pursue this generalization here.
\end{remark}

\subsubsection{Volterra semimartingale drivers}

Let us now extend our geometric Volterra signature to semimartingale drivers, more generally we can extend Definition~\ref{def:volterra-signature} to the class of \emph{Volterra semimartingales} adapted to some complete, filtered probability space $(\Omega,\cF,\mathbb{F}=(\mathcal{F}_t)_{t\in [0,T]},\mathbb{P})$. To this end, we will focus on Itô processes of the type \begin{equation}\label{eq:ito_process}
 X_t = X_0+\int_0^tb_s\dd s+ \int_0^t \sigma_s \dd B_s,    
\end{equation} where $B$ is a $d$-dimensional and $\mathbb{F}$-adapted Brownian motion, and $b$ and $\sigma$ are $\mathbb{F}$-adapted satisfying 
\begin{equation}\label{eq:coefficients_ito}
    b_t \in \mathbb{R}^d \quad \text{and}\quad \sigma_t \in \mathbb{R}^{d\times d} \quad \text{such that} \quad \sup_{0\leq t \leq T} \mathbb{E}[|b_t|^p+|\sigma_t\sigma_t^\top|^{p/2}]<\infty, 
\end{equation} for some $p \geq 2$. Such Itô processes are in particular $\mathcal{H}^p$-semimartingales in the sense of \cite[Chapter IV]{Protter2004}, that is
\[
\Vert X \Vert_{\cH^p}:= \Big \Vert \int_0^T |b_t|\dd t \Big \Vert_{L^p} + \Big \Vert \int_0^T|\sigma_t \sigma_t^\top|\dd t  \Big \Vert^{1/2}_{L^{p/2}}<\infty.
\]
In the sequel we shall work with the following class of Volterra processes, which naturally leads to a well-defined extension of the Volterra signature introduced in Definition~\ref{def:volterra-signature}.

\begin{definition}\label{def:V_semiMG}
We define by $\mathcal{V}^p_\mathbb{F}:=\mathcal{V}^p_\mathbb{F}([0,T];\mathbb{R}^m)$ the space of Volterra semimartingales starting from zero, that is \begin{equation}\label{eq:admissible_VP}
      \Delta^2 \ni (t,\tau) \mapsto Z_t^\tau=\int_0^tK(\tau,u) \dd X_u \in \mathbb{R}^m,  
\end{equation} 
where $X$ is a $d$-dimensional semimartingale of the form \eqref{eq:ito_process} satisfying \eqref{eq:coefficients_ito},  and $K\in \Ltwokernel$. Moreover, we define $\mathcal{V}^{\infty-}_\mathbb{F}([0,T];\mathbb{R}^m):= \bigcap_{p\geq 2} \mathcal{V}^p_\mathbb{F}([0,T];\mathbb{R}^m)$.
\end{definition}
The reason for working with the intersection of 
$\cV_\mathbb{F}^p$ spaces, rather than $\cV^\infty_\mathbb{F}$ is motivated in the following example. Throughout this paper, we will make use of the convention~\eqref{eq:convention}, which allows to extend constantly to processes on $[0,T]^2$ setting
\[Z_t^\tau:=Z_{t\wedge\tau}^\tau,
    \qquad t,\tau\in[0,T].
\]

\begin{example}\label{ex:SMG}
    The most relevant examples of Volterra semimartingales in  this article are solutions to stochastic Volterra equations. Indeed, suppose $K\in \Ltwokernel$, see \eqref{eq:kernel_norm}, and denote by $(Y_t)_{t\in [0,T]}$ a weak solution to \[ Y_t=\xi+\int_0^tb(Y_u)K(t,u)\dd u+ \int_0^t\sigma(Y_u)K(t,u)\dd B_u, \quad 0<t\leq T,
    \] which can for instance be obtained for continuous coefficients $(b,\sigma)$ with linear growth conditions, see \cite[Theorem 3.4]{jaber2019affinevolterraprocesses}. Defining the Itô process $\dd X_t:= b(Y_t)\dd t+ \sigma(Y_t)\dd B_t$, it was shown that the condition \eqref{eq:coefficients_ito} holds for all $p\geq 2$ in  \cite[Lemma 3.1]{jaber2019affinevolterraprocesses}. Therefore $Z_t^\tau -\xi := \int_0^tK(\tau,u)\dd X_u \in \mathcal{V}^{\infty-}_\mathbb{F}$, and by definition we have $Z_t^t=Y_t$. However, assuming  \eqref{eq:coefficients_ito} for $p=\infty$ would exclude unbounded coefficients, and therefore we will only work with $\cV_\mathbb{F}^{\infty-}$ in this paper.
\end{example}

The following lemma justifies the name \emph{Volterra semimartingales}, and will be used repeatedly in the sequel when  expanding equations driven by $Z \in \cV^p_\mathbb{F}$. Note, similar to the above, that we set $\cH^{\infty-} \assign \bigcap_{p \geq 2} \cH^{p}$.

\begin{lemma}\label{lem:ito}
    For any $Z \in \mathcal{V}_\mathbb{F}^p$ with $2\leq p <\infty$ or $p = \infty-$ as defined above, and  $0\leq \tau \leq T$, the process $(Z_t^{\tau})_{t\leq \tau}$
 defines a $\mathcal{H}^p$-semimartingale and \begin{equation}\label{eq:V_sup_norm}
     \Vert Z \Vert_{\mathcal{V}_\mathbb{F}^p}:= \sup_{0\leq \tau \leq T}\Vert Z^\tau \Vert_{\mathcal{H}^p}<\infty.
 \end{equation} Moreover, for any $f\in C^2(\RR^m;\RR)$, the following change-of-variales formula holds \begin{equation}
    \label{eq:ito_formula}
        f(Z^\tau_t) = f(Z_0)+ \sum_{i=1}^m\int_0^t\partial_{x_i}f(Z_u^\tau)\dd Z_u^{\tau,i} + \frac{1}{2}\sum_{i,j=1}^m \int_0^t \partial^{2}_{x_ix_j}f(Z_u^\tau) \dd [Z^{\tau,i},Z^{\tau,j}]_u.
    \end{equation}\end{lemma}
\begin{proof}
    Note that for fixed $0\leq \tau \leq T$, by \eqref{eq:coefficients_ito} we have \[
    \mathbb{E}\Big [\int_0^t|K(\tau,u)||b_u| \dd u +\Big (\int_0^t|K(\tau,u)|^2 |a_u| \dd u \Big )^{1/2} \Big ] \leq C \times \Big (\Vert K(\tau,\cdot) \Vert_{L^1}+\Vert K(\tau,\cdot) \Vert_{L^2}), 
    \] where $C= \sup_{0\leq t \leq T}\mathbb{E}[|b_t|+|\sigma_t\sigma_t^\top|^{1/2}]$. By \cite[Theorem 30 and the subsequent Corollary]{Protter2004}, the stochastic integrals $t \mapsto \int_0^tK(\tau,u)\dd X_u$ are well-defined, continuous  semi-martingales. In particular, standard applications of Jensen's inequality shows \begin{align*}
        \Vert Z^\tau \Vert_{\mathcal{H}^p}  & := \Big \Vert \int_0^T |K(\tau,t)b_t| \dd t \Big \Vert_{L^p} + \Big \Vert \int_0^T |K(\tau,t)\sigma_t\sigma_t^\top K(\tau,t)^\top |\dd t \Big \Vert_{L^{p/2}}^{1/2} \\ & \leq 2\big (\Vert K(\tau,\cdot)\Vert_{L^1}+ \Vert K(\tau,\cdot) \Vert_{L^2} \big )\sup_{0\leq t \leq T} (\Vert b_t\Vert_{L^p} + \Vert \sigma_t\sigma_t^\top \Vert_{L^{p/2}}).
    \end{align*} Thus, since $K\in L^{\infty,2}$, the first bracket is finite uniformly in $\tau$, together with the assumption \eqref{eq:coefficients_ito} it follows that $\Vert Z \Vert_{\mathcal{V}^p}<\infty$. Finally, since $Z^\tau_\cdot$ is a semi-martingale, equation \eqref{eq:ito_formula} is a consequence of Itô's formula.
\end{proof}
In what follows we will define an Itô Volterra signature for Volterra semimartingales $Z\in \cV^{\infty-}_{\mathbb{F}}$. As it turns out, in order to study stability (that is, continuity) of the map $Z\mapsto \bZ$, we have to slightly strengthen the topology~\eqref{eq:V_sup_norm}, requiring Hölder continuity in the upper parameter $\tau$. As we shall see in the remark below, this is naturally the case for a large class of potentially singular kernels.

\begin{definition}\label{def:V_hoelder_space}For any $0<\beta\leq 1$ and $p \geq 2$, we define the space $\cV_\mathbb{F}^{\beta,p}=\cV_\mathbb{F}^{\beta,p}([0,T];\mathbb{R}^m)$ consisting of $\beta$-Hölder Volterra semimartingales starting from zero, that is \begin{equation*}\label{eq:hölder_VP}
    \cV_\mathbb{F}^{\beta,p} = \left \{ Z \in \cV_\mathbb{F}^{p}: \Vert Z \Vert_{(\beta,p)}:= \sup_{0 \leq \tau \leq T} \Vert Z^{\tau} \Vert_{\cH^p}+ \sup_{\tau < \eta} \frac{\Vert Z^\tau-Z^\eta\Vert_{\cH^p}}{|\tau-\eta|^{\beta}}\right \}.
\end{equation*} As before, we also define $\cV^{\beta,\infty-}_\mathbb{F}([0,T];\RR^m):= \bigcap_{p\geq 2}\cV_\mathbb{F}^{\beta,p}([0,T];\RR^m)$,  together with the projective limit topology  induced by the metric  \begin{equation}\label{eq:projective_Vp}
d_{\beta}(Z,Y) = \sum_{q\geq 2}2^{-q}\frac{\Vert Z-Y \Vert_{(\beta,q)}}{1+\Vert Z-Y \Vert_{(\beta,q)}}.
\end{equation} Finally, for $\beta=0$ we will  make the abuse of notation $\cV^{0,p}_\mathbb{F}:= \cV^{p}_\mathbb{F}$.
\end{definition}

\begin{remark}
   To get a feeling for the additional Hölder regularity in the upper-parameter, we can repeat the first equation in the proof of Lemma~\ref{lem:ito}, but for the difference $Z^\tau-Z^{\eta} $, to see that $$\Vert Z^\tau-Z^\eta \Vert_{\cH^p} \leq C \big (\Vert K(\tau,\cdot)-K(\eta,\cdot)\Vert_{L^1}+ \Vert K(\tau,\cdot)-K(\eta,\cdot) \Vert_{L^2} \big ).$$
    Therefore we see, a sufficient condition for $\Vert Z \Vert_{(\beta,p)}<\infty$, is $\beta$-Hölder regularity of the $L^2$-norm, that is \[
    \Vert K(\tau,\cdot)-K(\eta,\cdot) \Vert_{L^2} \leq C|\tau-\eta|^{\beta}.
    \]
   Such a condition is natural for kernels with fractional-type singularity, that is $|K(\tau,t)| \lesssim (\tau-t)^{\alpha-1}$ for $\alpha>1/2$. We also note that this condition is common in the literature of Volterra SDEs, see for instance~\cite[Equation (2.5)]{jaber2019affinevolterraprocesses}.
\end{remark}
For readers familiar with stochastic Taylor expansions based on Itô iterated integrals, see \cite{Kloeden1991}, it will not come as a surprise that the expansion of Volterra Itô equations requires an augmentation by quadratic variation terms. In our framework, this amounts to enhancing the Volterra driver by 
\begin{equation}\label{eq:Volterra_bracket}
    (t,\tau_1,\tau_2) \mapsto [Z^{\tau_1},Z^{\tau_2}]_t \in \RR^{m\times m}
\end{equation} In the following definition, also recall  the definition of the extended alphabet $\AA$  in~\eqref{eq:enhanced_alpha}. \begin{definition}\label{def:Ito_Volterra signature}
    For any Volterra semimartingale $Z \in \cV_\mathbb{F}^p([0,T];\RR^m)$, we define for $(s,t)\in \Delta^2$ the Itô Volterra signature 
\begin{equation}\label{eq:VSIG_general_ito}
\mathbf Z_{s,t}: {\Lambda}_\AA \rightarrow \mathbb{R}, \qquad (i_1\cdots i_n,\bsigma) \mapsto 
\bZ_{s,t}^{i_1\cdots i_n,\bsigma}
:=
\int_{\Delta_{s,t}^n}
\dd \cZ_{r_1}^{r_{\bsigma(1)},i_1}\cdots \dd \cZ_{r_n}^{r_{\bsigma(n)},i_n},
\end{equation} where we use the convention $r_{n+1}:=t$, and \begin{equation*}
    \cZ_{r_\ell}^{r_{\bsigma(\ell)},i_\ell} := \begin{cases}
        Z_{r_\ell}^{r_{\bsigma_1(\ell)},i} & \text{ if } i_\ell = i,  \\ 
        [Z^{r_{\bsigma_1(\ell)},i},Z^{r_{\bsigma_2(\ell)},j}]_{r_\ell} & \text{ if } i_\ell= [i,j].
    \end{cases}
\end{equation*} Moreover, we linearly extend $\bZ$ to the vector space ${\cH}_\AA= \mathrm{span}_\mathbb{R}({\Lambda}_\AA)$.
\end{definition}
\begin{example}\label{ex:illu_example_continued}
    We continue Example~\ref{ex:illustration_shuffle}, that is $w=v=i$. An application of Itô's formula \eqref{eq:ito_formula} shows now $$(\bZ^{i}_{s,t})^2 = 2 \int_{\Delta^{(2)}_{s,t}}\dd Z_{u_1}^{t,i}\dd Z_{u_2}^{t,j}+\int_s^t \dd [Z^{t,i}]_u = 2\bZ_{s,t}^{ii,\sigma}+ \int_s^t \dd [Z^{t,i}]_u.$$
    We see that the right-hand side is included in the augmented Volterra signature of Definition~\ref{def:Ito_Volterra signature} above  $\cZ=(Z,[Z,Z])$. In the reminder of this section, we aim to make this relation more precise by establishing a quasi-shuffle identity.  
\end{example}
\begin{remark}
    We see that Definition~\ref{def:Ito_Volterra signature} extends Definition~\ref{def:volterra-signature}, in the sense that if we restrict to words in $\cW_\cA$, we have $\bZ^{w,\bsigma}_{s,t}=\bZ_{s,t}^{w,\bsigma_1}$ for $(w,\bsigma_1)\in \cW_\cA^n\times \Sigma^n$, which corresponds to the definition in \eqref{eq:VSIG_general}.
\end{remark}

In the following theorem we show that the Itô Volterra signature has moments of all orders and is continuous on the underlying space of Volterra semimartingales. The proof is presented in Appendix~\ref{sec:proofs}. 

\begin{theorem}\label{lem:finite_moments}
Let $\beta \in [0,1]$ and $Z\in\cV_{\mathbb F}^{\beta,\infty-}$.
Then the Itô Volterra signature $\bZ$ satisfies
\[
    \sup_{(s,t)\in\Delta^2}
    \mathbb E\big[|\bZ_{s,t}^{\mathbf a}|^p\big]
    <\infty \quad \text{for all }  \mathbf{a} \in\Lambda_\AA \text{ and } p\geq 1.
\]
Moreover, for $0<\beta \leq 1$ and every $\mathbf a\in\Lambda_\AA$ and $(s,t)\in \Delta^2$, the map
\[
    \bZ^\mathbf{a}: \cV_{\mathbb{F}}^{\beta,\infty-} \longrightarrow L^{\infty-}, \qquad Z\longmapsto
    \bZ_{s,t}^{\mathbf a}
\]
is continuous
with respect to the projective limit topologies. More precisely, it is locally Lipschitz into $L^p$ for every $p\geq 1$. 
\end{theorem}

We conclude this section with the algebraic analogue to the shuffle-identity in Proposition~\ref{prop:shuffle_identity} for Volterra paths with Lipschitz drivers, the proof is postponed to Appendix~\ref{sec:proofs}.
\begin{proposition}\label{prop:quasi_shuffle} For any $Z\in\cV^{\beta,\infty-}_{\mathbb F}$ and $\mathbf a,\mathbf b\in\Lambda_\AA$, we have almost surely
$$
    \bZ_{s,t}^{\mathbf a}\bZ_{s,t}^{\mathbf b}
    =
    \bZ_{s,t}^{\mathbf a\star\mathbf b}, \qquad 0 \leq s \leq t \leq T,$$
where $\star$ denotes the quasi-shuffle product from Definition~\ref{def:Sigma_quasi_shuffle}.
\end{proposition}

Similar as in the smooth case, we end this section with an universality result on the Itô Volterra signature on the space of semimartingales. We note that, compared to the pathwise and locally uniform convergence in  Corollary~\ref{cor:universal}, here we consider a weaker, probabilistic universality result,  see also the discussion in Remark~\ref{rem:V_RP_universal}. To do so, we consider time-augmented Volterra semimartingales, that is we restrict to \[
\widehat{\cV}_\mathbb{F}^{\beta,p}:= \Big \{ Z\in \cV_\mathbb{F}^{\beta,p}: Z_t^{\tau,1}=t \Big \}.
\]

\begin{corollary}\label{cor:SMG_universality} Let $Z\in \widehat{\cV}^{\beta,p}_\mathbb{F}$ for some $\beta\in [0,1]$ and $p\geq 2$. Then, for any Borel measurable function $F:C(\Delta^2_{0,T};\mathbb{R}^m)\rightarrow \mathbb{R}$, we can find a sequence $(\ell_n)_{n\in \mathbb{N}}$ in $\cH_{\AA}$, such that  \[
\lim_{n\to \infty}\bZ_{0,T}^{\ell_n} = F(Z)\quad  \text{ almost surely.}
\]
\end{corollary}
\begin{proof}

As shown in the proof of Proposition~\ref{prop:quasi_shuffle},
the Volterra process $(t,\tau)\mapsto Z_t^\tau$ admits an
almost surely jointly continuous modification. We may therefore
work on the canonical probability space
$(E,\mathcal B(E),\mathbb P_Z)$, where
$E=C(\Delta_{0,T}^2;\mathbb R^m)$,
$\mathcal B(E)$ denotes the Borel $\sigma$-algebra induced
by the supremum norm, and $\mathbb P_Z$ is the law of $Z$. We first want to establish
$\mathcal B(E)=\sigma(\mathcal S^\ell:\ell\in\Lambda_\AA)$,
where $\mathcal S^\ell:E\to\mathbb R$ are Borel measurable
representatives of the geometric Volterra signature coordinates $\bZ^\ell_{0,T}$, which can be defined as follows. For each $\ell\in\Lambda_\AA$, consider Borel
measurable approximations $I_n^\ell:E\to\mathbb R$ such that
$I_n^\ell(Z)\to\bZ_{0,T}^\ell$ in probability.
Passing to a subsequence $(n_k)_{k\in \mathbb{N}}$, we can assume this
convergence holds almost surely. We then define
$\mathcal S^\ell(x):=\lim_{k\to\infty}I_{n_k}^\ell(x)$
whenever the limit exists and is finite, and set
$\mathcal S^\ell(x):=0$ otherwise.
In particular, 
$\mathcal S^\ell$ are  Borel measurable maps satisfying
$\mathcal S^\ell(Z)=\bZ_{0,T}^\ell$ almost surely.
As a consequence, we have 
$\sigma(\mathcal S^\ell:\ell\in\Lambda_\AA)
\subseteq\mathcal B(E)$. For the other direction, we can first restrict to $\ell_{n,m}^i=(w_{n,m}^i,\bsigma)$ of the type 
\[
    w_{n,m}^i
    :=
    \underbrace{1\cdots1}_{n\text{ times}}\,
    i\,
    \underbrace{1\cdots1}_{m\text{ times}}, \quad \bsigma(k)=k+1, \qquad \mathcal{S}^{w_{n,m}^i,\bsigma}(Z)
    =
    \int_0^T
        \frac{(T-t)^{m-1}}{(m-1)!}
        \left(
            \int_0^t \frac{s^n}{n!}\,\dd Z_s^{t,i}
        \right)\dd t
    \qquad\text{a.s}. 
\]
We observe that the last integral can be defined for any continuous   $z\in E$, and we denote it by $\varphi_{n,m}^i(z)$. Now suppose $z,y \in E$ are such that
$\varphi_{n,m}^i(z)=\varphi_{n,m}^i(y)$ for all $i\in \cA$ and $n,m \geq 1$, and write  $D:=z-y$. Setting $F_n(t)= \int_0^t \frac{s^n}{n!}\dd D_s^{t,i}$, the previous identity we know
\[
    \int_0^T
        \frac{(T-t)^{m-1}}{(m-1)!}\,F_n(t)\,\dd t
    =0,
    \quad n,m\geq1.
\]
But since polynomials are dense in $L^2$ and $F_n$ is continuous, it follows that $F_n(t)=0$ for all $n\geq1$ and almost everywhere. Using integration by parts
and again density of polynomials also shows 
\[
    F_n(t)
    = \frac{1}{(n-1)!}
    \int_0^t s^{n-1}(D_t^{t,i}-D_s^{t,i})\,\dd s
    =0 \quad \forall n\geq1
    \quad\Longrightarrow\quad
    D_s^{t,i}=D_t^{t,i}
    \quad\forall s\in[0,t],
\]
which holds $\dd t$-almost everywhere,
 using that $D^t_\cdot$ is in particular continuous. Finally since $D_0^{t,i}=0$, we conclude that $D_s^{t,i}=0$. Thus $\Phi(z)=(\varphi_{n,m}^i(z): i \in \cA, \, n,m \geq 1 )$ is injective and Borel measurable. By the Lusin–Souslin theorem, it follows that $\cB = \sigma(\varphi_{n,m}^i:n,m\geq 1,\, i \in \cA )\subseteq \sigma(\cS^\ell:\ell \in \Lambda_{\AA})$, and thus $\cB=\sigma(\cS^\ell:\ell \in \Lambda_{\AA})$. We can now conclude the result by the functional monotone class theorem. In particular, define $\cM= \{\cS^{\ell}: \ell \in \cH_{\AA} \}$, which forms an unital algebra thanks to the quasi-shuffle identity in Proposition~\ref{prop:quasi_shuffle}, and let $\cY$ be the closure of $\cM$ under convergence in probability. Since  $\sigma(\cM)= \cB$, by the monotone class theorem~\cite[Theorem 6.3]{jacod2004probability} it follows that $\cY$ contains every $\cB$-random variable, which concludes the proof.
\end{proof}

\begin{remark}\label{rem:V_RP_universal}
The universality obtained in Corollary~\ref{cor:SMG_universality}
is of probabilistic nature: for a fixed Volterra semimartingale
$Z$, we can approximate measurable functionals
$F(Z_t^\tau:0\leq t\leq\tau\leq T)$ almost surely by
linear functionals of $\bZ$. This is in contrast to the usual
universality results for signatures, such as
Corollary~\ref{cor:universal}, which establish locally uniform
approximation on the space of driving paths, and thus in
particular pathwise approximation for stochastic drivers
on compact sets. The main ingredient missing here is continuity of the
signature lift $Z\mapsto\bZ$ on the underlying path space.
A possible approach would be to introduce a suitable
Volterra rough path space, in which $Z$ is enhanced to
a rough path $\mathbb{Z}$ and the signature $\bZ$ is obtained
as a continuous transformation of $\mathbb Z$.
Provided that the signature coordinates also separate
points of this enhanced space, locally uniform universality
for continuous functionals could then be obtained by
the Stone-Weierstrass theorem, thanks to the quasi-shuffle property in Proposition~\ref{prop:quasi_shuffle}. A natural starting point would be the Volterra rough path
spaces introduced in~\cite{HT21,HTW23,HTW24,BK23}, constructed based on the singularity of the kernel $K$.
However, building our geometric Volterra signature
as a pathwise transformation of a Volterra
rough path is beyond the scope of this work.
\end{remark}

\section{Taylor expansions for Volterra equations}\label{sec:Volterra_expansions}
In this section we present the main expansions to Volterra equations of the type\begin{equation}\label{eq:volterra_eq_2}
Y_t = \xi + \sum_{i=1}^m\int_0^t V_i(Y_s)\dd Z_s^{t,i}, \quad 0 \leq t \leq T, \qquad \xi \in \mathbb{R}^k
\end{equation} where  $V_i:\RR^k \rightarrow \RR^k$ are smooth vector-fields $V_i\in C_b^\infty(\RR^k;\RR^k)$. We begin with $Z\in \Vone$ introduced in  Section~\ref{sec:smooth_VP}.
\subsection{Series representations for Volterra CDEs}\label{sec:cde_expansion}
In the smooth setting, we work with the alphabet of single letters $\cA=\{1,\dots,m\}$. Let us briefly recall here the well-established case $Z_t^\tau \equiv X_t$, for some smooth path $X:[0,T]\rightarrow \RR^m$. In this case, recalling the first-order operator notation from~\eqref{eq:classical_VF_level_1}, a simple change of variables shows that $$f(Y_t)=f(\xi)+\sum_{l=1}^k\sum_{i=1}^m\int_0^tV_{i}^l(Y_s)\partial _l f(Y_s)\dd X^{i}_s = f(\xi)+ \sum_{i=1}^m \int_0^t \mathbf{V}^i(f)(Y_s)\dd X_s^i.$$
Doing the same for $\mathbf{V}^i(f):\RR^k \rightarrow \mathbb{R}$ inside the integral on the right-hand side, leads to $$f(Y_t)=f(\xi)+ \sum_{i=1}^m\mathbf{V}^i(f)(\xi)(X_t^i-X_0^i)+ \sum_{i_1,i_2=1}^m\int_0^t\int_0^u \mathbf{V}^{i_1}(\bV^{i_2}(f))(Y_s) \dd X_s^{i_1}\dd X_u^{i_2}.$$ Repeating this procedure naturally leads to the standard expansion~\eqref{eq:classical_expansion}, with respect to the operator recursion $\bV^{i_1}\circ \cdots \circ  \bV^{i_n}$ in \eqref{eq:classical_VF_iteration}.

Returning to the Volterra framework, let $Y$ be the solution to \eqref{eq:volterra_eq_2} with respect to some Volterra paths $Z \in \Vone$ (see Definition~\ref{def:VP}). We will write $Y_t= Y_t^t$, where $$\Delta^2 \ni (t,\tau) \mapsto Y_t^\tau = \xi+ \sum_{i=1}^m\int_0^tV_i(Y_u)\dd Z_u^{\tau,i}.$$ Now since $Z \in \cV^1$, the map $u \mapsto Y_u^t$ is Lipschitz on $[0,t]$ and we can apply a change of variables in the lower variable to see \begin{equation}\label{eq:cov_smooth}
    f(Y_t) = f(\xi)+ \sum_{l=1}^k\int_0^t \partial_lf(Y_s^t) \dd Y_s^{t,l} = f(\xi)+ \sum_{l=1}^k \sum_{i=1}^m\int_0^t V_{i}^l(Y_s^s)\partial_lf(Y_s^t) \dd Z_{s}^{t,i}.
\end{equation}
In particular, recalling the definition of the first order vector-field operators~\eqref{eq:vector_fields_first_level}\begin{equation}\label{eq:exp_first_order}
f(Y_t)=f(\xi)+ \sum_{i=1}^m \int_0^t\mathbf{V}^{i}(f)(Y_s^s,Y_s^t) \dd Z_{s}^{t,i}.
\end{equation}
Applying another change of variables now to $s \mapsto \mathbf{V}^{i}(f)(Y_s^s,Y_s^t)$, we find \begin{align*}
    \mathbf{V}^i(f)(Y_s^s,Y_s^t) & = \mathbf{V}^i(f)(\xi)+ \sum_{l=1}^k \int_0^s \partial_l \mathbf{V}^i(f)(Y_u^s,Y_u^t) \dd Y_u^{s,l}+ \sum_{l=1}^{k} \int_0^s\partial_{k+l}\mathbf{V}^i(f)(Y_u^s,Y_u^t) \dd Y_u^{t,l} \\ &
    = \mathbf{V}^i(f)(\xi)+ \sum_{l=1}^k \sum_{j=1}^m \int_0^t \partial_l \mathbf{V}^i(f)(Y_u^s,Y_u^t) V_{j}^l(Y^u_u) \dd Z_{u}^{s,j} \\ & \qquad \qquad \qquad \, + \sum_{l=1}^{k}\sum_{j=1}^m \int_0^s\partial_{k+l}\mathbf{V}^i(f)(Y_u^s,Y_u^t)V_{j}^l(Y_u^u)\dd Z_u^{t,j}.
\end{align*}
Now, with our definition of vector-fields at hand (see  Definition~\ref{def:vector_fields_iteration}), its recursive property \eqref{eq:recursion_formula} for single letters reveals the integrands $\bV^{(ij,\sigma_1)}(f)(Y_u^u,Y_u^s,Y_u^t)$ resp. $\bV^{(ij,\sigma_2)}(f)(Y_u^u,Y_u^s,Y_u^t)$ on the right-hand side, for the admissible maps $\sigma_1(k)=k+1$ and $\sigma_2(k)\equiv 3$ (see Definition~\ref{def:admissible-maps}). On the other-hand, when plugging the last identity into the first order expansion \eqref{eq:exp_first_order}, we get two remainder integrals with integrators $\dd Z^{s,j}_{u_1}\dd Z^{t,i}_{u_2}$ and $\dd Z_{u_1}^{t,j}\dd Z_{u_2}^{t,i}$. We can label them by the same elements $(ij,\sigma_1)$ and $(ij,\sigma_2)$ in $\Lambda_\cA^2$. In summary, since $\{\sigma_1,\sigma_2\}$ are the only admissible maps on level $2$, we have \begin{equation}\label{eq:N1_expansions}
f(Y_t)=f(\xi)+ \sum_{i=1}^m \mathbf{V}^i(f)(\xi)\bZ_{0,t}^i+ \sum_{(ij,\sigma)\in \Lambda^2_\cA} \int_{\Delta^2_{0,t}} \mathbf{V}^{(ij,\sigma)}(f)(Y_{r_1}^{r_1},Y_{r_1}^{r_2},Y_{r_1}^t) \dd Z^{r_{\sigma(1)},i}_{r_1} \dd Z_{r_2}^{r_{\sigma(2)},j}.
\end{equation}
This procedure can be generalized to any level $n\in \mathbb{N}$, which leads to our  series representation for Volterra CDEs.

\begin{theorem}\label{thm:main_cde_expansion}
    Let $Z\in \Vone$ and  $V_1,\dots,V_m\in C_b^\infty(\RR^k;\RR^k)$, and suppose $Y$ solves the Volterra CDE \eqref{eq:volterra_eq_2}. Then, for any $f\in C_b^{N+1}(\RR^k;\RR)$ and $N\in \mathbb{N}$, we have \begin{equation}\label{eq:taylor_volterra_cde}
    f(Y_t) =f(\xi) + \sum_{n=1}^N\sum_{\mathbf{a}\in \Lambda_\cA^n} \bV ^{\mathbf{a}}(f)(\xi){\mathbf{Z}}^{\mathbf{a}}_{0,t} + R_\cA^{N}(t),
\end{equation} where the remainder is given by \begin{equation}\label{eq:reminder_taylor}
R_\cA^{N}(t)= \sum_{(w,\sigma)\in \Lambda_\cA^{N+1}} \int_{\Delta_{0,t}^{(N+1)}}\bV^{(w,\sigma)}(f)(Y_{r_1}^{r_1},Y_{r_1}^{r_2},\dots,Y_{r_1}^{t}) \prod_{l=1}^{N+1}\dd Z_{r_l}^{r_{\sigma(l)},i_l}, \qquad r_{N+2}:=t.
\end{equation}
\end{theorem}
\begin{proof}
    We prove the claim by induction on $N\in\mathbb{N}$.
    The case $N=1$ was treated at the beginning of this section;
    see \eqref{eq:N1_expansions} and its derivation.
    Suppose now that \eqref{eq:taylor_volterra_cde} holds for $N-1$.
    By linearity, it suffices to show that
    \[
        R_\cA^{N-1}(t)
        =
        \sum_{\mathbf a\in\Lambda_\cA^N}
        \bV^{\mathbf a}(f)(\xi)\bZ_{0,t}^{\mathbf a}
        +
        R_\cA^N(t).
    \]

    For each $(w,\sigma)\in\Lambda_\cA^N$, set $r_{N+1}:=t$.
    Applying the change-of-variables formula on $[0,r_1]$ and
    substituting the defining expression for the differentials
    of $Y$ yields
    \begin{equation}\label{eq:proof_onestep_taylor}
        \begin{aligned}
            &\bV^{(w,\sigma)}(f)
            \bigl(
                Y_{r_1}^{r_1},\dots,Y_{r_1}^{r_N},Y_{r_1}^{t}
            \bigr)
            \\
            &\quad =
            \bV^{(w,\sigma)}(f)(\xi)
            +
            \sum_{e=1}^{N+1}
            \sum_{l_1=1}^k
            \sum_{i_1=1}^m
            \int_0^{r_1}
            \bigl(
                \partial_{x_e^{l_1}}\bV^{(w,\sigma)}(f)
            \bigr)
            \bigl(
                Y_u^{r_1},\dots,Y_u^{r_N},Y_u^t
            \bigr)
            V_{i_1}^{l_1}(Y_u^u)
            \,\dd Z_u^{r_e,i_1}.
        \end{aligned}
    \end{equation}
    Inserting \eqref{eq:proof_onestep_taylor} into $R_\cA^{N-1}(t)$
    and shifting the integration variables by one, we identify
    the sum over $l_1$ using the recursive
    identity~\eqref{eq:recursion_formula}. This gives
    \[
        R_\cA^{N-1}(t)
        =
        \sum_{\mathbf a\in\Lambda_\cA^N}
        \bV^{\mathbf a}(f)(\xi)\bZ_{0,t}^{\mathbf a}
        +
        \widetilde R^N(t),
    \]
    where
    \begin{align*}
        \widetilde R^N(t)
        &=
        \sum_{i_1\cdots i_{N+1}\in\cW_\cA^{N+1}}
        \sum_{\sigma\in\Sigma^N}
        \sum_{e=1}^{N+1}
        \int_{\Delta_{0,t}^{N+1}}
        \bV^{(i_1\cdots i_{N+1},\sigma_+^e)}(f)
        \bigl(
            Y_{r_1}^{r_1},\dots,
            Y_{r_1}^{r_{N+1}},Y_{r_1}^{t}
        \bigr)
        \prod_{l=1}^{N+1}
        \dd Z_{r_l}^{r_{\sigma_+^e(l)},i_l},
    \end{align*}
    with $r_{N+2}:=t$.
    Here $\sigma_+^e$ is the single-letter version of the
    extension defined in \eqref{eq:crafting}. But by \eqref{eq:sets_recursion}, restricted to the alphabet $\cA$,
    every map in $\Sigma^{N+1}$ occurs exactly once as
    $\sigma_+^e$, with $\sigma\in\Sigma^N$ and
    $e\in\{1,\dots,N+1\}$. Re-indexing the last expression therefore
    gives $\widetilde R^N(t)=R_\cA^N(t)$, which completes the induction.
\end{proof}
\begin{corollary}\label{cor:linear}Suppose now $Z\in \cV^1([0,T];\RR^m)$ such that $\tau \mapsto Z^\tau$ is constant, and write $X_t\equiv Z_t^\tau$. In this case, under the same assumptions on $(f,V_1,\dots,V_m)$ as in Theorem~\ref{thm:main_cde_expansion}, we have \begin{equation}
        f(Y_t)=f(\xi)+ \sum_{n=1}^N \sum_{w=i_1\cdots i_n}  \bV^{i_1\cdots i_n}(f)(\xi)\mathrm{Sig}(X)^w_{0,t} + R_N(t),
    \end{equation} where we write $ \bV^{i_1\cdots i_n}(f):=\bV^{(i_1)} \circ \cdots \circ \bV^{(i_n)}(f)$.
\end{corollary}
\begin{proof}
    This is an immediate consequence of the relation~\eqref{eq:relation_VF_classical} in Proposition~\ref{prop:vector_field_algebra}, the fact that the geometric Volterra signature coincides with the classical Volterra signature if $Z$ is constant in the upper-parameter by Definition~\ref{def:volterra-signature}, and by the expansion~\eqref{eq:taylor_volterra_cde}.
\end{proof}

\begin{remark}\label{rem:linear_exp}
If the vector-fields $V_1,\dots,V_m$ are linear and $f(x)=x_i$,  then the expansion~\eqref{eq:taylor_volterra_cde} simplifies to the $i$-th component of the  expansion derived in \cite[Proposition 2.38]{HHPT25}, with respect to the linear Volterra signature~\eqref{eq:classical_VSIG}. To see this, at level $n \in \mathbb{N}$ observe that whenever $\sigma \notin \Sigma^n_{\mathrm{con}}=\{ \sigma \in \Lambda^n: \sigma(k)=k+1\}$, there exists an $ i \in \{ 3,\dots,{n+1}\}$, such that $|\{\sigma(r)=i \}| >1$. In this case, by linearity of $V$ and $f$, the quantities  \eqref{eq:elementary_vectorfield} vanish and thus the expansions can be reduced over $\Sigma_{\mathrm{con}}$, and we can conclude by  Remark~\ref{rem:classical_VSIG}.
\end{remark}

\begin{remark}
In Theorem~\ref{thm:main_cde_expansion}, we expand the difference
$f(Y_t)-f(Y_0)$. The same argument can be started at any time
$0<s\leq t$. A minor but important distinction is that, for fixed terminal
parameter $t$, the expansion on $[s,t]$ starts from $f(Y_s^t)$ rather than
from $f(Y_s)=f(Y_s^s)$. Consequently, the difference $f(Y_t)-f(Y_s)$
contains an additional term accounting for the dependence on the past. Moreover, the dependence on the upper parameters no longer allows for a
multiplicative separation of the vector-field operators and the iterated
integrals. Indeed, when expanding around time zero, the coefficients satisfy
$\mathbf{V}^{(w,\sigma)}(f)(Y_0^{r_1},\dots,Y_0^{r_n},Y_0^t)
=\mathbf{V}^{(w,\sigma)}(f)(\xi)$ and can therefore be factored out of the
integrals. When the expansion is started at time $s$, however, the
coefficients
$\mathbf{V}^{(w,\sigma)}(f)(Y_s^{r_1},\dots,Y_s^{r_n},Y_s^t)$
depend non-trivially on the upper parameters. Nevertheless, one obtains the
generalized expansion
\[
    f(Y_t)-f(Y_s)
    =f(Y_s^t)-f(Y_s)
    +\sum_{n=1}^N\sum_{\mathbf{a}\in \Lambda_\cA^n}
    [\mathbf{V}_s(f)\star \mathbf{Z}]^\mathbf{a}_{s,t}
    +R_\cA^N(s,t),
\]
where, for $\mathbf a=(i_1\cdots i_n,\sigma)$, we use the notation
\[
    [\mathbf{V}_s(f)\star \mathbf{Z}]^\mathbf{a}_{s,t}
    :=
    \int_{\Delta^n_{s,t}}
    \mathbf{V}^{\mathbf{a}}(f)
    (Y_s^{r_1},\dots,Y_s^{r_n},Y_s^t)
    \dd Z_{r_1}^{r_{\sigma(1)},i_1}
    \cdots
    \dd Z_{r_n}^{r_{\sigma(n)},i_n}.
\]
In the classical case $Z_t^\tau\equiv X_t$, we have
$Y_s^\tau=Y_s$ for every $\tau\geq s$. Hence the additional term
vanishes and the vector-field coefficient can again be factored out, yielding
\[
    [\mathbf{V}_s(f)\star \mathbf{Z}]^{(w,\sigma)}_{s,t}
    =
    \mathbf{V}^{(w,\sigma)}(f)(Y_s)
    \mathrm{Sig}(X)^w_{s,t}.
\]
\end{remark}

\subsection{Series representations for Volterra SDEs}\label{sec:stochastic_expansion} 
In this section we generalize Theorem~\ref{thm:main_cde_expansion} to the class of  Volterra semimartingales $Z\in \cV_\mathbb{F}^p$ introduced in Definition~\ref{def:V_semiMG}. 
In order to derive a series representation for strong solutions to the general Volterra SDE~\eqref{eq:volterra_eq_2}, we follow the same recipe as in the smooth case, but moving from $\cA$ to the enhanced alphabet $\AA$. 
The reason is easy to appreciate, since we will use the Itô-formula from Lemma~\ref{lem:ito}, which applied to the lower variable of $(Y_t^\tau)_{(t,\tau)\in \Delta^2}$ (see also Example~\ref{ex:SMG}) reads \begin{equation*}
    f(Y_t)=f(\xi)+ \sum_{l=1}^k \sum_{i=1}^m \int_0^tV_{i}^l(Y^s_s)\partial_{l}f(Y^t_s) \dd Z_{s}^{t,i} + \frac{1}{2}\sum_{p,q=1}^k \sum_{i,j=1}^m\int_0^tV_{i}^pV_{j}^q(Y^s_s)\partial_{pq}^{2}f(Y^t_s)\dd [Z^{t,i},Z^{t,j}]_s.
\end{equation*}
We have already prepared the enhancement of the Volterra path to $\cZ=(Z,[Z])$ in Definition~\ref{def:Ito_Volterra signature}, with respect to the the enhanced alphabet $\mathbb{A}$ in \eqref{eq:enhanced_alpha}. Together with the vector-field notation from Definition~\ref{def:vector_fields_iteration}, the one-step expansion can be written in the  compact form \begin{equation*}
    f(Y_t)=f(\xi) + \sum_{a\in \mathbb{A}}\int_0^t\mathbf{V}^{a}(f)(Y_s^s,Y_s^t)\dd \cZ_s^{t,a}.
\end{equation*}
Similar as in the smooth case, we can iterate this procedure, leading to the following  stochastic Taylor expansion of the Volterra equations.

\begin{theorem}\label{thm:main_sde_expansion}
    Consider a Volterra semimartingale $Z\in \mathcal{V}_\mathbb{F}^{\infty-}([0,T];\RR^m)$ (see  Definition~\ref{def:V_semiMG}) and  $V_1,\dots,V_m\in C_b^\infty(\RR^k;\RR^k)$. For the strong solution $Y$ to \eqref{eq:volterra_eq_2}, and any $f\in C_b^{2N+2}(\RR^k;\RR)$ with $N\in \mathbb{N}$, we have \begin{equation}\label{eq:taylor_volterra_sde}
    f(Y_t) =f(\xi) + \sum_{n=1}^N\sum_{\mathbf{a}\in {\Lambda}_\AA^n} \bV ^{\mathbf{a}}(f)(\xi){\mathbf{Z}}^{\mathbf{a}}_{0,t} + R_{\AA}^{N}(t),
\end{equation} where the remainder is given by \[
R_{\AA}^{N}(t)= \sum_{(w,\boldsymbol{\sigma})\in {\Lambda}_\AA^{N+1}} \int_{\Delta_{0,t}^{(N+1)}}\bV^{(w,\boldsymbol{\sigma})}(f)(Y_{r_1}^{r_1},Y_{r_1}^{r_2},\dots,Y_{r_1}^{t}) \prod_{l=1}^{N+1}\dd \cZ_{r_l}^{r_{\bsigma(l)},i_l}, \qquad r_{N+2}:=t,
\] and $\bZ$ denotes the lift defined in Definition~\ref{def:Ito_Volterra signature}.
\end{theorem}

\begin{proof}
    We argue again by induction on $N$, following the same recipe as the proof of Theorem~\ref{thm:main_cde_expansion}. The case $N=0$ follows from the one-step expansion at the
    beginning of this section. Suppose that the assertion holds
    for $N-1$. For each $(w,\bsigma)\in\Lambda_\AA^N$, set $r_{N+1}:=t$.
    Applying Itô's formula~\eqref{eq:ito_formula} on $[0,r_1]$
    to
    $\bV^{(w,\bsigma)}(f)
    (Y_u^{r_1},\dots,Y_u^{r_N},Y_u^t)$
    and substituting the differentials of $Y$ and their quadratic
    covariations yields
    \begin{equation}\label{eq:proof_onestep_taylor_2}
        \begin{aligned}
            &\bV^{(w,\bsigma)}(f)
            \bigl(
                Y_{r_1}^{r_1},\dots,Y_{r_1}^{r_N},Y_{r_1}^{t}
            \bigr)
            \\
            &\quad =
            \bV^{(w,\bsigma)}(f)(\xi)
            \\
            &\qquad+
            \sum_{e=1}^{N+1}
            \sum_{p=1}^k
            \sum_{i=1}^m
            \int_0^{r_1}
            V_i^p(Y_u^u)
            \bigl(
                \partial_{x_e^p}\bV^{(w,\bsigma)}(f)
            \bigr)
            \bigl(
                Y_u^{r_1},\dots,Y_u^{r_N},Y_u^t
            \bigr)
            \,\dd Z_u^{r_e,i}
            \\
            &\qquad+
            \sum_{e,h=1}^{N+1}
            \sum_{p,q=1}^k
            \sum_{i,j=1}^m
            \int_0^{r_1}
            \mathsf V_{[i,j]}^{p,q}(Y_u^u)
            \bigl(
                \partial_{x_e^p}\partial_{x_h^q}
                \bV^{(w,\bsigma)}(f)
            \bigr)
            \bigl(
                Y_u^{r_1},\dots,Y_u^{r_N},Y_u^t
            \bigr)
            \,\dd\cZ_u^{r_e,r_h,[i,j]}.
        \end{aligned}
    \end{equation}
    Here
    $\cZ^{r_e,r_h,[i,j]}=[Z^{r_e,i},Z^{r_h,j}]$,
    and the factor $1/2$ from Itô's formula is included in
    $\mathsf V_{[i,j]}^{p,q}$ by~\eqref{eq:generalized_vfields}. Inserting \eqref{eq:proof_onestep_taylor_2} into
    $R_{\AA}^{N-1}(t)$ and shifting the integration variables
    by one, we apply the recursive
    identity~\eqref{eq:recursion_formula} to both terms.
    The first-order term corresponds to a single new letter
    with one attachment index, while the second-order term
    corresponds to a double new letter with two attachment indices.
    Thus,
    \[
        R_{\AA}^{N-1}(t)
        =
        \sum_{\mathbf a\in\Lambda_\AA^N}
        \bV^{\mathbf a}(f)(\xi)\bZ_{0,t}^{\mathbf a}
        +
        \widetilde R^N(t),
    \]
    where
    \begin{align*}
        \widetilde R^N(t)
        &=
        \sum_{(i_2\cdots i_{N+1},\bsigma)\in\Lambda_\AA^N}
        \sum_{i_1\in\AA}
        \sum_{e\in[N+1]^{|i_1|}}
        \int_{\Delta_{0,t}^{N+1}}
        \bV^{(i_1\cdots i_{N+1},\bsigma_+^e)}(f)
        \bigl(
            Y_{r_1}^{r_1},\dots,
            Y_{r_1}^{r_{N+1}},Y_{r_1}^{t}
        \bigr)
        \prod_{l=1}^{N+1}
        \dd\cZ_{r_l}^{r_{\bsigma_+^e(l)},i_l},
    \end{align*}
    with $r_{N+2}:=t$ and $\bsigma_+^e$ as defined
    in~\eqref{eq:crafting}. Again by \eqref{eq:sets_recursion}, each
    $(w,\btau)\in\Lambda_\AA^{N+1}$ is obtained uniquely by
    choosing its first letter $i_1$, a compatible tail
    $(i_2\cdots i_{N+1},\bsigma)\in\Lambda_\AA^N$,
    and $e\in[N+1]^{|i_1|}$ such that
    $\btau=\bsigma_+^e$.
    Re-indexing the last expression therefore gives
    $\widetilde R^N(t)= R_{\AA}^N(t)$,
    which completes the induction.
\end{proof}

\subsection{Extension to path-dependent functionals}
\label{sec:functional_expansions}
The expansions presented in the previous sections can, under suitable
regularity assumptions, be extended to path-dependent functionals $f$.
For $K\equiv\mathrm{Id}$, this was done in~\cite{Litterer2011OnAC},
building on Dupire's functional It\^o calculus~\cite{dupire2019functional},
and we closely follow here their presentation. We consider
non-anticipative functionals
$$
    F:[0,T]\times C([0,T];\RR^k)\to\RR,
$$
and, for every fixed $t\in[0,T]$, aim to expand
$F(t,Y^t_{\cdot\wedge t})$. Here we write
$Y^t_{\cdot\wedge s}:=(Y^t_{u\wedge s})_{u\in[0,T]}$ for
$0\leq s\leq t$: the superscript $t$ denotes the fixed upper
Volterra parameter, while stopping in the lower variable is
indicated in the subscript. The central idea is again that,
thanks to Example~\ref{ex:SMG}, the process
$(Y_u^t)_{0\leq u\leq t}$ is a continuous semimartingale for
every fixed $t$, whenever the driver $Z$ is a Volterra
semimartingale according to Definition~\ref{def:V_semiMG}.
In particular, the functional It\^o formula
of~\cite{dupire2019functional} applies to
$s\mapsto F(s,Y^t_{\cdot\wedge s})$ for sufficiently regular
$F$, with $t$ fixed, and we shall now briefly recall the required
calculus.

Denote by $D([0,T];\RR^k)$ the space of c\`adl\`ag paths,
and write $\omega_{\cdot\wedge s}$ for the path stopped at
time $s$. We denote the space of stopped paths by
$\cS_k:=\{(s,\omega_{\cdot\wedge s}):s\in[0,T],\
\omega\in D([0,T];\RR^k)\}$, equipped with the semi-metric
$$
    d_\infty
    \bigl(
        (s,\omega_{\cdot\wedge s}),
        (r,\eta_{\cdot\wedge r})
    \bigr)
    :=
    |s-r|
    +
    \|\omega_{\cdot\wedge s}-\eta_{\cdot\wedge r}\|_\infty.
$$
A functional $F:[0,T]\times D([0,T];\RR^k)\to\RR$ is called
\emph{non-anticipative} if
$F(s,\omega)=F(s,\omega_{\cdot\wedge s})$ for all
$(s,\omega)$, and can therefore be identified with a map
$F:\cS_k\to\RR$. We assume throughout that all functionals
are progressively measurable with respect to the canonical
filtration. Whenever the following limits exist, 
the horizontal derivative and, for $p=1,\dots,k$, the vertical
derivative in direction $e_p$ are defined by
$$
    \partial_0F(s,\omega)
    :=
    \lim_{h\downarrow0}
    \frac{F(s+h,\omega_{\cdot\wedge s})-F(s,\omega)}{h},
    \qquad
    \partial_pF(s,\omega)
    :=
    \lim_{h\to0}
    \frac{
        F\bigl(
            s,\omega_{\cdot\wedge s}
            +h e_p\mathbf 1_{[s,T]}
        \bigr)
        -F(s,\omega)
    }{h}.
$$
For $p,q\in\NN_0$, we denote by $C_b^{p,q}(\cS_k)$ the
space of functionals for which every ordered composition
of at most $p$ horizontal and at most $q$ vertical derivatives
exists and is bounded and continuous with respect to
$d_\infty$; see, for instance,
\cite[Section~3]{Litterer2011OnAC} for more details.

Returning to our framework, we enhance our alphabet by the time-letter $0$, writing $\AA_0:=\AA\cup\{0\}$. The iterated integrals in Definition~\ref{def:Ito_Volterra signature} are then simply extended by 
$\dd\cZ_r^{(\alpha,\beta),0}:=\dd r$, and thereby extend the definition of $\bZ$ to $\Lambda_{\AA_0}$. The composition of vector-fields in Definition~\ref{def:vector_fields_iteration} is then also naturally extended, where we include now the path-derivatives. First, we extend the definition of $\mathsf{V}$ in~\eqref{eq:generalized_vfields} to $\AA_0$, by setting $\mathsf{V}_0=1$. To extend the composition of vector-fields in~\eqref{eq:v_field_operator_generalized}, we will use for a time letter $i_r=0$ the
convention $l_r=0$, whereas for single and double letters
the indices $l_r$ are chosen as before. Then, for any $\mathbf a=(i_1\cdots i_n,\bsigma)\in\Lambda_{\AA_0}^n$, in direct analogy with~\eqref{eq:v_field_operator_generalized} we have
$$
    \mathbf V^{\mathbf a}(F)(r,\mathbf x,\omega)
    =
    \sum_{l_1}\cdots\sum_{l_n}
    \prod_{q=1}^n
    \left(
        \prod_{\substack{(a,h):\,\bsigma_h(a)=q\\
                         1\leq h\leq |i_a|}}
        \partial_{l_a^h}
    \right)
    \mathsf V_{i_q}^{l_q}(x_q)
    \left(
        \prod_{\substack{(a,h):\,\bsigma_h(a)=n+1\\
                         1\leq h\leq |i_a|}}
        \partial_{l_a^h}
    \right)
    F(r,\omega),
$$ where $\bx=(x_1,\dots,x_n)$ and on the 
terminal functional $F$, $\partial_0$ denotes the horizontal, while $\partial_p$, $p=1,\ldots,k$, denote the vertical
derivatives. With this notation at hand, the fact that $Y\in \cV_\mathbb{F}^\infty$ and the functional Itô formula~\cite[Theorem 9]{Litterer2011OnAC} immediately yields the following change-of-variables in our setting.

\begin{lemma}\label{lem:functional_ito}
Let $Z\in\cV^{\infty-}_{\mathbb F}$ and let $Y$ be a solution
to~\eqref{eq:volterra_eq_2}. For any
$F\in C_b^{1,2}(\cS_k)$ and fixed $t\in[0,T]$, we have
almost surely
$$
    F(t,Y^t_{\cdot\wedge t})
    =
    F(0,\boldsymbol\xi)
    +
    \sum_{a\in\AA_0}
    \int_0^t
    \mathbf V^a(F)
    \bigl(s,Y_s^s,Y^t_{\cdot\wedge s}\bigr)
    \dd\cZ_s^{t,a},
$$
where $\boldsymbol\xi$ denotes the constant path with value
$\xi$.
\end{lemma}

Finally, under sufficient regularity assumptions on $F$, the
functional It\^o formula can be iterated exactly as in
Theorem~\ref{thm:main_sde_expansion}, leading to a finite series
representation. The crucial observation is that the vector-field
factors appearing during the iteration are themselves compatible
with the functional calculus. Indeed, the maps
$(s,\omega)\mapsto V_i^p(\omega_s)$ are non-anticipative
functionals whose horizontal derivative vanishes, while their
vertical derivatives coincide with the usual derivatives of
$V_i^p$. The proof of the following result then follows from exactly the same arguments as the proof of Theorem~\ref{thm:main_sde_expansion}, relying on the Lemma~\ref{lem:functional_ito}.

\begin{theorem}
\label{thm:functional_volterra_expansion}
Let $Z\in\cV^{\infty-}_{\mathbb F}$,
$V_1,\dots,V_m\in C_b^\infty(\RR^k;\RR^k)$, and let $Y$ be a
solution to~\eqref{eq:volterra_eq_2}. Fix $N\in\NN$ and suppose
that $F\in C_b^{N+1,2N+2}(\cS_k)$. Then, for every fixed
$t\in[0,T]$, we have almost surely
$$
    F(t,Y^t_{\cdot\wedge t})
    =
    F(0,\boldsymbol\xi)
    +
    \sum_{n=1}^N
    \sum_{\mathbf a\in\Lambda_{\AA_0}^n}
    \bV^{\mathbf a}(F)(0,\xi,\boldsymbol{\xi})\bZ_{0,t}^{\mathbf a}
    +
    R_{\AA_0}^N(t),
$$
with remainder term given by
$$
R_{\AA_0}^N(t)
    =
    \sum_{\mathbf a\in\Lambda_{\AA_0}^{N+1}}
    \int_{\Delta_{0,t}^{N+1}}
    \bV^{\mathbf a}(F)
    \bigl(
        r_1,
        Y_{r_1}^{r_1},\ldots,Y_{r_1}^{r_{N+1}},
        Y^t_{\cdot\wedge r_1}
    \bigr) \prod_{\ell=1}^{N+1}
    \dd\cZ_{r_\ell}^{r_{\bsigma(\ell)},i_\ell}.
$$
\end{theorem}

\begin{remark}\label{rem:endpoint_functional}
For an endpoint functional $F(t,\omega)=f(\omega_t)$, the
horizontal derivative vanishes and the vertical derivatives agree
with the ordinary derivatives of $f$. Consequently, all terms
containing the time letter vanish and
Theorem~\ref{thm:functional_volterra_expansion} reduces to
Theorem~\ref{thm:main_sde_expansion}.
\end{remark}

\begin{remark}
    \label{rem:V=1}
    Suppose $m=k$ and $V_i\equiv e_i$, where $e_i$
    denotes the $i$th canonical basis vector in
    $\mathbb{R}^m$. In this case, we have
    $Y=\xi+Z$ in~\eqref{eq:volterra_eq_2}, and it
    follows directly from the definition that
    $\mathbf{V}^{(i_1\cdots i_n,\bsigma)}(F)\neq 0$
    only if $\bsigma(\ell)=n+1$ for all
    $1\leq\ell\leq n$.
    Indeed, by definition, we have
    $\mathsf{V}_i^p=\delta_{ip}$ and
    $\mathsf{V}_{[i,j]}^{p,q}
    =\frac{1}{2}\delta_{ip}\delta_{jq}$,
    and the relevant vector-field operators $\bV$
    and signature components in the general expansion
    are given by
    \[
        \bV^{w}(F)(r,\bx,\omega)
        =2^{n-|w|}
        \partial_{i_1}\cdots\partial_{i_n}F(r,\omega),
        \qquad
        \bZ_{0,t}^{(w,\bsigma)}
        =
        \int_{\Delta_{0,t}^{n}}
        \dd\cZ_{t_1}^{t,i_1}\cdots
        \dd\cZ_{t_n}^{t,i_n},
    \]
    where $w=i_1\cdots i_n\in\cW_{\AA_0}^n$
    and $\partial_{[i,j]}:=\partial_{ij}$.
    Note that the latter is simply the classical
    It\^o signature of the augmented path
    $(\cZ^t_u:0\leq u \leq t)$. This is not surprising, since
    for fixed $t$, the process $\cZ^t_\cdot$ is a
    continuous semimartingale, and the same expansion
    can be obtained by iterating the functional
    It\^o formula for $F(t,\xi+Z^t)$, which can be viewed as the
    It\^o counterpart of the functional stochastic
    Taylor expansion in~\cite{Litterer2011OnAC}.
\end{remark}

\subsection{Expected expansions and absolute convergence}\label{sec:expected_expansion}

The goal of this section is to apply the Volterra signature expansion from Theorem~\ref{thm:main_sde_expansion} to expand expectations $\mathbb{E}[f(Y_t)]$ for sufficiently regular $f$. Keeping Theorem~\ref{lem:finite_moments} in mind and having \eqref{eq:taylor_volterra_sde} at hand, one can readily write
\begin{align}\label{eq:expected_expansion}
    \mathbb{E}[f(Y_t)]= f(\xi)+ \sum_{n=1}^N\sum_{\mathbf{a}\in \Lambda_\AA^n}\mathbf{V}^\mathbf{a}(f)(\xi)\boldsymbol{\mu}_{0,t}^\mathbf{a}+\mathbb{E}[R^N_\AA(t)], 
\end{align}
whenever the \emph{expected Itô Volterra signature} $\boldsymbol{\mu}_{0,t}:=(\mathbb{E}[\mathbf{Z}_{0,t}^\mathbf{a}])_{\mathbf{a}\in \Lambda_\AA}$ is componentwise well-defined. 
Specifically, as the driving process $Z$ we will consider
\textit{Gaussian Volterra processes}, defined below as a subclass of $\cV_\FF^{\infty-}([0,T]; \RR^m)$, for which we will prove algebraic relations and a convergence result for the expansion of the expected value in \eqref{eq:expected_expansion}.
To define Gaussian Volterra processes, we partition the alphabet $\cA = \{1,\dots,m\}$ into two disjoint subsets
\begin{align*}
    \cA = \cB \; {\dot\cup} \; \cM,
\end{align*}
where $\cB$ and $\cM$ represent the bounded variation and martingale components, respectively.

\begin{definition} \label{def:brownian_volterra_process}
    Consider the alphabet $\cA = \cB \; {\dot\cup} \; \cM$ from above and let $K_1 \in L^{\infty,1}(\Delta^2; \RR^{|\cB|})$ and $K_2 \in L^{\infty,2}(\Delta^2; \RR^{|\cM| \times d})$ for $d \in \NN$. Further let $B$ be a $d$-dimensional Brownian motion. Then a process $Z: \Delta^2 \to \RR^{m}$ is called a \textit{Gaussian Volterra process} if 
    we have the representation
    \begin{align*}
        Z^{\tau,i}_t 
        =
        \begin{cases}
            \displaystyle \int_0^tK_1^i(\tau,s)\dd s, & i \in \cB, \\
            \displaystyle \int_0^tK^i_2(\tau,s)\dd B_s, & i \in \cM,
        \end{cases}
    \end{align*}
    where we recall the notation $K^i_2(\tau,s)\dd B_s:= \sum_{j=1}^d K^{ij}_2(\tau,s)\dd B_s^j$ for $i \in \cM$.
\end{definition}
We then can split Volterra SDEs~\eqref{eq:volterra_eq_2} into bounded variation and martingale parts as follows
\begin{equation}\label{eq:volterra_eq_3}
    Y_t = \xi + \sum_{i\in \cB} \int_0^t V_i(Y_s) K_1^i(t,s)\dd s+\sum_{i\in \cM} \int_0^t V_i(Y_s)K_2^i(t,s)\dd B_s.
\end{equation}  

\begin{remark}\label{rem:standard}
    A commonly studied class of stochastic Volterra equations is of the form
    \begin{equation}\label{eq:SDE_restricted}
        Y_t
        =
        \xi
        + \int_0^t K_\mu(t,s)\mu(Y_s)\dd s
        + \int_0^t K_\sigma(t,s)\sigma(Y_s)\dd B_s,
    \end{equation}
    where $\xi\in\RR^k$, $B$ is an $\RR^d$-dimensional Brownian motion, $\mu:\RR^k\to\RR^k$ and $\sigma:\RR^k\to\RR^{k\times d}$ are continuous, and $K_\mu,K_\sigma:\Delta^2\to\RR^{k\times k}$ are sufficiently regular kernels.
    Such equations can be cast into the form \eqref{eq:volterra_eq_3} by a componentwise separation of the vector fields from the corresponding Brownian Volterra drivers of Definition~\ref{def:brownian_volterra_process}. We work with \eqref{eq:volterra_eq_3} from the outset, since this separation considerably simplifies notation and subsequent algebraic analysis.
    For instance, the affine Volterra process \cite{jaber2019affinevolterraprocesses}
    \begin{align*}
        Y_t
        &=
        Y_0
        + \int_0^t K(t,s)(b_0+b_1Y_s)\dd s
        + \int_0^t K(t,s)\sqrt{a_0+a_1Y_s}\dd B_s,
    \end{align*}
    with a scalar \emph{Riemann--Liouville kernel}
    \[
        K(t,s)=\frac{(t-s)^{\alpha-1}}{\Gamma(\alpha)},
    \]
    fits directly into this framework by taking one bounded-variation and one martingale component with kernel $K$, together with the corresponding vector fields
    \[
        V_1(y)=b_0+b_1y,
        \qquad
        V_2(y)=\sqrt{a_0+a_1y}.
    \]
    More generally, polynomial Volterra processes will be treated in detail in Subsection~\ref{sec:fPVP}.
\end{remark}

\subsubsection{Expected Itô Volterra signature for Gaussian Volterra processes}
By Theorem \ref{lem:finite_moments}, all local martingales in the expansion~\eqref{eq:taylor_volterra_sde} are true martingales and
the respective expected Itô Volterra signature terms 
vanish. This motivates the definition of the restricted alphabet
\begin{align*}
    \AA_E := \cB \cup [\cM,\cM]= \cB \cup \{ [i,j]:i,j \in \cM\} \subseteq \AA,
\end{align*}
where the letter $E$ again refers to expectation. Consequently, define the sets $\cW_{\AA_E}^n$, $\cW_{\AA_E}$, $\Lambda_{\AA_E}^n$ and $\Lambda_{\AA_E}$ for $n \in \NN$ analogously to the discussion following Lemma \ref{lem:algebra}.
The martingale term cancellations are then directly expressed as follows.

\begin{lemma} \label{lem:expected_vsig}
    Consider a Gaussian Volterra process $Z$ from Definition \ref{def:brownian_volterra_process} and consider its Itô-Volterra signature $\bZ$ from Definition \ref{def:Ito_Volterra signature} with componentwise expectation $\bmu$. Then for $(s,t) \in \Delta^2$ we have
    \begin{align*}
        \bmu_{s,t}^{\boldsymbol{a}} =
        \begin{cases}
            0, \qquad & \boldsymbol{a} \in \Lambda_{\AA} \setminus \Lambda_{\AA_E} \\
            \bZ^{\boldsymbol{a}}_{s,t}, \qquad & \boldsymbol{a} \in \Lambda_{\AA_E}
        \end{cases}
    \end{align*}
    In particular $\bZ_{s,t}^{\boldsymbol{a} }$ is deterministic for all $\boldsymbol{a} \in \Lambda_{\AA_E}$ and of the form
    \begin{align*}
        \bZ^{\boldsymbol{a}}_{s,t}
        =
        \int_{\Delta_{s,t}^n}
        \prod_{\ell=1}^n K^{i_\ell}(r_{\bsigma(\ell)}, r_\ell) \dd r_\ell,
    \end{align*}
    where we use the convention $r_{n+1}:=t$, and 
    \begin{align*}
        K^{i_\ell}(r_{\bsigma(\ell)}, r_\ell)
        \assign
        \begin{cases}
            \displaystyle K_1^{i}(r_{\bsigma_1(\ell)}, r_\ell), & \text{ if } i_\ell = i, \\
            \displaystyle  (K_2(r_{\bsigma_1(\ell)}, r_\ell)K_2(r_{\bsigma_2(\ell)}, r_\ell)^\top)^{ij}, & \text{ if } i_\ell= [i,j].
        \end{cases}
    \end{align*}
\end{lemma}
\begin{proof}
    Note that whenever $\boldsymbol{a} = (w,\bsigma)\in \Lambda_{\AA} \setminus \Lambda_{\AA_E}$ with $w = i_1\cdots i_n$ there must at least be one letter $i_\ell$ with $i_\ell = i$ and $i \in \cM$ or $i_\ell = [i,j]$ with at least one letter in $\cB$. By Definition \ref{def:Ito_Volterra signature} and with Theorem \ref{lem:finite_moments} we must have $\bmu^{\boldsymbol{a}}_{s,t} = 0$. Otherwise, for $\boldsymbol{a} = (w,\bsigma)\in \Lambda_{\AA_E}$ all components in Definition \ref{def:Ito_Volterra signature} are deterministic since the quadratic variation of the martingale components are deterministic in the Gaussian Volterra process case. This finishes the proof.
\end{proof}

The quasi-shuffle identity in Proposition~\ref{prop:quasi_shuffle} describes the multiplication of Itô Volterra signature coordinates.
Upon taking expectations, terms containing single martingale letters vanish, leaving only bounded-variation letters and deterministic quadratic-variation brackets.
Consequently, the expected Itô Volterra signature of Gaussian Volterra processes is supported on $\Lambda_{\AA_E}$, while retaining the quasi-shuffle structure in its kernel assignments.
This leads to the following \emph{expected shuffle} identity, which combines the ordinary shuffle $\shuffle$ of words with the quasi-shuffle $\star$ of their kernel assignments.

\begin{proposition} \label{prop:expected_shuffle}
    Let $\boldsymbol{a}=(w,\bsigma)\in\Lambda_{\AA}^n$ and $\boldsymbol{b}=(v,\btau)\in\Lambda_{\AA}^m$. Recall the notion $\bsigma \; \hat \star_\eta \;\btau$ for any $\eta \in \mathrm{Sh}(n,m) \subset \mathrm{QSh}(n,m)$ from Definition \ref{def:Sigma_quasi_shuffle} and the $\eta$-shuffle from Definition \ref{def:Sigma_shuffle}. Define the \emph{expected shuffle} as
    \begin{align*}
        \boldsymbol{a} \shuffle_E \boldsymbol{b}
        = (w,\bsigma) \shuffle_E (v,\btau) \assign \sum_{\eta \in \mathrm{Sh}(n,m)} (w \shuffle_\eta v, \bsigma \; \hat \star_\eta \;\btau).
    \end{align*}
    Then the expected signature satisfies
    \begin{align*}
        \bmu^{\boldsymbol{a}}_{s,t} \bmu^{\boldsymbol{b}}_{s,t} 
        = 
        \bmu^{\boldsymbol{a} \shuffle_E \boldsymbol{b}}_{s,t}.
    \end{align*}
    In particular, outside the set $\AA_E$ both sides are equal to zero. 
\end{proposition}
\begin{proof}
    Suppose first that w.l.o.g.\ $\boldsymbol{a} \in \Lambda_\AA \setminus \Lambda_{\AA_E}$. In this case the left hand side is equal to zero since $\bmu^{\boldsymbol{a}}_{s,t} = 0$. By assumption $w$ contains at least one bracket letter of the form $[i,j]$ with $i,j$ not both in $\cM$ or a single letter $k \in \cM$. Any shuffle $\eta \in \mathrm{Sh}(n,m)$ might only change the position of this letter. Hence, each $w \shuffle_\eta v$ also contains this letter and makes the left hand side equal to zero.

    Assume that $\boldsymbol{a}, \boldsymbol{b} \in \Lambda_{\AA_E}$. The important part to recall is that $\bmu = \bZ$ on $\Lambda_{\AA_E}$. Hence we have by Proposition \ref{prop:quasi_shuffle}
    \begin{align*}
        \bmu^{\boldsymbol{a}}_{s,t} \bmu^{\boldsymbol{b}}_{s,t} 
        =
        \bZ^{\boldsymbol{a}}_{s,t} \bZ^{\boldsymbol{b}}_{s,t}
        = 
        \bZ^{\boldsymbol{a} \star \boldsymbol{b}}_{s,t},
    \end{align*}
    and it remains to show $\bZ^{\boldsymbol{a} \star \boldsymbol{b}}_{s,t} = \bmu^{\boldsymbol{a} \shuffle_E \boldsymbol{b}}_{s,t}$. However, as one can easily see $w \star_\eta v = 0$ for any $\eta \in \mathrm{QSh}(n,m) \setminus \mathrm{Sh}(n,m)$ since the letters of $w$ and $v$ are in $\AA_E$. Hence we are left with the definition of $\shuffle_E$ and the proof is finished.
\end{proof}

\begin{example}
    Consider the two (expected) Itô Volterra signature terms (of a suitable one-dimensional Gaussian Volterra process $Z$)
    \begin{align*}
        \int_{\Delta_{0,t}^1} K(t,r_1) \dd r_1, \qquad \int_{\Delta_{0,t}^2} K(t,r_2) K(r_2,r_1) K(t,r_1) \dd (r_1,r_2),
    \end{align*}
    corresponding to $\boldsymbol{a} = (w, \bsigma)$ with $w = 1$ and $\bsigma(1) = (2,2)$, and $\boldsymbol{b} = (v, \btau)$ with $v = [2,2]1$ and $\btau(1) = (2,3)$ and $\btau(2) = (3,3)$, respectively. As can be seen from Definition \ref{def:volterra-words}, the concatenation $(\sigma_1\tau_1, \sigma_2\tau_2)$ takes values $(4,4), (3,4), (4,4)$ for inputs $e = 1, 2, 3$, respectively. 
    Noting that there are three admissible shuffles in $\mathrm{Sh}(1,2)$, the term $\bmu^{\boldsymbol{a}\shuffle_E \boldsymbol{b}}_{0,t}$ evaluates to
    \begin{align*}
        \int_{\Delta_{0,t}^3} K(t,r_3) K(t,r_2) K(r_3,r_1) K(t,r_1) \dd r 
        &+
        \int_{\Delta_{0,t}^3} K(t,r_3) K(r_3,r_2) K(t,r_2) K(t,r_1) \dd r\\
        &+
        \int_{\Delta_{0,t}^3} K(t,r_3) K(t,r_2) K(r_2,r_1) K(t,r_1) \dd r,
    \end{align*}
    with $\dd r = \dd (r_1,r_2,r_3)$, and as some calculations show, this sum precisely equals $\bmu^{\boldsymbol{a}}_{0,t}\bmu^{ \boldsymbol{b}}_{0,t}$.
\end{example}

\begin{remark}\label{rem:fawcett}
    Suppose $\cB=\varnothing$ and set $M_r=Z_r^T=\EE[Z_T^T\mid\cF_r]$. Write $\mathrm{Sig}(M)_{s,T}\in T((\RR^{\cM}))$ for its Stratonovich signature. For $w=i_1\cdots i_{2n}\in\cW_{\cM}^{2n}$, let
    $\mathfrak b(w)=([i_1,i_2]\cdots[i_{2n-1},i_{2n}],\bsigma)$ with $\sigma_1(k)=\sigma_2(k)=n+1$.
    Since the covariations of $M$ are deterministic, only complete adjacent-pair contractions survive expectation. Each pair $i,j$ contributes $\tfrac12\,\dd[M^i,M^j]$ by Itô--Stratonovich conversion, so
    \begin{equation}\label{eq:frozen_expected_signature}
        \EE[\mathrm{Sig}(M)_{s,T}^{w}]
        =2^{-n}\int_{\Delta_{s,T}^n}\prod_{k=1}^n\dd[M^{i_{2k-1}},M^{i_{2k}}]_{r_k}
        =2^{-n}\bmu_{s,T}^{\mathfrak b(w)}.
    \end{equation}
    Odd levels have zero expectation. The logarithm of this expected signature is a Magnus series over symmetric two-tensors~\cite[Section~6.1.1, Example~6.1]{friz2022unified}. On the indices $\mathfrak b(w)$, Proposition~\ref{prop:expected_shuffle} is the shuffle character identity of its group-like exponential. For constant $K_2$, the evaluation maps drop out and Fawcett's formula~\cite{Fawcett2003} reads $\bmu_{s,T}=\exp_\otimes\!\bigl((T-s)\sum_{i,j\in\cM}(K_2K_2^\top)^{ij}[i,j]\bigr)$ as a series over $\AA_E$.
\end{remark}

\subsubsection{A convergence result for the expected expansion}

In the following, let $Z$ be a Gaussian Volterra process, let $V = (V_1,\dots,V_m) \in C(\RR^k;\RR^k)^m$ be continuous vector fields and let $\xi \in \RR^k$ be the initial condition. We consider the stochastic Volterra equation
\begin{align} \label{eq:volterra_eq_4}
    Y_t = \xi + \sum_{i=1}^m \int_0^t V_i(Y_s) \dd Z^{t,i}_s
    = \xi + \sum_{i\in \cB} \int_0^t V_i(Y_s) \dd Z^{t,i}_s
   + \sum_{i\in \cM} \int_0^t V_i(Y_s) \dd Z^{t,i}_s
\end{align}
and assume that the solution $Y$ exists. 
The regularity assumptions on the vector fields remain minimal for now, but will be strengthened later on. Any such solution process $Y$ can arise as a strong or weak solution of the underlying Volterra SDE. In the former case, we are given a probability space and a Brownian motion and find a corresponding solution $Y$. In the latter case, the space and Brownian motion are part of the solution. In either case we are left with a process $Y$ satisfying \eqref{eq:volterra_eq_4}.

The upcoming theorem establishes the expected series expansion for Volterra SDEs of the form \eqref{eq:volterra_eq_4}. Recall once more that $\bmu^{\boldsymbol{a}} = 0$ for $\boldsymbol{a} \in \Lambda_\AA \setminus \Lambda_{\AA_E}$ and $\bmu^{\boldsymbol{a}} = \bZ^{\boldsymbol{a}}$ for $\boldsymbol{a} \in \Lambda_{\AA_E}$, where we use the fact that $\cB$-components and the quadratic variation of $\cM$-components of the Volterra path $Z$ are deterministic.
This has the implication that we only require the enhanced vector field from Definition \ref{def:vector_fields_iteration} to be smooth for indices in $\AA_E$ and not on the full alphabet $\AA$. Importantly, this includes the choice $\sigma(y) = \sqrt{\psi(y)}$ for $\psi$ sufficiently smooth in equation \eqref{eq:SDE_restricted}. As a consequence, the following result holds in particular  for affine and polynomial Volterra equations \cite{jaber2019affinevolterraprocesses} and \cite{abi2024polynomial}. 
Concretely, the  theorem provides an explicit series representation of $\EE[f(Y_t)]$ and establishes absolute convergence of the series under a set of assumptions on the input $(f,V,K)$. These conditions have to be verified on a case-by-case basis. While they appear rather abstract at this level of generality, we verify them in the next section for polynomial specifications.

\begin{Assumption}\label{ass:absolut_convergence_assumption}
Let $\omega:\RR^k \rightarrow [1,\infty)$ be a continuous weight function. For $n \in \NN$ and any function $P: (\RR^k)^n \to \RR$ we define 
\begin{align*}
    \|P\|_\omega 
    := 
    \sup_{x_1,\dots, x_n \in \RR^k}
    \frac{|P(x_1,\dots,x_n)|}{\max_{j \leq n} \omega(x_j)},
\end{align*}
Let $f \in C^\infty(\RR^k, \RR)$ and consider the following set of assumptions:
\begin{itemize}
    \item[(i)]\label{ass:one} There exist positive constants $C,A>0$ and $\rho \geq 0$ such that for $n \in \NN_0$ we have
    \begin{equation}\label{eq:vf_assumption}
         \sum_{\mathbf{a}\in  \Lambda_{\AA_E}^n}\Vert \bV^\mathbf{a}(f)\Vert_{\omega}\leq CA^n(n!)^\rho,
    \end{equation}
    where $\bV^{\boldsymbol{a}}(f)(x) = f(x)$ for $\mathbf{a}\in \Lambda_{\AA}^0$.
    \item[(ii)]\label{ass:two} For any $n\in \mathbb{N}$ we have \begin{equation}
        \label{eq:MG_assumption}
\Vert \bV^{(iv,\bsigma)}(f)\Vert_{\omega}<\infty, \qquad \forall (iv,\bsigma)\in {\Lambda}_{\AA}^n \text{ where }
i\in\cM,\text{ and } v\in \cW^{n-1}_{\AA_E}.
    \end{equation}
    \item[(iii)]\label{ass:three}
    It holds that $\sup_{0\leq t \leq \tau \leq T}\mathbb{E}[\omega(Y_t^\tau)^2]<\infty$. 
    \item [(iv)]\label{ass:four} There exists $\eta>\rho$ and some constant $\widetilde{C}=\widetilde{C}_{\eta,T}>0$ such that \[
    |\bmu_{0,T}^\mathbf{a}| \leq \widetilde{C}^n\frac{1}{(n!)^\eta}, 
    \]
    for all $\mathbf{a} \in {\Lambda}_{\AA_E}^n$ with $n \in \NN$ such that $\mathbf{V}^\mathbf{a}(f)\neq 0$.
\end{itemize}
\end{Assumption}

\begin{theorem}\label{thm:expansion_expected}
    Assume that the enhanced vector-field $(\mathsf{V}_{\boldsymbol{a}})_{\boldsymbol{a} \in \AA_E}$ from Definition \ref{def:vector_fields_iteration} is smooth. 
    Suppose that Assumption \ref{ass:absolut_convergence_assumption} (i),(ii) and (iii) hold for smooth $f\in C^\infty(\RR^k;\RR)$.
    Then we have
    \begin{equation}\label{eq:exp_vsig_expansion}
        \EE[f(Y_t)] = f(\xi) + \sum_{n=1}^N\sum_{\mathbf{a}\in \Lambda_{\AA_E}^n} \bV ^{\mathbf{a}}(f)(\xi){\boldsymbol{\mu}}^{\mathbf{a}}_{0,t} + R^N_{\AA_E}(t),
    \end{equation} where the remainder is given by
    \[
        {R}^{N}_{\AA_E}(t)
        = 
        \sum_{(w,\bsigma)\in \Lambda_{\AA_E}^{N+1}} \int_{\Delta_{0,t}^{N+1}}
        \EE[\bV^{(w,\bsigma)}(f)(Y_{r_1}^{r_1},\dots,Y_{r_1}^{t})] 
        \prod_{l=1}^{N+1}\dd \cZ_{r_l}^{r_{\bsigma(l)},i_l}, \quad r_{N+2}:=t.
    \]
    Moreover, if Assumption \ref{ass:absolut_convergence_assumption} (iv) is additionally in force, the resulting series is absolutely convergent and we have
    \begin{equation}\label{eq:full_exp_vsig_expansion}
        \EE[f(Y_t)] = f(\xi) + \sum_{n\geq 1}\sum_{\mathbf{a}\in \Lambda_{\AA_E}^n} \bV ^{\mathbf{a}}(f)(\xi)\bmu^{\mathbf{a}}_{0,t}.
    \end{equation}
\end{theorem}

\begin{proof}
    Let $n \in \NN_0$ and $0 \leq t \leq \tau_1 \leq \dots \leq \tau_{n+1} \leq T$. Then we have by (i), (ii) and (iii)
    \begin{align} \label{eq:exp_exp_1}
        \EE[|\bV^{\boldsymbol{a}}(f)(Y_t^{\tau_1}, \dots ,Y_t^{\tau_{n+1}})|^2] 
        \leq
        (n+1)
        \|\bV^{\boldsymbol{a}}(f)\|_\omega^2
        \sup_{0 \leq t \leq \tau \leq T}\EE[\omega(Y_t^{\tau})^2] < \infty.
    \end{align}
    for all $\boldsymbol{a} \in \Lambda_{\AA_E}^n$ and $\boldsymbol{a} = (iv,\bsigma) \in \Lambda^n_\AA$ where $i \in \cM$ and $v \in \cW^{n-1}_{\AA_E}$.
    This fact will be used to show true martingality later on and we continue to show the result by induction over $N$ with an application of Proposition \ref{prop:vector_field_algebra}.

    The case $N = 0$ can be directly seen from an application of Itô's formula similar to the beginning of Section \ref{sec:stochastic_expansion}. Noting that the quadratic variation of $\cB$- and $\cM$-components of $Z$ vanishes, we have
    \begin{equation*}
        f(Y_t)=f(\xi) + \sum_{a\in \cM}\int_0^t\mathbf{V}^{a}(f)(Y_s^s,Y_s^t)\dd \cZ_s^{t,a} +\sum_{a\in \mathbb{A}_E}\int_0^t\mathbf{V}^{a}(f)(Y_s^s,Y_s^t)\dd \cZ_s^{t,a}.
    \end{equation*}
    Using \eqref{eq:exp_exp_1} and that $K_2 \in L^{\infty,2}([0,T];\RR^{|\cM| \times d})$ shows that the expectation of the sum over $\cM$ vanishes and the case $N=0$ is shown.
    
    For the general case we apply Itô's formula similar to the proof of Theorem \ref{thm:main_sde_expansion} to each integrand in the remainder term $R^{N-1}_{\AA_E}$. Let $(w,\bsigma) \in \Lambda_{\AA_E}^N$ and set $r_{N+1} := t$.  
    Applying Itô's formula on $[0,r_1]$
    to $\bV^{(w,\bsigma)}(f)(Y_{r_1}^{r_1},\dots,Y_{r_1}^t)$ yields 
    \begin{equation*}
        \begin{aligned}
            &\bV^{(w,\bsigma)}(f)
            \bigl(
                Y_{r_1}^{r_1},\dots,Y_{r_1}^{t}
            \bigr)
            \\
            &\quad =
            \bV^{(w,\bsigma)}(f)(\xi)
            \\
            &\qquad+
            \sum_{e=1}^{N+1}
            \sum_{p=1}^k
            \sum_{i=1}^m
            \int_0^{r_1}
            \mathsf{V}_i^p(Y_u^u)
            \bigl(
                \partial_{x_e^p}\bV^{(w,\bsigma)}(f)
            \bigr)
            \bigl(
                Y_u^{r_1},\dots,Y_u^t
            \bigr)
            \,\dd Z_u^{r_e,i}
            \\
            &\qquad+
            \sum_{e,h=1}^{N+1}
            \sum_{p,q=1}^k
            \sum_{i,j=1}^m
            \int_0^{r_1}
            \mathsf V_{[i,j]}^{p,q}(Y_u^u)
            \bigl(
                \partial_{x_e^p}\partial_{x_h^q}
                \bV^{(w,\bsigma)}(f)
            \bigr)
            \bigl(
                Y_u^{r_1},\dots,Y_u^t
            \bigr)
            \,\dd\cZ_u^{r_e,r_h,[i,j]}.
        \end{aligned}
    \end{equation*}
    Note again that $\ZZ^{r_e,r_h,[i,j]}$ represents the quadratic variation of $Z^{r_e,i}$ and $Z^{r_h,j}$ and that the factor of $1/2$ from Itô's formula is included in $\mathsf{V}_{[i,j]}^{p,q}$. Using Proposition \ref{prop:vector_field_algebra} this expression can be rewritten precisely as
    \begin{equation} \label{eq:exp_exp_2}
    \begin{aligned}
        &\bV^{(w,\bsigma)}(f)
        \bigl(
            Y_{r_1}^{r_1},\dots,Y_{r_1}^{t}
        \bigr) \\
        &\hspace{50pt}=
        \bV^{(w,\bsigma)}(f)(\xi) +
        \sum_{\ell \in \AA} \sum_{e \in [N+1]^{|\ell|}} \int_0^{r_1} \bV^{(\ell w,\bsigma_+^e)}(f)
        \bigl(
            Y_{u}^{u},\dots,Y_{u}^{t} 
        \bigr)
        \dd \ZZ^{r_{e}, \ell}_u
    \end{aligned}
    \end{equation}
    where $r_e$ is understood to be $(r_{e_1},r_{e_2})$ for $|e|=2$.
    Upon taking expectations, the quadratic variation of $\cB$- and $\cM$-components vanish and the $\cM$-components are true martingales by \eqref{eq:exp_exp_1}.
    Plugging this into the remainder $R_{\AA_E}^{N-1}$ we are thus left with
    \begin{align*}
        R_{\AA_E}^{N-1}(t) = \sum_{\boldsymbol{a} \in \Lambda_{\AA_E}^N}  \bV^{\boldsymbol{a}}(f)(\xi) \bmu^{\boldsymbol{a}}_{0,t} + \widetilde R_N(t),
    \end{align*}
    where, upon relabeling, 
    \begin{align*}
        \widetilde R_N(t)
        &=
        \sum_{(i_2\cdots i_{N+1},\bsigma)\in\Lambda_{\AA_E}^N}
        \sum_{i_1\in\AA_E}
        \sum_{e\in[N+1]^{|i_1|}}
        \int_{\Delta_{0,t}^{N+1}}
        \EE[\bV^{(i_1\cdots i_{N+1},\bsigma_+^e)}(f)
        \bigl(
            Y_{r_1}^{r_1},\dots,Y_{r_1}^{t}
        \bigr)]
        \prod_{l=1}^{N+1}
        \dd\cZ_{r_l}^{r_{\bsigma_+^e(l)},i_l},
    \end{align*}
    with $r_{N+2} := t$. By \eqref{eq:sets_recursion}, each
    $(w,\btau)\in\Lambda_{\AA_E}^{N+1}$ is obtained uniquely by
    choosing its first letter $i_1$, a compatible tail
    $(i_2\cdots i_{N+1},\bsigma)\in\Lambda_{\AA_E}^N$,
    and $e\in[N+1]^{|i_1|}$ such that
    $\btau=\bsigma_+^e$.
    Re-indexing the last expression therefore gives
    $\widetilde R_N(t)=R^N_{\AA_E}(t)$,
    which completes the induction. 
    
    Assume now that Assumption \ref{ass:absolut_convergence_assumption} (iv) is in force.
    Note that for any $n \geq 1$ we have
    \begin{align*}
        \sum_{\boldsymbol{a} \in \Lambda_{\AA_E}^n}|\bV^{\boldsymbol{a}}(f)(\xi)| \leq \sum_{\boldsymbol{a} \in \Lambda_{\AA_E}^n}\omega(\xi) \|\bV^{\boldsymbol{a}}(f)\|_\omega \leq \omega(\xi)C A^n (n!)^\rho.
    \end{align*}
    Combining this with the assumed uniform bound on the expected Volterra signature and $\eta > \rho$ gives
    \begin{align*}
        \sum_{\boldsymbol{a} \in \Lambda_{\AA_E}^n} |\bV ^{\mathbf{a}}(f)(\xi)\bmu^{\mathbf{a}}_{0,t}| \leq
        \widetilde{C}\frac{1}{(n!)^\eta} \omega(\xi)C A^n (n!)^\rho = \omega(\xi)C\widetilde{C} A^n \frac1{(n!)^{\eta - \rho}}
    \end{align*}
    establishing absolute convergence. This finishes the proof.
\end{proof}

\begin{remark}
    It is worth noting that, under a slight modification of Assumption~\ref{ass:absolut_convergence_assumption}, the same absolute convergence holds for the CDE expansion presented in Section~\ref{sec:cde_expansion}, see   Theorem~\ref{thm:main_cde_expansion}. Indeed, replacing $|\bmu^\mathbf{a}|$ in the growth condition (iv) by $|\bZ^\mathbf{a}|$, a similar argument leads to absolute convergence of the series~\eqref{eq:taylor_volterra_cde}. In particular, choosing $\omega\equiv 1$, in the  deterministic setting (ii) and (iii) always holds true, and as we shall see in the next section, assumption (iv) includes natural classes of  integrable kernels $K$.
\end{remark}

\subsection{Case study: fractional polynomial processes}\label{sec:fPVP}
We close this section with a concrete case study for the  class of polynomial  Volterra SDEs, for which one can guarantee absolute convergence for the expected series expansions~\eqref{eq:full_exp_vsig_expansion}. By means of Theorem~\ref{thm:expansion_expected}, this reduces to the verification of Assumption~\ref{ass:absolut_convergence_assumption}. The first condition (i) is essentially a regularity requirement for the input $(f,V)$, preventing the differential operators $\mathbf{V}(f)$ from growing to quickly. Similar conditions (for $\omega\equiv 1$) are standard already in the early literature of series representation when $K\equiv 1$, but for the coefficients $|\bV^{(i_1)}\circ \cdots \circ \bV^{(i_n)}|$, see Proposition~\ref{prop:vector_field_algebra}. Classical examples are \cite[Proposition II.2]{fliess1981fonctionnelles} in the deterministic,  or \cite[Corollary 1]{arous1989flots} in the stochastic setting. The second conditions (ii)-(iii) are technical requirements ensuring true martingality of the stochastic integrals in the expansions, so that they vanish under expectation. On  the other-hand, the condition (iv) concerns the decay of the expected Volterra signature, which can be guaranteed already with simple integrability assumptions for the kernels, which we illustrate in the following lemma.
\begin{lemma}\label{lem:kernel_condition}
    Fix $\rho \geq 0$. Suppose that there exists some $q > 2$ with $q-2 > \rho q$ such that $K_1$ and $K_2$ from Definition \ref{def:brownian_volterra_process} satisfy $K_1 \in L^{\infty, q/2}(\Delta^2;\RR^m)$
    and $K_2 \in L^{\infty, q}(\Delta^2;\RR^{m \times d})$. Then Assumption~\ref{ass:absolut_convergence_assumption} (iv) holds with $\eta= \frac{q-2}{q}.$
\end{lemma}
\begin{proof}
    Note that we are in the Brownian Volterra setting from Definition \ref{def:brownian_volterra_process} and this translates to 
    \begin{align*}
        \ZZ_t^{\tau_1,\tau_2,i_k} 
        =
        \begin{cases}
            \int_0^t K_1^i(\tau_1,s) \dd s, \quad & i_k = i \in \{1,\dots,m\} \\
            \int_0^t (K_2(\tau_1,s) K_2(\tau_2,s)^\top)_{ij}  \dd s, \quad & i_k = [i,j] \text{ with } i,j \in \{m,\dots,2m\}
        \end{cases}
    \end{align*}
    in Definition \ref{def:Ito_Volterra signature}. 
    For simplicity, by $|K(r_{\bsigma(l)}, r_l)|$ for $l \in \{1,\dots,n\}$ we mean $|K_1(r_{\bsigma_1(l)}, r_l)|$ if $|i_l| = 1$ and $|K_2(r_{\bsigma_1(l)}, r_l)||K_2(r_{\bsigma_2(l)}, r_l)|$ if $|i_l| = 2$. 
    Denote the constant
    \begin{align*}
        c := \max\big\{\|K_1\|_{L^{\infty, q / 2}}^{q / 2}, \|K_2\|_{L^{\infty,q}}^{q}\big\},
    \end{align*}
    then one can proceed by induction to show that
    \begin{align*}
        \int_{\Delta^n_{0,t}} \prod_{l=1}^n |K(r_{\bsigma(l)}, r_l)|^{\frac q2} \dd r_l \leq c^n.
    \end{align*}
    Thus, applying Hölder's inequality we have 
    \begin{align*}
        |\boldsymbol{\mu}_{0,t}^{\boldsymbol{a}}|
        \leq
        \mathrm{Vol}(\Delta_{0,t}^n)^{\frac{q-2}{q}} c^n
        \leq
        (n!)^{\frac{2-q}{q}} \big( c t^{\frac{q-2}{q}}\big)^n.
    \end{align*}
    Since $\rho q < q-2$ by assumption this finishes the proof.
\end{proof}
In the following we consider fractional polynomial processes, that is for any $\alpha\in(1/2,1)$, we consider weak solutions (see \cite{abi2024polynomial}) to
\begin{equation}\label{eq:fractional_PVP}
Y_t= \xi+ \int_0^tb(Y_s)(t-s)^{\alpha-1}\dd s + \int_0^t \sigma(Y_s)  (t-s)^{\alpha-1}  \dd B_s, \quad 0\leq t \leq T,
\end{equation}
where $\xi\in \mathbb{R}^k$, $B$ is an $m$-dimensional Brownian motion, and $b:\mathbb{R}^k\rightarrow \mathbb{R}^k$, $\sigma:\mathbb{R}^k \rightarrow \mathbb{R}^{k\times m}$ are such that \begin{equation}\label{eq:poly_spec}
    b\in \mathrm{Pol}_1, \qquad a:=\sigma\sigma^\top \in \mathrm{Pol}_2.
\end{equation}  Weak solutions to \eqref{eq:fractional_PVP} belong to the class of polynomial Volterra processes introduced in \cite{abi2024polynomial}. In the following result, we show that for any polynomial $f\in\mathrm{Pol}_p(\mathbb{R})$ of degree at most $p\in\mathbb{N}$, the expansion \eqref{eq:exp_vsig_expansion} converges absolutely, and we provide a precise decay rate for the remainder as $N\to\infty$.

\begin{proposition}
Let $Y$ be a weak-solution to \eqref{eq:fractional_PVP} for some $\alpha \in (1/2,1)$, and fix $f\in \mathrm{Pol}_p(\mathbb{R}^k)$ for some $p\in \mathbb{N}$. Then (i)-(iv) from Assumption~\ref{ass:absolut_convergence_assumption} hold true, and thus $\mathbb{E}[f(Y_T)]$ admits an absolutely convergent series representation \eqref{eq:full_exp_vsig_expansion}. Moreover, we have the following remainder estimate \begin{equation}\label{eq:remainder_PVP}
    \sup_{0\leq t \leq T}|R^N_{\AA_E}(t)| \leq C_{f,V,T,\xi} \frac{D^{N+1}}{(N+1)!^{2\alpha-1}},
\end{equation} where $D>0$ is a constant independent of $N$.
\end{proposition}

\begin{proof}
    Let us choose the weight-function $\omega(x):=(1+|x|)^p$ in Assumption~\ref{ass:absolut_convergence_assumption}. Condition (iii) then follows directly from \cite[Lemma 3.1]{abi2024polynomial}. On the other-hand, condition (iv) holds since the kernel $K(t,s)=(t-s)^{\alpha-1}$ is in $L^{\infty,q}$ for all $1\leq q \leq \frac{1}{1-\alpha}$. In particular, by Lemma~\ref{lem:kernel_condition}, condition (iii) holds for $\eta = 2\alpha-1 >0$. Finally, we will now show that (i)-(ii) hold with $\rho=0$, such that from the proof of Theorem~\ref{thm:expansion_expected} it follows that \[
    \sup_{0\leq t \leq T}|R^N_{\AA_E}(t)| \leq C \frac{D^{N+1}}{(N+1)!^{2\alpha-1}}
    \] for two constants $C,A>0$, only depending on $f,V,T,\xi$ but not on $N$. In this case, the series~\eqref{eq:full_exp_vsig_expansion} is indeed absolutely convergent by the ratio test. Thus we are left proving Assumption~\ref{ass:absolut_convergence_assumption} (i)-(ii). For this we observe that $(w,\bsigma) \in \Lambda^n_{\AA_E}$, the map $(x_1,\dots,x_{n+1})\mapsto \bV^{(w,\bsigma)}(f)(x_1,\dots,x_n)$ is a polynomial of degree at most $p$. Indeed, recalling that $\AA_E = \cB \cup  [\cM,\cM]$, by Definition~\ref{def:vector_fields_iteration} we have $\mathsf{V}_i \in \mathrm{Pol}_1$ if $i \in \cB$, and $\mathsf{V}_i \in \mathrm{Pol}_2$ otherwise, so that $\bV^\mathbf{a}(f)$ is a polynomial, and by~\eqref{eq:v_field_operator_generalized} its degree is at most the degree of $f$. Moreover, recalling the recursive identity~\eqref{eq:recursion_formula} from Proposition~\ref{prop:vector_field_algebra}, for any $v\in \cW^{n-1}_{\AA_E}$ and $i \in \AA_E$ and $\bsigma^e_+$ with $\bsigma \in (\Sigma^{n-1})^2$ and $e \in [N+1]^{|i|}$, we have  \begin{equation}\label{eq:VF_one_step_proof}\bV^{(iv,\bsigma_+^e)}(f)(x_1,\dots,x_{n+1}) = \begin{cases}
        \sum_{l_1 \in [k]} \mathsf{V}^{l_1}_i(x_1) \partial_{x_{l_1}^e} \bV^{(w,\bsigma)}(x_2,\dots,x_{n+1}) & |i|=1 \\
     \sum_{l_1 \in [k]^{2}}\mathsf{V}_{[i_1,i_2]}^{l_1}(x_1)\partial_{x_e^{l_1^1}x_h^{l_1^2}}\mathbf{V}^{(w,\boldsymbol{\sigma})}(f)(x_2,\dots,x_{n+1}) & |i|=2.
    \end{cases} \end{equation}
    Recall from the same proof that $\btau \in (\Sigma^n)^2$ is of the form $\btau = \bsigma_+^e$ for some $\bsigma \in (\Sigma^{n-1})^2$. As a consequence, if $i\in \cM$ as in Assumption~\ref{ass:absolut_convergence_assumption} (ii), we know $(\mathsf{V}^{l_1}_i)^2 \in \mathrm{Pol}_2$ and $\partial_{x_e^{l_1}}\bV^{(w,\bsigma)}(f) \in \mathrm{Pol}_{p-1}$, and thus $$|V_i(x)| \leq C(1+|x|), \qquad |\partial_{x_e^{l_1}}\bV^{(w,\bsigma)}(x_2,\dots,x_{n+1})| \leq K\max_{2\leq u \leq n+1}(1+|x_i|)^{p-1},$$ so that (ii) holds. To show (i), for general polynomials $P(x)=\sum_{\alpha \in \cI}c_\alpha x^\alpha$, we introduce the coefficient norm $\Vert P \Vert_{c}:= \sum_{\alpha \in \cI}|c_\alpha|$. Combined with the weight-function $\omega(x)=(1+|x|)^p$, for any polynomial $P$ of degree $p$ we have $\Vert P \Vert_{\omega} \leq C\Vert P \Vert_c$. Since $\bV(f)$ are polynomials restricted to $\Lambda_{\AA_E}$, it suffices to show $(i)$ with respect to the coefficient norm. Using once again that \[
{\Lambda}_\AA^{n+1} = \bigcup_{\substack{\mathbf{e} \in \mathbb{A}_{n+1},\, i\in \mathbb{A} \\ |\mathbf{e}|=|i|}} \Big \{ (iw,\boldsymbol{\sigma}^{\mathbf{e}}_+): (w,\boldsymbol{\sigma}) \in {\Lambda}_\AA^{n} \Big \},
\] together with \eqref{eq:VF_one_step_proof} and the fact that $\mathsf{V}_i$ are polynomials for $i \in \AA_E$, we have  \begin{align*}
        \sum_{\mathbf{a}\in \Lambda_{\AA_E}^{n}} \Vert \bV^\mathbf{a}(f)\Vert_c & = \sum_{e \in \mathbb{A}_{N+1}} \sum_{i \in \mathbb{A}, |i|=|e|} \sum_{(w,\bsigma) \in \Lambda_{\AA_E}^{n-1}}\Vert \bV^{(iw,\bsigma_+^e)}(f)\Vert_{c}
        \\ & \leq C \sum_{(w,\bsigma) \in \Lambda_{\AA_E}^{n-1}}\Big \{\sum_{e \in [N+1]}\sum_{i\in \cB} \sum_{l_1\in [k]}\Vert \partial_{x_e^{l_1}}\bV^{(w,\bsigma)}(f)\Vert_{c} \Vert \mathsf{V}_i^{l_1} \Vert_c \\ & \hspace{3cm}+ \sum_{e,h \in [N+1]^2}\sum_{[i,j]\in [\cM,\cM]} \sum_{l_1\in [k]^2}\Vert \partial_{x_e^{l^1_1}x_h^{l^2_1}}\bV^{(w,\bsigma)}(f)\Vert_{c} \Vert \mathsf{V}_{[i,j]}^{l_1} \Vert_c \Big \} \\& \leq C\max_{i \in \AA_E, l\in [k]^{|i|}} \Vert \mathsf{V}^l_i\Vert_c \sum_{(w,\bsigma) \in \Lambda_{\AA_E}^{n-1}}\Big \{\sum_{e \in [N+1]}\sum_{l_1\in [k]}\Vert \partial_{x_e^{l_1}}\bV^{(w,\bsigma)}(f)\Vert_{c}  \\ & \hspace{6cm}+ \sum_{e,h \in [N+1]^2} \sum_{l_1\in [k]^2}\Vert \partial_{x_e^{l^1_1}x_h^{l^2_1}}\bV^{(w,\bsigma)}(f)\Vert_{c} \Big \},
    \end{align*} where the constant $C$ changes from the second to the last inequality. To conclude the proof, we use the following elementary coefficient norm estimates for polynomials of order $p\in \mathbb{N}$ \[
    \sum_\nu \Vert \partial_{x_\nu}P\Vert_c\leq p \Vert P \Vert_c, \qquad \sum_{\nu,\mu}\Vert \partial_{x_\nu x_\mu} P\Vert_c\leq p(p-1)\Vert P\Vert_c.
    \]
    Combining the last two estimates, we have \begin{equation}\label{eq:one_step_poly}
        \sum_{\mathbf{a}\in \Lambda_{\AA_E}^{n}} \Vert \bV^\mathbf{a}(f)\Vert_c \leq A \sum_{\mathbf{a} \in \Lambda_{\AA_E}^{n-1}} \Vert \bV^\mathbf{a}(f) \Vert_c, \qquad A=Cp(p-1)\max_{i \in \AA_E, l\in [k]^{|i|}} \Vert \mathsf{V}^l_i\Vert_c.
    \end{equation}
    Iterating the last inequality, we have \begin{equation*}
        \sum_{\mathbf{a}\in \Lambda_{\AA_E}^{n}} \Vert \bV^\mathbf{a}(f)\Vert_c \leq A^n \Vert f \Vert_{c},
    \end{equation*}
which shows Assumption~\ref{ass:absolut_convergence_assumption} (i) for $\rho=0$ and $C=\Vert f\Vert_c$. Hence the proof is finished.
\end{proof}

\section{Application: Weak kernel approximation schemes}\label{sec:application}
In this section we present a first application of the series expansions developed in this work. 
In many Volterra SDE models of interest, the underlying kernel is singular, as is the case for the fractional kernel considered in Section~\ref{sec:fPVP}. 
This makes even elementary numerical tasks, such as the simulation of the corresponding process, challenging. 
A common strategy to address this issue is to replace the original kernel $K$ by a family of simpler approximating kernels $(K^\varepsilon)_{\varepsilon>0}$, for which simulation becomes a standard numerical task. 

In this section we provide a general methodology, to choose from a parametric class $\{K^\theta:\theta \in \Theta \}$, a kernel $K^{\theta^\star}$ in a way to minimize the \emph{weak-error} of the underlying Volterra equations. More precisely, let us denote by $Y^K$ the solution to a Volterra SDE \eqref{eq:volterra_eq_4}, with respect to the underlying kernel $K$. Given a parametric class of kernels, we define the weak error \begin{equation}\label{eq:weak_error}
    \mathcal{E}_T(\theta;f):= \big | \mathbb{E}[f(Y_T^K)]-\mathbb{E}[f(Y_T^{K^\theta})] \big |, \qquad f \in \mathcal{F}
\end{equation} where $\mathcal{F}\subset C^\infty(\RR^k;\RR)$ is a suitable class of test functions. Under  Assumption~\ref{ass:absolut_convergence_assumption}, by Theorem~\ref{thm:expansion_expected} we may represent the weak error as \begin{equation}\label{eq:matching_signature_vague}
    \mathcal{E}_T(\theta;f) \leq \sum_{n\geq 1}\sum_{\mathbf{a}\in \Lambda_{\AA_E}^n}\bV(f)(\xi) \mathcal{I}(\theta;\mathbf{a}), \qquad \mathcal{I}(\theta;\mathbf{a}):= |\bmu^\mathbf{a}_{0,T}(K)-\bmu^\mathbf{a}_{0,T}(K^\theta)|.
\end{equation}
The advantage of the series representation is that, the dependence structure $K\mapsto Y^K$ is precisely isolated and separated in particular from the dependence on the vector-fields, which coincide for both the original equation and its approximation. In the spirit of cubature methods for SDEs \cite{Lyons2004}, it is therefore natural to choose the kernel, that is the parameter $\theta \in \Theta$, such that the expected Volterra signatures are matched/minimized up to a certain degree $L\in \mathbb{N}$, that is  \begin{equation}\label{eq:Vsig_matching}
    \theta_L^\star \in \operatorname{argmin}_{\theta \in \Theta} \sum_{n=0}^L\sum_{\mathbf{a} \in \Lambda_{\AA_E}^n} \cI(\theta;\mathbf{a}). 
\end{equation}
\subsection{Numerical illustration for fractional Wright-Fisher models}\label{sec:numerics_jacobi}
We illustrate the proposed method for a special case of the fractional polynomial Volterra processes from Section~\ref{sec:fPVP}, namely one-dimensional \emph{fractional Wright-Fisher} processes on $[0,1]$:
\begin{equation}
    \label{eq:wright-fisher}
    Y_t = \xi + \lambda \int_0^t \frac{(t-s)^{\alpha-1}}{\Gamma(\alpha)}(b-Y_s)\dd s
    + \nu\int_0^t \frac{(t-s)^{\alpha-1}}{\Gamma(\alpha)}
    \sqrt{Y_s(1-Y_s)} \dd B_s, \quad 0\leq t \leq T,
\end{equation}
where $\xi, b\in [0,1]$ and $\lambda,\nu > 0$. 
\begin{remark}
Classical Wright-Fisher diffusions ($\alpha=1$) are used to describe allele-frequency fluctuations and arise as
large-population limits of discrete-generation reproduction
models originating in the work of Fisher and Wright
\cite{fisher1930genetical,wright1931evolution}.
Seed-bank extensions incorporate dormancy, allowing inactive individuals to
rejoin the reproducing population after many generations
\cite{blath2016seedbank}.
Such models can have scaling limits admitting non-Markovian Volterra
representations, providing a biological motivation for such Volterra
Wright-Fisher models; see also \cite{jiao2023seedbank}. For a general well-posedness result, we refer to \cite[Theorem 2.7]{abi2024polynomial}. 
\end{remark}

We choose sum-of-exponential kernels as the parametric class, that is
\begin{equation}\label{eq:SOE_classs}
\mathscr{K}_N:=\Big \{K^\theta(t,s)=\sum_{i=1}^N\beta_ie^{-\lambda_i(t-s)}:
\theta=\big ((\beta_1,\lambda_1),\dots,(\beta_N,\lambda_N)\big ) 
\in (\mathbb{R}\times [0,\infty))^N \Big \}.
\end{equation}
This choice is natural for two reasons: First, Volterra SDEs with kernel $K\in \mathscr{K}_N$ admit a \emph{Markovian lift} (see, e.g., \cite[Proposition 2.1]{alfonsi2024approximation}), that is, $Y^K$ can be represented via an $N$-dimensional diffusion, and one can exploit standard analytical and numerical methods for diffusions. On the other hand, for fractional kernels (and more generally for completely monotone kernels), $\mathscr{K}_N$ provides a natural approximation class since
\[
\frac{t^{\alpha-1}}{\Gamma(\alpha)}
= c_\alpha\int_0^\infty e^{-t\lambda} \lambda^{-\alpha}\dd \lambda
\approx \sum_{i=1}^N \beta_ie^{-\lambda_it},
\]
where the right-hand side can, for instance, be obtained by quadrature.

The goal of this numerical illustration is to show that our method produces a small weak error $\mathcal{E}_T(\theta^\star;f)$ in \eqref{eq:weak_error}, restricted to the monomials $\mathcal{F}= \{x \mapsto x^p:p\geq 1\}$. The reference moments of the polynomial Volterra process \eqref{eq:wright-fisher} are obtained by numerically\footnote{Details of the implementation can be found in \url{https://github.com/lucapelizzari/Moments-for-polynomial-Volterra-processes}.} solving the characterizing moment equation \cite[Theorem 3.3]{abi2024polynomial}. For a sum-of-exponential kernel, the moments are computed using the classical moment formula for polynomial diffusions, i.e., via a matrix exponential \cite[Theorem 2.7]{cuchiero2012polynomial}. 
\begin{figure}
    \centering
    \includegraphics[width=\linewidth]{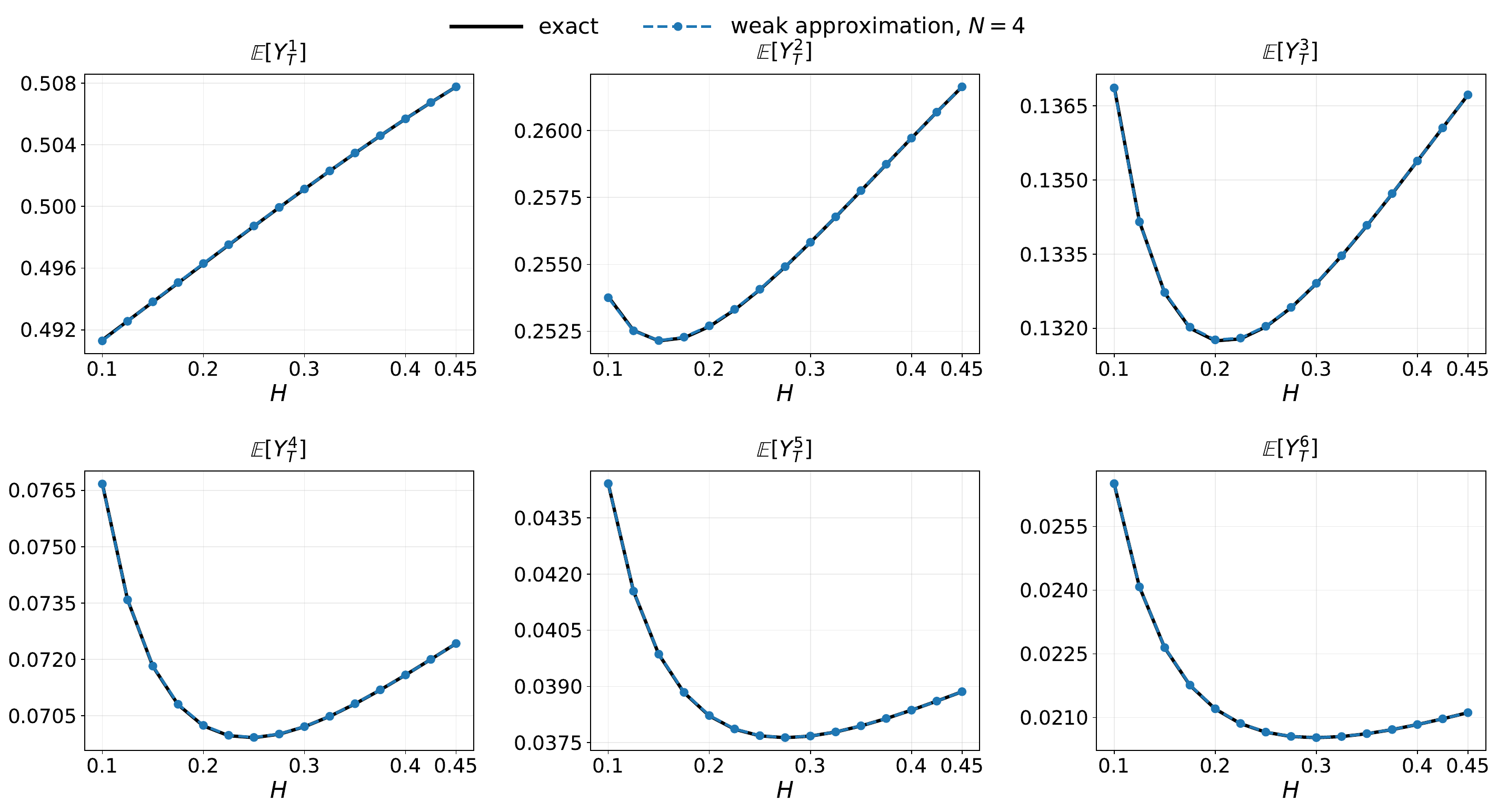}
    \caption{ 
The moments $\mathbb{E}[Y_T^p]$ for $p=1,\ldots,6$, as functions of $H=\alpha-1/2 \in (0,0.5)$. 
Solid lines denote the reference solution, obtained by numerically solving the moment equation from \cite{abi2024polynomial}, while dashed lines show the weak approximation obtained from our signature matching scheme for $\mathscr{K}_4$ in \eqref{eq:SOE_classs}.}
    \label{fig:raw_moments_approxi}
\end{figure}
Throughout this section, we fix the model parameters $T=1$, $\xi=0.3$, $b=0.6$, $\lambda=1.2$ and $\nu=0.2$, and vary the Hurst parameter $H=\alpha-1/2 \in (0,1/2)$. We select $\theta^\star$ by matching \eqref{eq:matching_signature_vague} up to level $L=2$, and all the kernel integrals appearing can be computed in closed form. More precisely, by definition of $\Lambda^2_{\AA_E}$, one easily checks that the following $12$ integrals are involved:
\[
\begin{cases}
    \displaystyle \int_0^T K(T,u)^a\dd u,
    & a \in \{1,2 \}, \\[0.8em]
    \displaystyle \int_{\Delta_{0,T}^2}
    K(T,s)^aK(s,r)^bK(T,r)^c\dd r \dd s,
    & a \in \{1,2\},
    (b,c) \in \{(1,0),(2,0),(0,1),(0,2),(1,1)\}.
\end{cases}
\]
For both the fractional kernel $K^H$ and the parametric kernel $K\in \mathscr K_N$, these integrals can be computed in closed form. Given a kernel $K$, we denote the above integrals by $I(K)=(I_j(K))_{j=1}^{S_L}$, and we select $\theta \in \Theta$ by numerically minimizing the loss
\begin{equation}\label{eq:kernel_matching_loss}
    \mathcal{L}_{S}(\theta)
    :=
    \sum_{j\in S}
    \left[
        \left(
            I_j(K^\theta)-I_j(K^H)
        \right)^2
        +
        \left(
            \frac{
                I_j(K^\theta)-I_j(K^H)
            }{
                I_j(K^H)+\varepsilon
            }
        \right)^2
    \right],
\end{equation}
where $\varepsilon=10^{-8}$ and $S$ denotes the set of integral coordinates matched at the chosen level $L$.

\begin{figure}
    \centering
    \includegraphics[width=\linewidth]{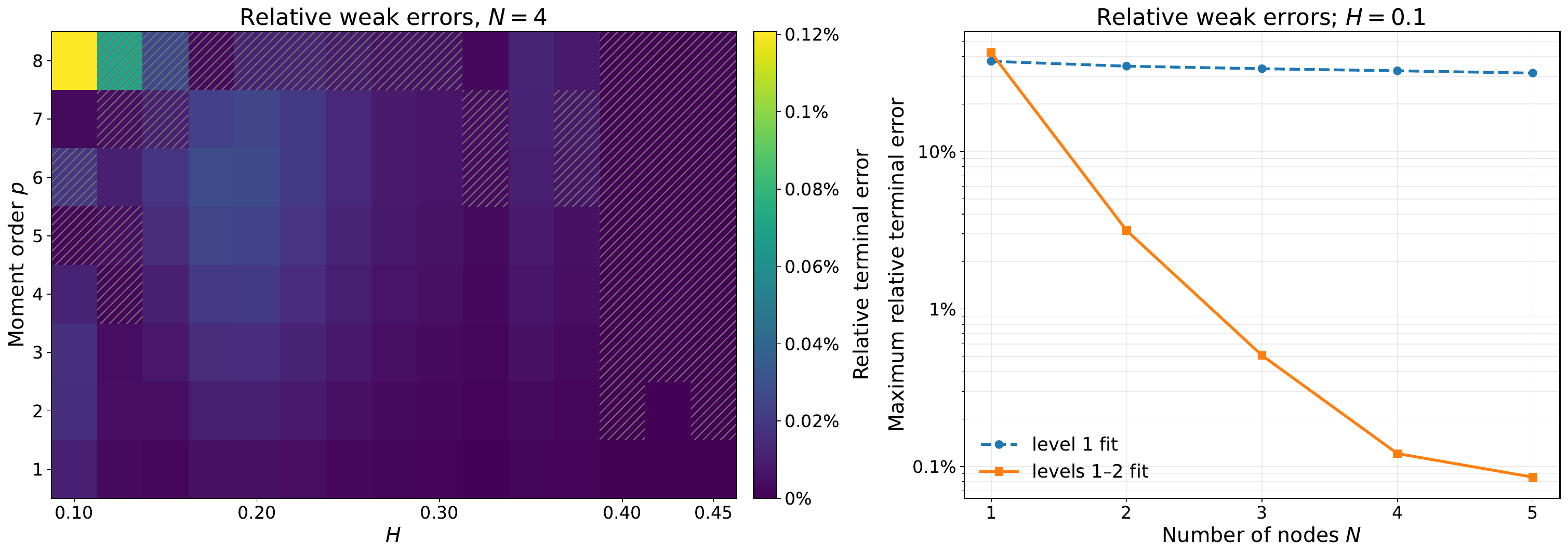}
    \caption{Left: relative errors of our method, matching the signature up to $L=2$ with $N=4$, as a function of the roughness parameter $H=\alpha-1/2$ and the moment order. Right: maximal relative error over the first eight moments at $H=0.1$ as a function of the number $N$ of exponential factors, comparing level-$L=1$ and level-$L=2$ matching.}
    \label{fig:heatmap}
\end{figure}

Figure~\ref{fig:raw_moments_approxi} illustrates the approximation directly at the level of the moments, for fixed $L=2$ and $N=4$. To the naked eye, the approximation with four exponentials is indistinguishable from the reference moments over all six moments and all considered Hurst parameters $H\in (0,0.5)$. A more systematic error quantification is presented in Figure~\ref{fig:heatmap}, where we plot the relative weak errors for moments up to order eight. The left panel displays the relative errors for the individual moments over the considered values of $H$, whereas the right panel summarizes the maximal relative error over the first eight moments as the number $N$ of exponential factors increases. We compare matching only the first-level integrals with matching both first- and second-level integrals.

\appendix

\section{Technical proofs Section~\ref{sec:prelim}}\label{sec:proofs}

\begin{proof}[Proof of Proposition~\ref{prop:vector_field_algebra}]
First, in the sum on the right-hand side of~\eqref{eq:recursion_formula}, we recall by definition
\begin{align*}
 & \bV^{(i_2\cdots i_n,\bsigma)}(f)(x_2,\dots,x_{n+1})
 \\
 & =
 \sum_{l_2 \in [k]^{|i_2|},\dots, l_n \in [k]^{|i_n|}}
 \prod_{q=2}^{n}
 \bigg (
 \prod_{\substack{(r,h): \bsigma_h(r)=q-1\\ 1\leq h\leq |i_{r+1}|}}
 \partial_{l^h_{r+1}}
 \bigg)
 \mathsf{V}_{i_{q}}^{l_{q}}(x_{q})
 \bigg(
 \prod_{\substack{(r,h): \bsigma_h(r)=n\\ 1\leq h\leq |i_{r+1}|}}
 \partial_{l^h_{r+1}}
 \bigg)
 f(x_{n+1}).
\end{align*}
Write $i_1:=\ell$. Suppose now $|\ell|=1$. By linearity, taking a partial derivative $\partial_{x_{e_1+1}^{l_1^1}}$ leads to the same expression but with one more derivative at position $q=e_1+1$. Multiplying the derivative of the summand with $\mathsf{V}_{i_1}^{l_1}(x_1)$ reads
\begin{align*}
    &\mathsf{V}_{i_1}^{l_1}(x_1)
    \partial_{x_{e_1+1}^{l_1^1}}
    \Big(
    \prod_{q=2}^{n}
    \Big(
    \prod_{\substack{(r,h): \bsigma_h(r)=q-1\\ 1\leq h\leq |i_{r+1}|}}
    \partial_{l^h_{r+1}}
    \Big)
    \mathsf{V}_{i_{q}}^{l_{q}}(x_{q})
    \Big(
    \prod_{\substack{(r,h): \bsigma_h(r)=n\\ 1\leq h\leq |i_{r+1}|}}
    \partial_{l^h_{r+1}}
    \Big)
    f(x_{n+1})
    \Big )
    \\
    & =
    \prod_{q=1}^n
    \Big(
    \prod_{\substack{(r,h): (\bsigma_+^e)_h(r)=q\\ 1\leq h\leq |i_r|}}
    \partial_{l^h_{r}}
    \Big)
    \mathsf{V}_{i_{q}}^{l_{q}}(x_{q})
    \Big(
    \prod_{\substack{(r,h): (\bsigma_+^e)_h(r)=n+1\\ 1\leq h\leq |i_r|}}
    \partial_{l^h_r}
    \Big)
    f(x_{n+1}),
\end{align*}
since $\bsigma_+^e$ is exactly $\bsigma$ shifted by one, with one more derivative added at position $e_1+1$. Exactly the same argument works when $|\ell|=2$ and we take two derivatives at $x_{e_1+1}$ and $x_{e_2+1}$, and $\bsigma_+^e$ adds two positions at $e_1+1$ and $e_2+1$. In summary, after summing the last identity over $l_1\in [k]^{|\ell|},\dots,l_n \in [k]^{|i_n|}$, we recover $\mathbf{V}^{(\ell i_2\cdots i_n,\bsigma_+^e)}(f)$ and the first identity \eqref{eq:recursion_formula} holds.
    
    We are left with the identity \eqref{eq:relation_VF_classical}, which we prove by induction. For $n=1$, there is a unique $\bsigma$ and
    \[
    \mathbf{V}^{(i,\bsigma)}(f)(\xi)
    =
    \sum_{l\in [k]^{|i|}}
    \mathsf{V}_{i}^l(\xi)\partial_l f(\xi)
    =
    \mathbf{V}^i(f)(\xi)
    \]
    directly from the definitions \eqref{eq:classical_VF_level_1} and \eqref{eq:v_field_operator_generalized}. By induction, assume the claim holds for $n-1$. By linearity of $\bV^{i_1}$ we have
    \begin{align*}
    \bV^{i_1} \circ \cdots \circ \bV^{i_n}(f)(\xi)
    &=
    \bV^{i_1}
    \big(
    \bV^{i_2} \circ \cdots \circ \bV^{i_n}(f)
    \big)(\xi)
    \\
    &=
     \sum_{\bsigma: (i_2\cdots i_n,\bsigma) \in \Lambda_\AA^{n-1}}
     \mathbf{V}^{i_1}
     \big(
     \mathbf{V}^{(i_2\cdots i_n,\bsigma)}(f)
     \big)(\xi).
    \end{align*}
    Now by definition and the product-rule, we have
    \begin{align}
    & \mathbf{V}^{i_1}
    \big(
    \mathbf{V}^{(i_2\cdots i_n,\bsigma)}(f)
    \big)
    \nonumber
    \\
    & =
    \sum_{l_1\in [k]^{|i_1|}}
    \sum_{l_2 \in [k]^{|i_2|},\dots,l_n \in [k]^{|i_n|}}
    \mathsf{V}_{i_1}^{l_1}
    \partial_{l_1}
    \left(
    \prod_{q=2}^{n}
    \left(
    \prod_{\substack{(r,h): \bsigma_h(r)=q-1\\ 1\leq h\leq |i_{r+1}|}}
    \partial_{l^h_{r+1}}
    \right)
    \mathsf{V}_{i_q}^{l_q}
    \left(
    \prod_{\substack{(r,h):\bsigma_h(r)=n\\ 1\leq h\leq |i_{r+1}|}}
    \partial_{l^h_{r+1}}
    \right)
    f
    \right)
    \nonumber
    \\
    & =
    \sum_{e\in [n]^{|i_1|}}
    \bV^{(i_1\cdots i_n,\bsigma_+^e)}(f),
    \label{eq:proof_vectorfield_id}
    \end{align}
    where the last equation follows from the recursive formula \eqref{eq:recursion_formula}. The last two observations combined show
    \begin{equation}
         \bV^{i_1} \circ \cdots \circ \bV^{i_n}(f)
         =
         \sum_{\bsigma: (i_2\cdots i_n,\bsigma) \in \Lambda_\AA^{n-1}}
         \sum_{e\in [n]^{|i_1|}}
         \mathbf{V}^{(i_1\cdots i_n,\bsigma_+^e)}(f).
    \end{equation}
    The claim follows if we can show
    \begin{equation}\label{eq:sets_recursion}
         \Big\{
         \bsigma:(i_1\cdots i_n,\bsigma)\in \Lambda_\AA^n
         \Big\}
         =
         \bigcup_{e\in [n]^{|i_1|}}
         \Big\{
         \btau_+^e:(i_2\cdots i_n,\btau)\in \Lambda_\AA^{n-1}
         \Big\}.
    \end{equation}
    The inclusion $\supseteq$ is clear, and for any $\bsigma$ from the left-hand side, set $e=\bsigma(1)-(1,1)\in \{1,\dots,n\}^2$ if $|i_1|=2$, and else $e=\bsigma_1(1)-1\in \{1,\dots,n\}$. Then choose $\btau$ such that $(i_2\cdots i_n,\btau)\in \Lambda_\AA^{n-1}$ and
    \[
        \btau(k)=\bsigma(k+1)-(1,1),
        \qquad k=1,\dots,n-1.
    \]
    Then by definition $\bsigma=\btau_+^e$, which finishes the proof.
\end{proof}

\medskip
\begin{proof}[Proof of Proposition~\ref{prop:shuffle_identity}]
        Suppose $\mathbf{a}=(i_1\cdots i_n, \sigma) \in \Lambda^n_\cA$ and $\mathbf{b}=(j_1\cdots j_p,\tau)\in \Lambda_\cA^p$, then by definition \begin{align*}
            \bZ_{s,t}^\mathbf{a}\bZ_{s,t}^\mathbf{b} & :=  \left (\int_{\Delta_{s,t}^n}
\dd Z_{r_1}^{r_{\sigma(1)},i_1}\cdots \dd Z_{r_n}^{r_{\sigma(n)},i_n}\right )\left (\int_{\Delta_{s,t}^p}
\dd Z_{r_1}^{r_{\tau(1)},j_1}\cdots \dd Z_{r_p}^{r_{\tau(p)},j_p}\right ) \\ &= \int_{\Delta_{s,t}^{n}\times  \Delta_{s,t}^{p}}\dd Z_{r_1}^{r_{\sigma(1)},i_1}\cdots \dd Z_{r_n}^{r_{\sigma(n)},i_n} \dd Z_{u_{n+1}}^{u_{n+\tau(1)},j_1}\cdots \dd Z_{u_{n+p}}^{u_{n+\tau(p)},j_p},
        \end{align*} with the conventions $r_{n+1}:=u_{n+p+1}:=t$.  Next, we can decompose the domain as follows \begin{align*}
            \Delta_{s,t}^{n} \times \Delta_{s,t}^{p} &= \bigcup_{\eta \in \mathrm{Sh}(n,p)} \{(r_1,\dots,r_{n+p}) \in [s,t]^{p+n}: r_{\eta^{-1}(1)}\leq \cdots \leq r_{\eta^{-1}(n+p)}\} \\ & =: \bigcup_{\eta \in \mathrm{Sh}(n,p)} \Delta_{s,t}^{(n+p),\eta}
        \end{align*} where the set $\mathrm{Sh}(n,p)$ was defined in Definition~\ref{def:classical_shuffle}.  Then it follows that \begin{align*}
            \bZ_{s,t}^\mathbf{a}\bZ_{s,t}^\mathbf{b} & = \sum_{\eta \in \mathrm{Sh}(n,p)} \int_{\Delta_{s,t}^{n+p,\eta}} \dd Z_{r_1}^{r_{\sigma(1)},i_1}\cdots \dd Z_{r_n}^{r_{\sigma(n)},i_n} \dd Z_{u_{n+1}}^{u_{n+\tau(1)},j_1}\cdots \dd Z_{u_{n+p}}^{u_{n+\tau(p)},j_p} \\ & =
            \sum_{\eta \in \mathrm{Sh}(n,p)} \int_{\Delta_{s,t}^{n+p}} \dd Z_{r_{\eta(1)}}^{r_{\tilde{\eta}(1)},i_1}\cdots \dd Z_{r_{\eta(n)}}^{r_{\tilde{\eta}(n)},i_n} \dd Z_{r_{\eta(n+1)}}^{r_{\tilde{\eta}(n+1)},j_1}\cdots \dd Z_{r_{\eta(n+p)}}^{r_{\tilde{\eta}(n+p)},j_p},
        \end{align*} where the second equality follows from a change of variables and using the notation \[
       \tilde{\eta}(k)= \begin{cases}
             \eta(\sigma(k)) & 1\leq k \leq n \text{ and } \sigma(k) \leq n
             \\ 
            \eta(n+\tau(k-n)) & n+1 \leq k \leq n+p \text{ and } \tau(k-n)\leq p. \\
            
            n+p+1 & \text{otherwise}.
            
        \end{cases}
        \] Now since $\eta$ is a permutation of $\{1,\dots,p+n\}$, we can order the integrators to $r_1,\dots,r_{n+p}$. For the superscript  we therefore transform $\tilde{\eta} \circ \eta^{-1}$, and one can observe (see Definition~\ref{def:Sigma_shuffle}) $$\tilde{\eta} \circ \eta^{-1} = \sigma \widehat{\shuffle}_\eta \tau.$$ Moreover, the letters in the superscript are transformed with $\eta^{-1}$ as well, that is \[
        \qquad h_k = \begin{cases}
            i_{\eta^{-1}(k)} & 1\leq \eta^{-1}(k)\leq n, \\ j_{\eta^{-1}(k)} & n+1 \leq \eta^{-1}(k) \leq n+m,
        \end{cases}
        \] which is simply the $k$ letter of the word $i_1\cdots i_n \shuffle_\eta j_1\cdots j_p$, see Definition~\ref{def:classical_shuffle}. Combining these observations, we have \begin{align}
            \bZ_{s,t}^\mathbf{a}\bZ_{s,t}^\mathbf{b} = \sum_{\eta \in \mathrm{Sh}(n,m)} \int_{\Delta_{s,t}^{n+m}} \prod_{l=1}^{n+m} \dd Z_{r_{l}}^{r_{(\sigma \widehat{\shuffle}_\eta \tau)(l)},h_l} & = \sum_{\eta \in \mathrm{Sh}(n,m)} \bZ_{s,t}^{w \shuffle_\eta v, \sigma \widehat{\shuffle}_\eta\tau}\nonumber  =\mathbf{Z}_{s,t}^{(w,\sigma)\shuffle (v,\tau)},
        \end{align} which finishes the proof.
    \end{proof}
    \medskip

\begin{proof}[Proof of Theorem~\ref{lem:finite_moments}]
It suffices to consider $p\geq2$, since the remaining
cases follow from the inclusion $L^2\subseteq L^p$.

\emph{Step 1:} First we extend the definition of the Itô Volterra signature to obtain a recursive definition. For fixed $(s,t)\in\Delta_{0,T}^2$, $n\in\mathbb N$,
 and a given word $w=i_1\cdots i_n$,
consider maps
\[
    \breve{\bsigma}:
    \{1,\ldots,n\}\longrightarrow\{2,\ldots,3n\}^2
\]
such that $\breve{\bsigma}(k)>k$ component-wise and
$\breve\sigma_1(k)=\breve\sigma_2(k)$ whenever
$|i_k|=1$.
For $\tau=(\tau_1,\dots,\tau_{2n}) \in [0,T]^{2n}$ we can define 
$\bZ_{s,t}^{w,\breve{\bsigma},\tau}$ as in
Definition~\ref{def:Ito_Volterra signature}, with the additional conventions $r_{n+j}:=\tau_j$ for $j=1,\ldots,2n.
$ Note that this generalizes Definition~\ref{def:Ito_Volterra signature}, when choosing for instance $\tau_1=t$ and restrict to $\breve \bsigma \in (\Sigma^n)^2$. In addition, this generalized object now has a recursive nature: 
\begin{equation}\label{eq:exp_vsig_lemma}
    \bZ_{s,t}^{w,\breve \bsigma,\tau}
    =
    \int_s^t
    \bZ_{s,u}^{i_1\cdots i_{n-1},
                 \breve{\bsigma}',\tau'(u,\tau)}
    \dd \cZ_u^{\tau_{\bsigma},i_n}, \qquad \tau_{\bsigma}:=\tau_{\bsigma(n)-n},
\end{equation}
where the inner integral can again be represented using at most 
$2(n-1)$ external parameters. Indeed, among $u,\tau_1,\ldots,\tau_{2n}$, we keep only the upper
parameters referenced by the first $n-1$ letters, and easily notice that at most $2(n-1)$ auxiliary parameters are required. Listing them in their inherited index order and replacing unused entries with $T$, if necessary, determines a new map $\breve{\bsigma}'$ and a parameters $\tau'(u,\tau)\in [0,T]^{2(n-1)}$, such that~\eqref{eq:exp_vsig_lemma} holds. 

\emph{Step 2:} In order to estimate the involved  stochastic integrals, let us recall some classical semimartingale estimates from \cite{Protter2004} and a general Kolmogorov-Centov result. First for any $v\geq 2$, recall the $\mathcal{S}^v$-norm $\|H\|_{\mathcal S^v}
    :=
    \|\sup_{t\leq T}|H_t|\|_{L^v},$ see \cite[Chapter V Section 2]{Protter2004}. Then, for any continuous semimartingale $X$ starting at zero and
a predictable process $H$, using the notation $H\cdot X=\int H\dd X$, we have the following inequalities 
\begin{equation}\label{eq:semimartingale_product_bounds}
    \|H\cdot X\|_{\mathcal H^v}
    \leq
    \|H\|_{\mathcal S^{2v}}\,
    \|X\|_{\mathcal H^{2v}},\qquad 
    \|[X,Y]\|_{\mathcal H^v}
    \leq
    \|X\|_{\mathcal H^{2v}}\,
    \|Y\|_{\mathcal H^{2v}},\qquad 
    \|X\|_{\mathcal S^v}
    \leq C_v\|X\|_{\mathcal H^v}.
\end{equation}
The first one is an application of the \'Emery inequality \cite[Chaper V Theorem 3]{Protter2004}, the second one is an immediate consequence of the Kunita-Watanabe inequality \cite[Chaper II Theorem 25]{Protter2004}, and the last one is precisely the maximal inequality \cite[Chapter V Theorem 2]{Protter2004}. For the Kolmogorov-Centsov continuity theorem, consider a multi-parameter stochastic process $[0,T]^d \ni \theta \mapsto F^\theta=(F_t^\theta)_{t\in [0,T]}$, which is continuous in $t$. Then, if $F$ is such that for $P\geq 1$ and $\alpha > 0$ we have  \begin{equation}\label{eq:kolmo_suff} \sup_{\tau} \Vert F^\tau \Vert_{S^P}+ \sup_{\tau \neq \eta}\frac{\Vert F^\tau-F^\eta \Vert_{\cS^P}}{|\tau-\eta |_{\infty}^{\alpha }} \leq C_0,\end{equation} then for any $0< \gamma < \alpha-d/P$ there exists a modification $F$, jointly continuous in time and the upper-parameters, such that 

\begin{equation}\label{eq:parameter_kolmogorov}
\left\|  \sup_{\theta,t}|F_t^\theta|
       +
       \sup_{\theta\neq\eta}
       \frac{\sup_{t\leq T}|F_t^\theta-F_t^\eta|}
            {|\theta-\eta|_\infty^\gamma}
     \right\|_{L^P}
    \leq
    CC_0.
\end{equation}

\emph{Step 3:} Next, we prove by induction over $n\geq 1$, that for every $p\geq 2$, we can find $q=q_n(p)$, such that \begin{equation}\label{eq:enhanced_vsig_moment_bound}
    \sup_{s\in[0,T]}
    \|\bZ_{s,\cdot}^{w,\breve \bsigma,\cdot }\|_{(\gamma_n,p)}
    \leq
    C\|Z\|_{\cV_{\mathbb F}^{\beta,q}}^{|w|}, \qquad w \in \cW^n_\AA, \quad \gamma_n=2^{1-n}\beta, \quad \breve \bsigma \text{ as in Step 1.}
\end{equation} For $n=1$, the process $t\mapsto \bZ_{s,t}^\tau$ is the increment of
$\cZ^{\tau_{\bsigma,i_1}}$, and if $|i_1|=1$ we use the definition of
$\cV_{\mathbb F}^{\beta,p}$, and if $i_1=[i,j]$ we can use bi-linearity of the quadratic variation and the second inequality in 
\eqref{eq:semimartingale_product_bounds}, to see 
\[
    \sup_{s\in[0,T]}
    \|\bZ_{s,\cdot}^{i_1,\breve \bsigma,\cdot }\|_{(\beta,p)}
    \leq
    C\|Z\|_{\cV_{\mathbb F}^{\beta,p}}, \qquad \sup_{s\in[0,T]}
    \|\bZ_{s,\cdot}^{[i,j],\breve \bsigma,\cdot }\|_{(\beta,p)}
    \leq
    C\|Z\|_{\cV_{\mathbb F}^{\beta,2p}}^2.
\]
Thus \eqref{eq:enhanced_vsig_moment_bound} holds for
$n=1$. Suppose the claim holds for $n-1$, and fix $p\geq2$. Recalling the recursion~\eqref{eq:exp_vsig_lemma}, we can write $\bZ_{s,t}^{w,\breve \bsigma,\tau} = (\bZ_{s,\cdot }^{i_1\cdots i_{n-1},
                 \breve{\bsigma}',\tau'(\cdot ,\tau)} \cdot \cZ^{\tau_{\bsigma},i_n})_t=: (H_s^\tau \cdot \cZ^{\tau_{\bsigma},i_n})_t$, and therefore, using the inequalities in~\eqref{eq:semimartingale_product_bounds} once again, we have \begin{align*}
                     \Vert \bZ_{s,\cdot }^{w,\breve \bsigma,\tau}-\bZ_{s,\cdot }^{w,\breve \bsigma,\eta} \Vert_{\cH^p} & \leq \Vert (H_s^\tau-H_s^\eta)\cdot  \cZ^{\tau_{\bsigma},i_n} \Vert_{\cH^p}  + \Vert  H_s^\eta\cdot(\cZ^{\eta_{\bsigma},i_n}-\cZ^{\tau_{\bsigma},i_n})\Vert_{\cH^p} \\ &
                 \leq C \Vert  H_s^\tau-H_s^\eta \Vert_{\cS^{2p}} \Vert   \cZ^{\tau_{\bsigma},i_n} \Vert_{\cH^{2p}} +  \Vert  H_s^\eta \Vert_{\cS^{2p}}\Vert\cZ^{\eta_{\bsigma},i_n}-\cZ^{\tau_{\bsigma}} \Vert_{\cH^{2p}}.
                 \end{align*}
                 Using the base case $n=1$ and the fact $\gamma_n<\gamma_1=\beta$, there exists a $q_1(p)$ such that $$\Vert\cZ^{\eta_{\bsigma},i_n}-\cZ^{\tau_{\bsigma}} \Vert_{\cH^{2p}} \leq  C \Vert Z \Vert_{(\beta,q_1(2p))}^{|i_n|}|\tau-\eta|_\infty^{\gamma_n}, \qquad \Vert   \cZ^{\tau_{\bsigma},i_n} \Vert_{\cH^{2p}} \leq C \Vert Z \Vert_{(\beta,q_1(2p))}^{{|i_n|}}.$$
                 We can the show the claim if we find $Q_{n-1}$ larger enough, such that \begin{equation}
                     \label{eq:proof_claim}
                     \Vert  H_s^\tau-H_s^\eta \Vert_{\cS^{2p}}  \leq C \Vert Z\Vert_{(\beta,Q_{n-1})}^{|i_1\cdots i_{n-1}|}|\tau-\eta|_\infty^{\gamma_n}, \qquad \Vert  H_s^\eta \Vert_{\cS^{2p}} \leq C Z\Vert_{(\beta,Q_{n-1})}^{|i_1\cdots i_{n-1}|}.
                 \end{equation}
                 Indeed, in this case we choose $q_n(p)= \max (q_1(2p),Q_{n-1})$, and the claim follows combinind the last three equations. The proof of~\eqref{eq:proof_claim} would be straightforward, if $t\mapsto H_u^\tau(t)$ was a semimartingale, since one could directly invoke again~\eqref{eq:semimartingale_product_bounds} and the induction hypothesis again. However, since $\tau'(t,\tau)$ depends on time, this is not the case in general, and we circumvent this using the Kolmogorov-Centsov argument from Step 2. Note that we have, for any $\gamma>0$ 
                    \begin{align*}
                     \Vert H_s^\tau-H_s^\eta \Vert_{\mathcal{S}^{2p}}& := \Vert \sup_{0\leq t \leq T} |H_s^\tau(t)-H_s^\eta(t)| \Vert_{L^{2p}} \\ & \leq \Vert  \sup_{\mu,\theta \in [0,T]^{2(n-1)}} \frac{\sup_{0\leq t \leq T}|\bZ_{s,t }^{i_1\cdots i_{n-1},
                 \breve{\bsigma}',\mu}-\bZ_{s,t}^{i_1\cdots i_{n-1},
                 \breve{\bsigma}',\theta}||\tau'(t,\tau)-\eta'(t,\eta)|^\gamma}{|\mu-\theta|_\infty^{\gamma}} \Vert_{L^{2p}} \\ & \leq |\tau-\eta|_\infty^{\gamma} \Big \Vert  \sup_{\mu,\theta \in [0,T]^{2(n-1)}, \mu \neq \theta} \frac{\sup_{0\leq t \leq T}|\bZ_{s,t}^{i_1\cdots i_{n-1},
                 \breve{\bsigma}',\mu}-\bZ_{s,t}^{i_1\cdots i_{n-1},
                 \breve{\bsigma}',\theta}|}{|\mu-\theta|_\infty^{\gamma}} \Big \Vert_{L^{2p}},
                 \end{align*} where we simply have used that $|\tau'(t,\tau)-\eta'(t,\tau)| \leq |\tau-\eta|_\infty$ by definition. The goal is now to apply Kolmogorov-Centsov~\eqref{eq:parameter_kolmogorov} for $\gamma_n$, that is we have to show~\eqref{eq:kolmo_suff} for $\hat{H}^\mu_s(t)=\bZ_{s,t}^{i_1\cdots i_{n-1},\breve \bsigma',\mu}$. But since we are now evaluating at generic upper-parameters, we can invoke the induction hypothesis, together with the maximal inequality in~\eqref{eq:semimartingale_product_bounds}, that is \[
                 \sup_{\mu \neq \theta}\frac{\Vert \hat{H}^\mu-\hat{H}^\theta \Vert_{\cS^P}}{|\mu-\theta |_{\infty}^{\alpha}} \leq C \sup_{\mu \neq \theta}\frac{\Vert \hat{H}^\mu-\hat{H}^\theta \Vert_{\cH^P}}{|\mu-\theta |_{\infty}^{\alpha}} \leq C\Vert Z \Vert_{(\beta,q_{n-1}(P))}^{|i_1\cdots i_{n-1}}, \qquad \forall \alpha <\gamma_{n-1}.
                 \] Note that we can do exactly the same for $\sup_{\mu} \Vert \hat{H}^\mu \Vert_{\cS^P}$. Now if $P\geq 2p$ is such that $\gamma_n<\gamma_{n-1}-\frac{2(n-1)}{P}$, which is equivalent to $2(n-1)<\gamma_nP$, we can deduce~\eqref{eq:proof_claim} from~\eqref{eq:parameter_kolmogorov}.
                  using $Q_{n-1}=q_{n-1}(P)$, so that the claim follows.

\emph{Step 4:} 
It remains to prove continuity. Let $\mathbf a\in\Lambda_\AA$ and set $r=|\mathbf a|$. By multilinearity of stochastic integration and
bilinearity of quadratic covariation, the preceding
induction also applies to mixed iterated integrals
$\mathcal I_{s,t}^{\mathbf a}(Z^1,\ldots,Z^r)$,
where each ordinary letter corresponds to one input
and each bracket letter to two inputs. In particular,
the same estimates yield, for sufficiently large $q$,
\[
    \sup_{(s,t)\in\Delta_{0,T}^2}
    \big\|\mathcal I_{s,t}^{\mathbf a}
    (Z^1,\ldots,Z^r)\big\|_{L^p}
    \leq C\prod_{j=1}^r
    \|Z^j\|_{(\beta,q)}.
\]
Since
$\bZ_{s,t}^{\mathbf a}
=\mathcal I_{s,t}^{\mathbf a}(Z,\ldots,Z)$,
multilinearity and a telescoping sum give
\[
\begin{aligned}
&\sup_{(s,t)\in\Delta_{0,T}^2}
\|\bZ_{s,t}^{\mathbf a}
-\bY_{s,t}^{\mathbf a}\|_{L^p}\\
&\quad\leq
C\|Z-Y\|_{(\beta,q)}
\sum_{j=1}^r
\|Y\|_{(\beta,q)}^{j-1}
\|Z\|_{(\beta,q)}^{r-j}\\
&\quad\leq
C\bigl(\|Z\|_{(\beta,q)}
+\|Y\|_{(\beta,q)}\bigr)^{r-1}
\|Z-Y\|_{(\beta,q)}.
\end{aligned}
\]
Thus, on bounded subsets of
$\cV_{\mathbb F}^{\beta,q}$, the signature
coordinate is Lipschitz continuous, which finishes the proof.
\end{proof}
\begin{proof}[Proof of Proposition~\ref{prop:quasi_shuffle}]
Fix a process $Z\in\cV^{\beta,\infty-}_{\mathbb F}$ such that we prove the identity for every fixed $0\leq s\leq t\leq T$. 
First, similar to the proof of Theorem~\ref{lem:finite_moments}, we  introduce an auxiliary family of iterated integrals which is closed under removing the last integration variable. For $n \geq 1$ define 
\begin{align*}
    \breve\Sigma^n
    &:=
    \left\{
        \sigma:[n]\to[3n]:
        \sigma(\ell)>\ell \text{ for every }\ell\in[n]
    \right\},\\
    \breve\Lambda_\AA^n
    &:=
    \left\{
        (i_1\cdots i_n,\bsigma):
        \bsigma\in(\breve\Sigma^n)^2,\;
        \bsigma_1(\ell)=\bsigma_2(\ell)
        \text{ whenever }|i_\ell|=1
    \right\}.
\end{align*}
The indices $1,\ldots,n$ represent integration variables, while $n+1,\ldots,3n$ represent $2n$ external parameters. For
$\mathbf a=(i_1\cdots i_n,\bsigma)\in\breve\Lambda_\AA^n$, $0\leq s\leq u\leq T$, and $\boldsymbol\alpha=(\alpha_1,\ldots,\alpha_{2n})\in[u,T]^{2n}$, define
\begin{align*}
    I_{s,u}^{\mathbf a}(\boldsymbol\alpha)
    :=
    \int_{\Delta_{s,u}^n}
    \dd\cZ_{r_1}^{r_{\bsigma(1)},i_1}
    \cdots
    \dd\cZ_{r_n}^{r_{\bsigma(n)},i_n},
    \qquad
    r_{n+j}:=\alpha_j,\quad j\in[2n].
\end{align*}
The ordinary and quadratic-covariation integrators are
understood as in~\eqref{eq:VSIG_general_ito} and note that we understand $I_{s,u}^{\varnothing}=1$. Note that the external parameters are not necessarily ordered. 
We make the recursive nature of this extended signature explicit. Define for $n \geq 2$ the index set
\begin{align*}
    D_{\mathbf a}
    =\{d_1<\cdots<d_q\} :=
    \left\{
        \bsigma_h(\ell):
        \ell\in[n-1],\ h\in\{1,2\}
    \right\}
    \cap\{n,\ldots,3n\},
\end{align*}
consisting of the $\bsigma$-targets which become external after removing $n$ from the domain of $\bsigma$. Note that this in particular may include the target index $n$ itself. Observe that $|D_{\boldsymbol{a}}| = q \leq 2(n-1)$ and define 
\begin{align*}
    \bsigma_h^-(\ell)
    &:=
    \begin{cases}
        \bsigma_h(\ell),
        & \bsigma_h(\ell)<n,\\
        n-1+j,
        & \bsigma_h(\ell)=d_j,
    \end{cases}
    \qquad \ell\in[n-1],\quad h\in\{1,2\},\\
    \alpha_j^-(r)
    &:=
    \begin{cases}
        r, & d_j=n,\\
        \alpha_{d_j-n}, & d_j>n,
    \end{cases}
    \qquad j\in[q].
\end{align*}
If $q<2(n-1)$, we fill the unused external slots with the terminal time $T$. Note that we have $\boldsymbol\alpha^-(r)\in[r,T]^{2(n-1)}$ for $r \leq u$ and $\mathbf a^-:=(i_1\cdots i_{n-1},\bsigma^-)\in\breve\Lambda_\AA^{n-1}$. For $n=1$ we use the convention $\boldsymbol{a}^- = \varnothing$ and $\boldsymbol{\alpha}^-(r) = \varnothing$. Recalling that $\sigma(n)>n$ by definition, let us define 
\begin{align*}
    M_r^{\mathbf a}
    :=
    \cZ_r^{
    (\alpha_{\bsigma_1(n)-n}, \alpha_{\bsigma_2(n)-n}),\,i_n},
\end{align*}
such that we have the recursive structure 
\begin{align*}
    I_{s,u}^{\mathbf a}(\boldsymbol\alpha)
    =
    \int_s^u
    I_{s,r}^{\mathbf a^-}
        \bigl(\boldsymbol\alpha^-(r)\bigr)
    \dd M_r^{\mathbf a}.
\end{align*}

We show that this auxiliary object $I$ satisfies the quasi-shuffle property. Consider two elements $\mathbf a=(w,\bsigma)\in\breve\Lambda_\AA^n$ and $\mathbf b=(v,\btau)\in\breve\Lambda_\AA^m$, with external times $\boldsymbol\alpha\in[u,T]^{2n}$ and $\boldsymbol\beta\in[u,T]^{2m}$, respectively. Fix $\eta\in\mathrm{QSh}(n,m)$ with image $[k]$. For the integration variables $r_1<\cdots<r_k$, we define for $h\in\{1,2\}$
\begin{align*}
    \Theta_h^\eta(p)
    :=
    \begin{cases}
        r_{\eta(\bsigma_h(p))},
        & p\leq n,\ \bsigma_h(p)\leq n,\\
        \alpha_{\bsigma_h(p)-n},
        & p\leq n,\ \bsigma_h(p)>n,\\
        r_{\eta(n+\btau_h(p-n))},
        & p>n,\ \btau_h(p-n)\leq m,\\
        \beta_{\btau_h(p-n)-m},
        & p>n,\ \btau_h(p-n)>m.
    \end{cases}, \quad \boldsymbol\theta_e^\eta
    :=
    \left(
        \Theta_1^\eta(\eta_-(e)),
        \Theta_2^\eta(\eta_+(e))
    \right),
\end{align*}
for $e\in[k]$. Recall that $\eta_-(e)=\min\eta^{-1}(e)$ and $\eta_+(e)=\max\eta^{-1}(e)$ and define
\begin{align*}
    J_{s,u}^\eta
        (\mathbf a,\boldsymbol\alpha;
         \mathbf b,\boldsymbol\beta)
    :=
    \int_{\Delta_{s,u}^k}
    \dd\cZ_{r_1}^{\boldsymbol\theta_1^\eta,\ell_1}
    \cdots
    \dd\cZ_{r_k}^{\boldsymbol\theta_k^\eta,\ell_k} \quad \text{if }  w\star_\eta v=\ell_1\cdots\ell_k\neq0,  \quad \text{and } J_{s,u}^\eta =0 \text{ else.}
\end{align*}
One can verify that every nonzero $J_{s,u}^\eta$ is again an integral of the type introduced above. Indeed, each upper-variable that is evaluated within the integral, is strictly later than the variable itself, since $\eta$ preserves the order within each input word. Similarly $\boldsymbol{\theta}^\eta$ is valid since if $|l_e| = 1$, we have $|\eta^{-1}(e)| = 1$ by definition of $\star_\eta$ and hence $\eta_+(e) = \eta_-(e)$. Likewise, in this case the concatenation $wv$ at entry $\eta_+(e)$ either has length one and hence $\btau_h$ or $\bsigma_h$ respectively have the same component at this entry, or the shuffle $w \star_\eta v = 0$ and hence we would have $J^\eta_{s,u} = 0$. We will continue in proving the identity
\begin{align} \label{eq:proof_quasi_shuffle}
    I_{s,u}^{\mathbf a}(\boldsymbol\alpha)
    I_{s,u}^{\mathbf b}(\boldsymbol\beta)
    =
    \sum_{\eta\in\mathrm{QSh}(n,m)}
    J_{s,u}^\eta
        (\mathbf a,\boldsymbol\alpha;
         \mathbf b,\boldsymbol\beta)
\end{align}
by induction on $N = n+m$. If either word is empty, the identity follows from the convention $I^\varnothing_{s,u} = 1$ and that only admissible map $\eta$ on the right hand side is the identity.

Suppose that the identity holds whenever the sum of the input word lengths is strictly less than $N$ and take $n,m \geq 1$ with $n+m = N$. For fixed external times we introduce the notation
\begin{align*}
    A_r
    :=
    I_{s,r}^{\mathbf a}(\boldsymbol\alpha),
    \quad 
    U_r
    :=
    I_{s,r}^{\mathbf a^-} \bigl(\boldsymbol\alpha^-(r)\bigr),\quad 
    B_r
    :=
    I_{s,r}^{\mathbf b}(\boldsymbol\beta),
    \quad 
    V_r
    :=
    I_{s,r}^{\mathbf b^-}\bigl(\boldsymbol\beta^-(r)\bigr).
\end{align*}
By the above, we have the recursions $A_r=\int_s^r U_v\,\dd M_v$ and $B_r=\int_s^r V_v\,\dd N_v$, where $M=M^{\mathbf a}$ and $N=M^{\mathbf b}$ have fixed upper external times. By Itô's formula, we have
\begin{align*}
    A_uB_u
    =
    \int_s^u A_rV_r\,\dd N_r
    +
    \int_s^u U_rB_r\,\dd M_r
    +
    \int_s^u U_rV_r\,\dd[M,N]_r.
\end{align*}
To identify these three terms, we partition
$\mathrm{QSh}(n,m)$ according to the last fiber of respective maps $\eta$ in this set, i.e.\ according to $\eta^{-1}(k)$ for $k \in \{n \lor m, \dots, n+m\}$. 
Let $\eta$ be any such map with image $[k]$. By strict monotonicity on the two input blocks, the last fiber has exactly one of the forms $\{n\}$, $\{n+m\}$, or $\{n,n+m\}$. Denote the corresponding classes by $\mathcal Q_A(n,m)$, $\mathcal Q_B(n,m)$, and $\mathcal Q_{AB}(n,m)$, respectively.
Deleting and restoring the respective last fiber, we can identify these classes with quasi-shuffles of lower order, i.e. we have the bijections
\begin{align*}
    \mathcal Q_A(n,m)
    &\longleftrightarrow
    \mathrm{QSh}(n-1,m),\\
    \mathcal Q_B(n,m)
    &\longleftrightarrow
    \mathrm{QSh}(n,m-1),\\
    \mathcal Q_{AB}(n,m)
    &\longleftrightarrow
    \mathrm{QSh}(n-1,m-1).
\end{align*}
Starting with $\int_s^u A_rV_r\,\dd N_r$, the induction hypothesis applies to $A_rV_r$ for each $r \in [s,u]$, since its input lengths are $n$ and $m-1$. Hence we have
\begin{align*}
    \int_s^u A_rV_r\,\dd N_r 
    &=
    \sum_{\eta \in \mathrm{QSh}(n,m-1)} \int_{s}^u J^\eta_{s,r}
        (\mathbf a,\boldsymbol\alpha;
         \mathbf b^-,\boldsymbol\beta^-(r)) \dd N_r,
    =
    \sum_{\eta\in\mathcal Q_B(n,m)}
    J_{s,u}^\eta
        (\mathbf a,\boldsymbol\alpha;
         \mathbf b,\boldsymbol\beta),
\end{align*}
where the second equality is derived as follows: take any $\eta \in \mathrm{QSh}(n,m-1)$ with co-domain $[k-1]$, and note that $N$ has upper index $\beta_{\btau_1(m) - m}, \beta_{\btau_2(m) - m}, j_m$. The first two entries are exactly $\boldsymbol{\theta}^{\breve\eta}_k = (\Theta_1^{\breve\eta}(\breve \eta_-(k)), \Theta_2^{\breve\eta}(\breve \eta_-(k)))$ where $\breve \eta$ is chosen according to the bijection established above with $\breve \eta^{-1}(k) = \{n+m\}$. The letter $j_m$ coincides with the last letter of the shuffle $w \star_{\breve \eta} v$ by Definition \ref{def:classical_quasi_shuffle}. Lastly, note that $\boldsymbol{b}^-$ and $\boldsymbol{\beta}^-(r)$ were only relabeled and the upper indices were not changed apart from removing the last one, which is now appended back again by $N$. This shows the equality. 
The exact same argument shows the equality
\begin{align*}
    \int_s^u U_rB_r\,\dd M_r
    =
    \sum_{\eta\in\mathcal Q_A(n,m)}
    J_{s,u}^\eta
        (\mathbf a,\boldsymbol\alpha;
         \mathbf b,\boldsymbol\beta).
\end{align*}
Lastly, we identify the quadratic variation term. If either last letter of $w$ and $v$ is a bracket letter, we have $[M,N] = 0$. In this case, $J_{s,u}^\eta = 0$ as well for each $\eta \in \cQ_{AB}(n,m)$ by Definition \ref{def:classical_quasi_shuffle}. Otherwise, let the last letters be $i_n,j_m\in\cA$ such that 
\begin{align*}
    M 
    =
    Z^{
    \alpha_{\bsigma_1(n)-n}\,i_n}, \qquad
    N
    =
    Z^{
    \beta_{\btau_1(m)-m}\,j_m}.
\end{align*}
Again the induction hypothesis implies 
\begin{align*}
    \int_s^u U_rV_r\,\dd[M,N]_r
    =
    \sum_{\eta\in\mathrm{QSh}(n-1,m-1)} \int_s^u
    J^\eta_{s,r}
        (\mathbf a^-,\boldsymbol\alpha^-(r);
         \mathbf b^-,\boldsymbol\beta^-(r)) \dd \ZZ^{(\alpha_{\bsigma_1(n)-n}, \beta_{\btau_1(m)-m}),[i_n,j_m]}_r.
\end{align*}
To identify the right hand side, we take any $\eta \in \mathrm{QSh}(n-1,m-1)$ with co-domain $[k-1]$ and accordingly define $\breve \eta \in \cQ_{A,B}(n,m)$ such that $\breve \eta^{-1}(k) = \{n,n+m\}$. Similarly to the above $\boldsymbol{\theta}^{\breve\eta}_k$ agrees precisely with $(\alpha_{\bsigma_1(n)-n}, \beta_{\btau_1(m)-m})$ and the last letter $[i_n,j_m]$ is the last letter of $w \star_{\breve \eta} v$ by Definition \ref{def:classical_quasi_shuffle}.
The inner upper indices were again only relabeled and not changed apart from removing their last entry, which is now appended back again by $[M,N]$. This identifies 
\begin{align*}
    \int_s^u U_rV_r\,\dd[M,N]_r
    =
    \sum_{\eta\in\mathcal{Q}_{A,B}(n,m)} 
    J^\eta_{s,u}
        (\mathbf a,\boldsymbol\alpha;
         \mathbf b,\boldsymbol\beta).
\end{align*}
Hence the identify \eqref{eq:proof_quasi_shuffle} holds. 

Finally, take $\mathbf a=(w,\bsigma)\in\Lambda_\AA^n$ and $\mathbf b=(v,\btau)\in\Lambda_\AA^m$ from the original index set. Regard them as auxiliary indices using only
the first external slot and set all external parameters
equal to the fixed time $t$ and take $u=t$. Then we have $I_{s,t}^{\mathbf a}(t,\ldots,t)=\bZ_{s,t}^{\mathbf a}$ and similarly for $\mathbf b$. To identify the right-hand side of \eqref{eq:proof_quasi_shuffle} fix some surjection $\eta:[n+m] \twoheadrightarrow [k]$, extend it by
$\eta(n+m+1)=k+1$, and set $r_{k+1}:=t$. Using the original concatenation of evaluation maps with a common terminal root, the target rule above becomes $\Theta_h^\eta(p) =r_{\eta((\sigma_h\tau_h)(p))}$. Consequently,
\begin{align*}
    \boldsymbol\theta_e^\eta
    =
    (
        r_{(\eta\circ(\sigma_1\tau_1)\circ\eta_-)(e)},
        r_{(\eta\circ(\sigma_2\tau_2)\circ\eta_+)(e)}
    )
    =
    (
        r_{(\bsigma\widehat\star_\eta\btau)_1(e)},
        r_{(\bsigma\widehat\star_\eta\btau)_2(e)}
    ).
\end{align*}
Thus every nonzero auxiliary term on this diagonal is
exactly the coordinate prescribed by our quasi-shuffle
product on $\Lambda_\AA$, which finishes the proof.
\end{proof}
\section*{Acknowledgements}
This research was funded in whole or in part by the Austrian Science Fund (FWF) Y 1235.
\bibliographystyle{plain}
\bibliography{bib}
\end{document}